\documentclass[12pt,reqno]{amsart}

\usepackage[final
]{showkeys}

\usepackage{AKstyle}

\numberwithin{equation}{section}

\usepackage{caption}
\usepackage[labelfont=rm]{subcaption}

\usepackage[pagebackref=true, colorlinks]{hyperref}

\hypersetup{pdffitwindow=true,linkcolor=blue,citecolor=blue,urlcolor=blue,menucolor=blue}

\usepackage{comment}

\begin{document}

\author{Alex Kontorovich}
\thanks{Kontorovich is partially supported by NSF grant DMS-1802119, BSF grant 2020119, and the 2020-2021 Distinguished Visiting Professorship at the National Museum of Mathematics.}
\email{alex.kontorovich@rutgers.edu}
\address{Department of Mathematics, Rutgers University, New Brunswick, NJ, and National Museum of Mathematics (MoMath), New York, NY}
\author{Xin Zhang}
\thanks{Zhang is partially supported by ECS grant 27307320 and NSFC grant 12001457.}
\email{xzhang@maths.hku.hk}
\address{Department of Mathematics, The University of Hong Kong, Pokfulam, Hong Kong}

\title
{On Length Sets of Subarithmetic Hyperbolic Manifolds}

\begin{abstract}
We formulate the Asymptotic Length-Saturation Conjecture on the length sets of closed geodesics on hyperbolic  manifolds whose fundamental groups are subarithmetic, that is, contained in an arithmetic group. We prove the first instance of the conjecture for punctured, Zariski dense covers of the modular surface.
\end{abstract}
\date{\today}
\maketitle

\section{Introduction}

By the length spectrum of  a hyperbolic manifold $M$, we mean the set of lengths of closed geodesics on $M$, with multiplicity. As is well-known, closed geodesics on $M$ correspond to hyperbolic conjugacy classes of its fundamental group 
$$
\G=\pi_1(M)<\Isom(\bH^n)\cong PO(n,1),
$$ 
and lengths of the former are a simple function of traces 
of the latter. It is also classical to study the length set, that is, the set of lengths of closed geodesics, now counted {\it without} multiplicity; again, this is intimately related to the set $\cT(\G)$ of traces (without multiplicity) of hyperbolic conjugacy classes of $\G$. In this paper, we initiate a detailed study of the latter for (sub)arithmetic groups, from the viewpoint of local-global theory.

To motivate our main results, we begin with a few illustrative examples. 

\noindent
{\bf Example 1:} Consider the Hecke $(2,3,\infty)$ triangle group, or rather, its orientation preserving double cover, the modular group $\G=\SL_2(\Z)$. The trace set $\cT(\G)$ of the latter is elementarily seen to be all of $\Z$, as for any desired integer $t$, one simply expresses $t$ as $t=a+d$ and factors $bc=ad-1$ to make a matrix $\mattwos abcd\in\G$ having trace $t$. This is because $\G$ is an arithmetic (or better yet, congruence) group, and hence {\it any} solution to $ad-bc=1$ over $\Z$ gives an element. 
We can compare these facts with the count for the number of points in $\G$ in an archimedean ball $B_N$ of radius $N$: as with any lattice in $\SL_2(\R)$, we crudely have
\be\label{eq:GBNrough}
\#\G\cap B_N  \ \sim \ c N^2,
\ee
for some constant $c>0$, which means that the average  number of times that a particular integer $t\asymp N$ arises as a trace of a matrix in $B_N$ is of order $N^2/N=N$. 
But this does not take into account the fact that trace is a conjugacy class invariant.
For $t>2$, let $H(t)$ denote the number of conjugacy classes of elements in $\G$ with trace $t$. As is well-known (see Chowla-Cowles-Cowles \cite{ChowlaCowlesCowles1980}, though this was surely known to Selberg and likely before), $H(t)$ is equal to $h(t^2-4)$, where $h(D)$ is the classical class number, that is,  the number of equivalence classes of binary quadratic forms of discriminant $D$ (not necessarily primitive). 
By the Prime Geodesic Theorem we have (compare to \eqref{eq:GBNrough}):
\be\label{Htrough}
\sum_{t<N}H(t) \ \sim \ {N^2\over2\log N},
\ee
and so a ``typical'' value of $H(t)$ (for $t\asymp N$) is more like $N/\log N$, rather than $N$. (Note that the fundamental unit $\epsilon_D$ for discriminant $D=t^2-4$ is about as small as possible, $\epsilon_D=(t+\sqrt D)/2$, and hence this class number is as large as possible, of size about $\sqrt D$).  The discrepancy in counting matrices versus counting conjugacy classes makes sense, as the archimedean size of elements under conjugation grow exponentially (the stabilizer group of a conjugacy class is a discrete subgroup of some $\SO(1,1)$), so $(\log N)$-many matrices of size $N$ are grouped together. This is a minor issue here, but will play a major role in the next example.

\

\noindent
{\bf Example 2:} Now consider the Hecke $(2,5,\infty)$ triangle group, or rather its double cover, the group $\G$ generated by 
$$
\G=\<\mattwo01{-1}0,\mattwo 1\phi01\>,
$$
where $\phi=(1+\sqrt 5)/2$ is the golden mean. (Recall that $\SL_2(\Z)$ has similar generators, except with $\phi$ replaced by $1$.) The group $\G$ is nonarithmetic, but it is {\it subarithmetic}, being a subgroup of the arithmetic group $\widetilde\G:=\SL_2(\Z[\phi])$. The latter does not act discretely on $\bH$, but {\it is} a lattice in $\SL_2(\R)\times \SL_2(\R)$, where it acts by the Galois conjugate in the second factor. (The fact that $\G$ is a lattice in the first factor, while the second factor is non-compact,  is one way to see its nonarithmeticity.) In fact $\G$ is a {\it thin} group (see, e.g. \cite{Kontorovich2014}), as the set of $\Z$-points of its Zariski closure is exactly $\widetilde\G$, and it has infinite index in the latter. 
The set of traces of $\widetilde\G$ is again elementary to determine; it is the full order $\cO=\Z[\phi]$, for the same reason as in Example 1. But now we may ask, which $t\in \cO$ are also traces of $\G$?

The asymptotic count \eqref{eq:GBNrough} of {\it matrices} in a ball is still of order $N^2$. But 
 now, since $\cO\cong \Z\oplus\phi\Z$ is a quadratic ring, the number of  $t\in \cO$ with norm at most $N$ is also roughly $N^2$. Therefore the average number of matrices in $B_N$ having a given trace $t\in\cO$ is a positive constant.
 
But what happens when we group by conjugacy classes? (In Example 1, this caused the average count to drop by a factor of $\log N$, but here we don't have this factor to spare!) 
Let $H_\G(t)$ denote the number of conjugacy classes of $\G$ having trace $t\in\cO$. As $\G$ is a lattice in $\SL_2(\R)$, we still have the Prime Geodesic Theorem (see \eqref{Htrough}), that
\be\label{eq:HGt}
\sum_{\N( t) < N } H_\G(t) \  \sim \ {N^2\over 2\log N},
\ee
where $\N:\cO\to\Z$ is the algebraic norm. Therefore there can't be more than $O(N^2/\log N)$ elements $t\in\cO$ which actually arise as traces in $\G\cap B_N$, and thus the density of those that do arise is zero! 
While we can't say much about the class number $H_\G(t)$, in every conjugacy class that does arise, there should be about $\log N$ matrices of size $N$, as before. So when counting matrices, even though the ``average'' multiplicity is bounded, what's really going on is that $100\%$ of the time, the multiplicity is exactly zero, and very rarely there are somewhat large (of size at least $\log N$) multiplicities. See also 
recent work of McMullen \cite{McMullen2020} in this direction. We remark that the number of elements of $\widetilde\G$ in a ball $B_N$ is roughly $N^4$, and $H_{\widetilde\G}(t)$ is roughly of order $t^2$.

\

\noindent
{\bf Example 3:} For our last example, consider the Hecke $(2,7,\infty)$ triangle group, or rather its double cover, the group $\G$ generated by 
$$
\G=\<\mattwo01{-1}0,\mattwo 1\eta01\>,
$$
where $\eta=\cos(2\pi/7)$. The ring $\cO=\Z[\eta]$ is cubic, and so the number of matrices of size $N$ in $\G$ is of order $N^2$, while the number of possible values of the trace up to $N$ is $N^3$, and hence it is clear that very few numbers in $\cO$ can occur as traces.

Returning to the general setting, in light of these examples, to be able to say anything about the length set of a general manifold $M$, we need some conditions. First we should assume that its fundamental group $\G$ it is subarithmetic (or else the traces could be completely random numbers), and let $\cO$ be the ring of integers of the trace field (that is, the field generated by the traces) of $\G$.

{\bf Obstruction 1:}
Let $\ga>0$ be the ``growth exponent'' of $\G$, in the sense that 
$$
\#\G\cap B_N = N^{\ga+o(1)}.
$$ 
(When $\G$ is a geometrically finite group and $\gd$ is the Hausdorff dimension of its limit set, then $\ga=2\gd$.)
As in Example 3 (and Example 2), to be able to study the length set, we  require that $\ga$ exceeds the rank of $\cO$.
One can think of this as an ``archimedean local obstruction.''

{\bf Obstruction 2:}
There are also potentially other local obstructions. Already in the case of a classical congruence group $\G(q):=\ker(\SL_2(\Z)\to\SL_2(q))$, only the numbers that are $2(\mod q)$ can arise as traces in $\G(q)$.
We say that $t\in\cO$ is {\bf admissible} if, for every ideal $\cI\subset\cO$, $t\in \cT(\G)\mod \cI$. 
(We remark that the set of admissibles  is ``easy'' to determine in practice, as follows from Strong Approximation for Zariski dense groups.)

{\bf Obstruction 3:}
There is one final archimedean local obstruction. Given any manifold, we can take a cover that destroys the systole; that is, what was the shortest closed geodesic need not remain closed under a cover, making the shortest length (that is, smallest trace) of such a moving target. So we should allow for finitely many ``small'' values of $\cO$ to not arise as traces.

We may now formulate our main conjectures.
\begin{Def}
With notation as above, we say that a subarithmetic hyperbolic group $\G$ {\bf length-saturates} if: every admissible $t\in\cO$ with sufficiently large norm arises in the trace set of $\G$.
\end{Def}
\begin{Def}
We say that $\G$ {\bf asymptotically length-saturates} if
\be\label{def:asymp}
{\#\cT(\G)\cap B_N \over \#\{t\in\cO\cap B_N \ : \ t \text{ is admissible\}}} \ \to \ 1,
\ee
as $N\to\infty$.
\end{Def}

Thus the modular group length-saturates, as in Example 1 (see also work of Marklof \cite{Marklof1996} studying distinct length sets of arithmetic 3-folds), while the Hecke $(2,5,\infty)$ group does not even asymptotically length-saturate (it fails Obstruction 1).

\begin{conj}[Asymptotic Length-Saturation]\label{conj:main}
Let $\G$ be a subarithmetic hyperbolic group, with growth exponent $\ga$ exceeding the rank of $\cO$. Then $\G$ asymptotically length-saturates.
\end{conj}

The stronger statement that with the same assumptions, $\G$ length-saturates is false. Indeed, already for certain cocompact arithmetic 2-folds corresponding to the norm-one elements of a quaternionic division algebra, the trace equation can cut out a ternary indefinite inhomogeneous quadric, which can exhibit infinitely many Brauer-Manin-type obstructions.

In this paper, we make the first progress towards \conjref{conj:main}, by proving asymptotic length-saturation for punctured, Zariski-dense, submodular groups.
\begin{thm}\label{thm:main1}
Any punctured, Zariski dense cover of the modular surface is asymptotically length-saturating. In fact, it is effectively so, in that the right hand side of \eqref{def:asymp} is $1+O(N^{-\vep})$ for some $\vep>0$.
\end{thm}

Note that we do not require geometric finiteness of $\G$, as we can pass, if needed, to a finitely generated but still Zariski dense and punctured subgroup of $\G$.\footnote{In dimension 2, as is the setting of \thmref{thm:main1}, geometric finiteness is equivalent to finite generation.}

\begin{rmk}
The question of length-saturation is closely related to the Local-Global Conjecture for Apollonian packings, Zaremba's Conjecture, and McMullen's Classical Arithmetic Chaos Conjecture (see, e.g. \cite{Kontorovich2013, Kontorovich2016, McMullenNotes} for discussions of these). Each of these problems amounts to understand the image of a linear form (which in the setting of this paper is the trace) of a Zariski dense subgroup (or sub-semigroup). In each of the previous cases, the expected multiplicity in a ball of size $N$ was some fixed positive power of $N$. 

Note that in our setting here, there is no restriction on the growth exponent $\ga$ of $\G$; indeed, the Hausdorff dimension $\gd$ of the limit set of $\G$, which can be as small as $\gd>1/2$ (due to the puncture, see \cite{Beardon1965}).
Counting with multiplicity, we have that:
$$
\#\{\g\in\G\cap B_N \ : \ \tr(\g)<N\}  \ = \ N^{2\gd+o(1)}.\qquad (N\to\infty)
$$
So the multiplicity of a typical trace $t\asymp N$ in the trace set $t\in\cT(\G)$ may be extremely small,
\be\label{eq:expMult}
\#\{\g\in\G\cap B_N \ : \ \tr(\g)=t\} \ \overset{?}{=} \ N^{2\gd-1+o(1)},
\ee
where $2\gd-1$ can be any quantity just above $0$, and yet  our methods produce a density-one set of traces in this setting.

Also note that the methods introduced in \cite{BourgainKontorovich2010} and applied to both the Zaremba \cite{BourgainKontorovich2014} and Apollonian \cite{BourgainKontorovich2014a} settings required the linear form to have a bilinear structure. That is, the linear form was of type: 
\be\label{eq:bilinear}
\g \mapsto \<v,\g w\>
\ee
 for some fixed vectors $v,w$. The trace is {\it not} of this form, and so the best one can currently say towards McMullen's conjecture is a  strong level of distribution \cite{BourgainKontorovich2018}. It is  not even currently known that  a positive proportion of numbers arises in the set of traces of a Zaremba-type semigroup (see \cite[\S3]{Kontorovich2016}).
For related work, see also the recent PhD thesis of Brooke Ogrodnik  \cite{Ogrodnik2021}.
\end{rmk}

\begin{rmk}
Returning to the setting of this paper, here are some further remarks: 
\begin{enumerate}
\item
 It is  sometimes possible to completely determine the trace set of $\G$, even if the latter is thin. For example, take the ``Lubotzky 1-2-3'' group, $\G=\<\mattwos 1301,\mattwos1031\>$. It is easy to see that every trace is $\equiv2(\mod9)$, and indeed the element $\mattwos1031\mattwos1301^n$ has trace $2+9n$, so all admissible traces are represented by this one arithmetic progression. 

\item 
Since we have assumed that $\G$ does contain a parabolic element, it is immediate that its trace set $\cT(\G)$ comprises a positive proportion of integers, since, as above, $\cT(\G)$ contains entire arithmetic progressions. 
Without assuming that $M$ is punctured,   current technology cannot not even produce  a positive proportion of traces!

 \item
 On the other hand, 
 an argument based on Furstenberg's topology on the integers shows that,
 if
there is even a single local-global failure (that is, an admissible $t$ not in  $\cT(\G)$), then
finitely many such arithmetic progressions cannot possibly cover even a density-1 subset of $\cT(\G)$. 
\end{enumerate}
\end{rmk}

\begin{figure}
\includegraphics[width=3in]{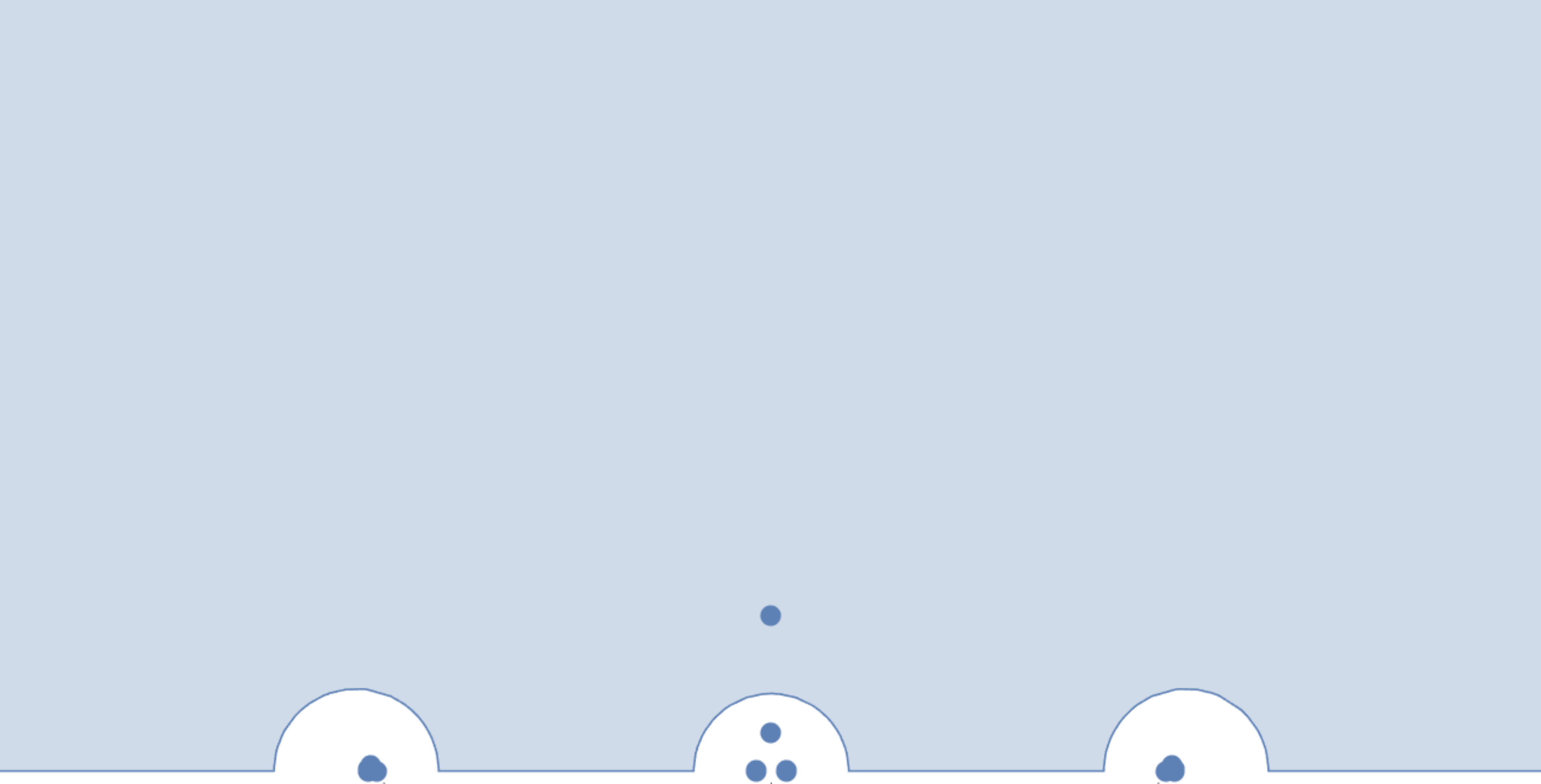}
\caption{The fundamental domain of $\G_0$, and a typical orbit.}
\label{fig:FunDom}
\end{figure}

 Note that it is easy to construct explicit  examples of groups with arbitrarily many local-global failures. For just one family of such, fix $m$ large, and consider the group $\G_0<\SL_2(\Z)$ generated by $S_1:=\mattwos0{-1}10$, 
$S_2:=\mattwos 1m01\mattwos0{-1}10\mattwos1{-m}01$, 
and
$S_3:=\mattwos 1{-m}01\mattwos0{-1}10\mattwos1{m}01$;
a fundamental domain for the action of $\G_0$ is shown in \figref{fig:FunDom}.
It is easy to see that the systole of these groups grows with $m$.
By strong approximation, there is some $q_0$ such that the reduction of $\G_0$ mod any prime $p\nmid q_0$ is onto. Let $P=P(m)$ be a very large prime coprime to $q_0$ and let $\G$ be the group generated by $\G_0$ and the translation $\mattwos 1P01$. Then the reduction of $\G$ mod any $q$ is clearly all of $\SL_2(q)$, so all numbers are admissible, while taking $P$ large enough does not create a shorter closed geodesic.

\subsection{Methods}\

We use the (orbital) circle method to access the trace set $\cT(\G)$. In fact, our methods apply not just to the trace function $\tr:\SL_2(\Z)\to \Z$ but to any linear form, $\sL$, say, on $\SL_2$ (and hence we do not group traces by conjugacy class).
It turns out (see \eqref{eq:sLdef}) that the trace function is the ``generic'' linear form.
The main theorem, from which \thmref{thm:main1} follows immediately, is the following.

\begin{thm}\label{thm:main2}
Let $\G$ be a geometrically finite, punctured, Zariski dense subgroup of $\SL_2(\Z)$, and let $\sL:\G\to\Z$ be any linear form. Let $\cA$ denote the admissible values of $\sL$; that is, $n\in\cA$ if and only if $n\in \G(\mod q)$, for all $q$. Then there is some $\gT>0$, so that:
$$
{\#\{n\in\sL(\G)\cap[1,N]\}\over \#\{n\in\cA\cap[1,N]\}} = 1+O(N^{-\gT}),
$$
as $N\to\infty$.
The implied constant is effective. 
\end{thm}

In the case that the linear form is bilinear (as in \eqref{eq:bilinear}) and the critical exponent of $\G$ is sufficiently close to $1$, the above theorem is proved by the second-named author in \cite{Zhang2020}.
One of the  key innovations in the present paper is to use the parabolic element in {\it two} ways to produce not only arithmetic progressions, but values of binary quadratic polynomials in the set of values of the linear form $\sL$.
By this we mean the following: given a fixed element $\mattwos abcd\in\G$ and a parabolic, say, $\mattwos1P01\in\G$, we can compute that
$$
\tr(\mattwos abcd \mattwos1P01^x)=a + d + c x P,
$$
which is a linear form in $x$, whereas, say,
$$
\tr(\mattwos1P01^x\mattwos abcd \mattwos1P01^y\mattwos abcd)=a^2+2 b c+d^2+(a + d)c P (x+y)+c^2  P^2xy
$$
is quadratic as a function of the pair $(x,y)$.
(Using three or more copies of the parabolic produces cubic or higher forms;   the added cost of increasingly larger coefficient sizes seems not to be advantageous for this problem.)

Then varying $x$ and $y$, and letting $\mattwos abcd$ run over certain regions of $\G$ in a ball of size $N$, we study the ``representation number'' $\cR_N(n)$ of the number of times that $n\asymp N$ occurs as a trace of such. In fact our construction of $\cR_N$ is more complicated, as we need to create multilinear forms at several scales for later estimates, see \secref{sec:setup}. 

Following the (orbital) circle method, we decompose $\cR_N(n)$ into a ``main'' term $\cM_N(n)$ and an ``error'', $\cE_N(n)$, where we integrate over the ``major'' and ``minor'' arcs, resp. 
We note again that the main term is expected to be of size a singular series times $N^{2\gd-1}$ (see \eqref{eq:expMult}), which may be an arbitrarily small (but fixed) positive power of $N$. So we do not have much room to get an error off of the main term!

As in some other applications of the orbital circle method, we are  only able to control the error terms in $L^2$, and hence produce only Asymptotic Length-Saturation, and not full Length-Saturation (which perhaps could be expected in this setting). Even this involves several novel techniques; there is a more standard analysis of cancellation in certain exponential sums and averages thereof, and there is also an appeal at some point to Hilbert's Nullstellensatz in effective form (see \secref{sec:Null}). 

In the major arcs, we use the work of Bourgain-Varju \cite{BourgainVarju2012} and Bourgain-Gamburd-Sarnak \cite{BourgainGamburdSarnak2011} for an archimedean spectral gap, together with infinite volume counting methods of \cite{BourgainKontorovichSarnak2010},  to obtain an estimate for the main term $\cM_N(n)$. But there are several surprises here as well!
It turns out that the singular series is a very short sum (of length $N^\vep$) which is trying to approximate a quadratic Dirichlet $L$-function at $1$, see \secref{sec:shortSum}. On GRH, this $L$-value can indeed be approximated by such a short sum, but our statement is unconditional! Since our error term estimate is anyway only as an ``average'' over $n$, we also average on the main term; that is, we show (see \thmref{thm:sENfSq}) that, for all but very few $n$'s, the approximation is valid. But then the have another problem: we need Siegel's bound to know that the singular series, which is now one such $L$-value, is not too small. Again, because we are stating only an average result, we show (see \thmref{thm:fSnSize}) that we can bound these $L$-values from below for all but an exceptional set of values of $n$, with {\it effective} constants, leading to the effective constants in \thmref{thm:main2}.
\pagebreak

\subsection*{Outline}\

We begin in \secref{sec:prelim} with the setup of the circle method, introducing the main representation function, and its decomposition into a main term and error, corresponding to major and minor arcs. 
The next two sections provide preparatory lemmata for the main arguments. We record in \secref{sec:preMaj} various infinite volume counting theorems in congruence towers, the savings off of such counts in progressions for {\it arbitrarily large} modulus (this is where Nullstellensatz is used), as well as the analysis of the singular series, which involves Weil and Burgess bounds. In \secref{sec:minTech}, we prepare various exponential sum estimates used in the minor arcs analysis, using more standard analytic techniques such as estimates for Kloosterman-type sums. These allow us to complete the major arcs analysis in \secref{sec:major} and then the minor arc analysis in \secref{sec:minor}.

\subsection*{Notation}\

We use the standard notation $e(x)$ we mean $e^{2\pi i x}$, and $e_q(x)=e(x/q)$. The notation $\sum_{r(q)}'$ means summing over $r(\mod q)$ with $(r,q)=1$. 
We use the symbols $X=O(Y)$ and $X\ll Y$ interchangeably, and by $X\asymp Y$, we mean $X\ll Y \ll X$. All implied constants,  unless specified otherwise, may depend at most on $\Gamma$ and the linear form $\sL$, which are treated as fixed.

\subsection*{Acknowledegements}\

We thank Jens Marklof, Curt McMullen, Lillian Pierce, Alan Reid, Zeev Rudnick, and Peter Sarnak for discussions and insightful suggestions about this work.

\newpage

\section{Preliminaries}\label{sec:prelim}

We henceforth take $\G<\SL_2(\Z)$ to be a given Zariski dense  group.
We consider a linear form
 $\sL:\SL_2(\Z)\to\Z$ which is not everywhere vanishing; explicitly, this means that
\be\label{eq:sLdef}
\sL \ : \  \mattwo abcd \  \mapsto \  Aa+Bb+Cc+Dd  \ = \  \tr\left[\mattwo abcd\mattwo ACBD\right],
\ee
and we assume that at least one coefficient $A,B,C,D$ is non-zero.
Note that $\sL$ is of bilinear type (see \eqref{eq:bilinear}) exactly when its ``discriminant'',
$$
\gD\ =\ \gD_\sL   \ := \ AD-BC
$$
vanishes.
After conjugation (possibly renaming $A,B,C,D$ above), we may assume that $\G$ contains the fixed parabolic element 
$$
\mattwo 1P01\in\G.
$$
We make a few further simplifying assumptions. 
\begin{itemize}
\item We may assume that $\gcd(A,B,C,D)=1$, since otherwise we can pull out a common factor.
\item By applying fixed elements of $\G$ inside $\sL$, we may assume that the coefficients $A,\ B,\ C$ and $D$ are all non-zero. Indeed, one has that:
$$
\sL\left(\g_1 \mattwo abcd \g_0\right)
=
\tr\left[\mattwo abcd 
\g_0
\mattwo ACBD
\g_1
\right]
,
$$
and by the Zariski density of $\G$, there exist  elements $\g_0,\g_1\in\G$ so that 
$
\g_0
\mattwos ACBD
\g_1
$
 has every entry non-zero.
\item
By Strong Approximation and passing to a finite index subgroup of $\G$ if necessary, we may assume that for all $(q_1,q_2)=1$, we have:
\be\label{eq:multG}
\G/\G(q_1q_2) \cong
\G/\G(q_1)
\times
\G/\G(q_2),
\ee
where
\be\label{eq:Gqdef}
\G(q) \ := \ \{\g\in\G : \g\equiv I(\mod q)\}
\ee
is the ``principal congruence'' subgroup of (the possibly thin group) $\G$.
Moreover, for all ``good'' primes $p$, we have that for $q=p^\ell$, the mod $q$ reduction is onto
$$
\G(q)\bk\G=\SL_2(q).
$$
\item
For a finite list of ``bad'' primes $p$ (including $p=2$), we have an exponent (``saturation level'') $k=k_p$ so that 
\be\label{eq:kpDef}
\G(p^k)\bk\G=\{I\},
\ee
and for $\ell>k$,
$\G(p^\ell)\bk\G$ is the full lift of the identity from $\SL_2(p^k)$ to $\SL_2(p^\ell)$. In particular, the parabolic element $\mattwos 1P01\in\G$ satisfies 
\be\label{eq:Pmodpk}
P\equiv 0(p^k)
\ee 
for all bad primes $p$.
\item
We may assume, by increasing the saturation densities $k_p$ if necessary, that 
\be\label{eq:kpLB}
k_p > L_{p,B},
\ee
where $L=L_{p,B}$ is determined by $p^L\| B$. (Here, as below, $B$ is the coefficient of $b$ in $\sL$, as defined in \eqref{eq:sLdef}.)
\end{itemize}

\subsection{Setup of the Circle Method}\label{sec:setup}\

For $\g\in\G$, we construct the shifted binary quadratic form:
\be\label{eq:linearform}
\ff_\g(x,y) \ := \ 
\sL\bigg(\mattwos1{Px}01 \g \mattwos1{Py}01\bigg),
\ee
so that, if $\g=\mattwos {a_\g}{b_\g}{c_\g}{d_\g}$, then
\be\label{eq:ffIs}
\ff_\g(x,y)
=
(A a_{\gamma }
+B b_{\gamma }
+C c_{\gamma }
+D d_{\gamma })
+(A c_{\gamma }
+B  d_{\gamma })
Px
+(B  a_{\gamma }
+D  c_{\gamma })
Py
+B   c_{\gamma }
P^2
x y
   .
\ee

Note that, for any integers $x,y\in\Z$ and any $\g\in\G$, the value of $\ff_\g(x,y)$ arises in $\sL(\G)$.
Let $N$ be the main growing parameter, and $T,X$ be parameters determined by: 
\be\label{eq:TX2eqN}
T=N^{1/100},\ 
X=N^{99/200}, \ \text{ so that }\
TX^2=N.
\ee
Decompose $T$ further at 
\be\label{eq:TT1T2}
T=T_1T_2,\ \text{ with } \
T_2 = T_1^\cC,
\ee
with $\cC$ a very large constant depending only on the spectral gap for $\G$, see \eqref{eq:cCdef}.

We now define the main ensemble $\sF_T$ as follows
\be\label{eq:sFTdef}
\sF_T := \left\{
 \g_1\cdot\g_2 =
\mattwo ab{c}{d} :
{\g_1,\g_2\in\G \atop
{\foh T_1\le \|\g_1\| < T_1 \atop
{\foh T_2\le \|\g_2\| < T_2 \atop
\frac1{100}T<a,b,c,d}}}
\right\}.
\ee
We show in \lemref{lem:cFTcount} that $\sF_T$  has cardinality $\asymp T^{2\gd}$.
This is a sub-multi-set of $\G$, as the product $\g_1\g_2$ may have multiplicity; that is, for a fixed $\xi\in\sF_T$,
\be\label{eq:sFTmult}
\sum_{\g\in\sF_T}\bo_{\{\g=\xi\}} \ll T_1^{2\gd}.
\ee
Fixing a smooth nonnegative bump function $\gU$ with $\supp\gU\subset[\foh,1]$,  define the main ``representation number''
$$
\cR_N(n)
\ := \
\sum_{\g\in \sF_T}
\sum_{x\in\Z}
\sum_{y\in\Z}
\gU\left({x\over X}\right)
\gU\left({y\over X}\right)
\bo_{\{\ff_\g(x,y)=n\}}
.
$$

We decompose $\cR_N$ into a ``main'' term and an error according to a (smoothed) major/minor arcs decomposition of the circle. To this end, let 
$$
\psi(t) \ := \ \max(1-|t|,0)
$$
be the ``tent'' function whose Fourier transform is the Fej\'er-type kernel:
$$
\widehat \psi(\xi) = {\sin^2 (\pi \xi)\over (\pi \xi)^2}.
$$
We fix parameters $Q_0, K_0$ to be determined as follows. We set
\be\label{eq:Q0K0def}
Q_0=N^{\ga_0}, \qquad
K_0=N^{\gk_0},
\ee
where the exponents satisfy
\be\label{eq:ga03gk0}
\gk_0=3\ga_0
\ee
and
\be\label{eq:ga0gk0}
5\ga_0+\gk_0<\gT.
\ee
Here $\gT$ is the minimum of the two values in \lemref{lem:Expansion} and \lemref{lem:Arch}.
Setting
$$
\ga_0=\gT/10, \qquad \gk_0=3\gT/10
$$
will satisfy all the criteria.

With these choices, let
$$
\Psi_{N,K_0}(\gb):=
\sum_{m\in\Z}
\psi((\gb+m)\tfrac{N}{K_0}),
$$
and define the ``major arcs'' weight function as:
\be\label{eq:fMdef}
\fM(\gt)
\ := \
\sum_{q<Q_0}
\sum_{r(q)}'
\Psi_{N,K_0}(\gt-\tfrac rq)
.
\ee
Then the ``main'' term is given by:
\be\label{eq:cMdef}
\cM_N(n) \ := \
\int_0^1
\fM(\gt)
\widehat{\cR_N}(\gt)
e(-n\gt)
d\gt,
\ee
and of course the error is
$$
\cE_N(n) \ := \
\int_0^1
(1-\fM(\gt))
\widehat{\cR_N}(\gt)
e(-n\gt)
d\gt,
$$
so that
$$
\cR_N(n) = \cM_N(n)+\cE_N(n).
$$

\newpage

\section{Major Arc Technical Estimates}\label{sec:preMaj}

We record here a number of technical estimates needed in the  analysis of the main term.

\subsection{Spectral Analysis and Counting}\

Let $\gd=\gd_\G$ be the Hausdorff dimension of the limit set of $\G$, and recall that $\gd>1/2$. By Patterson-Sullivan theory \cite{Patterson1976, Sullivan1984}, $\gd$ is related to the bottom eigenvalue $\gl_0=\gd(1-\gd)$ of the hyperbolic Laplacian acting on $L^2(\G\bk\bH)$. Work of Lax-Phillips \cite{LaxPhillips1982} shows that the spectrum of the latter below $1/4$ consists of finitely many eigenvalues. By  Bourgain-Varju and Bourgain-Gamburd-Sarnak there is a uniform spectral gap, in the following sense.

\begin{thm}[{\cite{BourgainVarju2012, BourgainGamburdSarnak2011}}]
There exists a ``spectral gap,'' 
\be\label{eq:gT0}
\gT_0=\gT_0(\G)>0
\ee 
so that, for all $q\ge1$, 
the eigenvalue $\gl_0$ is the bottom of the spectrum of $L^2(\G(q)\bk\bH)$, and all other eigenvalues are at least $\gl_\gT:=s_\gT(1-s_\gT)$, where $s_\gT:=\gd-\gT_0$.
Here $\G(q)$ is as defined in \eqref{eq:Gqdef}.
\end{thm}

Recalling the construction of $\sF_T$ from \eqref{eq:sFTdef}, we record the following counting results, which follow from now-standard techniques.
\begin{lemma}\label{lem:cFTcount}
As $T\to\infty$,
$$
\#\sF_T  \ \asymp \ T^{2\gd}.
$$
\end{lemma}
\pf
This follows from infinite volume counting methods in Zariski dense groups with $\gd>1/2$; see, e.g., \cite{BourgainKontorovichSarnak2010}.
\epf
\begin{lemma}\label{lem:Expansion}
There exists $\gT>0$ so that, for any $q\ge1$, $\g_0\in\G(q)\bk\G$, $|\gb|<1/X^2$, and $x,y\asymp X$, we have that
$$
\sum_{\g\in \sF_T\atop \g\equiv \g_0(\mod q)}
e(\gb \ff_\g(x,y))
=
\frac1{[\G:\G(q)]}
\sum_{\g\in \sF_T}
e(\gb \ff_\g(x,y))
+
O(|\sF_T|N^{-\gT}), 
$$
as $T\to\infty$.
\end{lemma}

\pf
The proof is the same as that of \cite[Theorem 1.14]{BourgainKontorovichSarnak2010}.
\epf

\begin{lemma}\label{lem:Arch}
There exists $\gT>0$ so that, for $x,y\asymp X$ and $n\asymp N$, we have:
$$
\sum_{\g\in \sF_T}
\int_\R
\psi(\gb\tfrac{N}{K_0})
e(\gb (\ff_\g(x,y)-n))
d\gb
\gg
{|\sF_T|\over K_0}
+
O(|\sF_T| N^{-\gT})
,
$$
as $T\to\infty$.
\end{lemma}

\pf
We first note that
$$
\int_\R
\psi(\gb\tfrac{N}{K_0})
e(\gb (\ff_\g(x,y)-n))
d\gb
=
{K_0\over N}
\widehat\psi((\ff_\g(x,y)-n){K_0\over N})\ge 0,
$$
and if $|\ff_\g(x,y)-n){K_0\over N}|<\frac12$, then $\widehat\psi(\cdot)>\frac25$. So we need to show the count:
$$
\sum_{\g\in \sF_T}
\bo_{\{|\ff_\g(x,y)-n|<\frac{N}{2K_0}\}}
\gg
{|\sF_T|\over K_0}
+
O(|\sF_T| N^{-\gT}).
$$
The latter follows from the same techniques as the proof of \cite[Theorem 1.15]{BourgainKontorovichSarnak2010}.
\epf

\subsection{Nullstellensatz}\label{sec:Null}\

\begin{thm}\label{thm:sFTqBnd}
Let $\gT_0$ be the spectral gap in \eqref{eq:gT0}.  Define $\cC$ by
\be\label{eq:cCdef}
\cC=3\times 10^8/\gT_0.
\ee
There exists an $\eta_0>0$ depending only on the spectral gap for $\G$, so that, for all $1\le q<N$, and all $r(\mod q)$, 
\be\label{eq:sFTqBnd}
\sum_{\g\in\sF_T}\bo_{\{c_\g \equiv r(\mod q) \}} \  \ll \ \frac1{q^{\eta_0}}|\sF_T|.
\ee
\end{thm}
\pf
We first drop the condition $\frac1{100}T<a,b, c,d$ from $\sF_T$ in \eqref{eq:sFTdef}, so that we need to count the number of $\g_1\asymp T_1$, $\g_2\asymp T_2$ so that the ``$c$'' entry of $\g_1\g_2$, 
$$
\<e_2,\g_1\g_2e_1\>\equiv r(\mod q)
,
$$ 
where $e_j$ are standard basis vectors.
This decomposes into two cases according to the size of $q$.

{\bf Case $q<T_2^{\gT_0/3}$:} 
In this case, we simply apply spectral theory in $\g_2$ while leaving $\g_1$ fixed, as follows. Break $\g_2$ into progressions mod $q$:
\beann
&&
\sum_{\g_1\asymp T_1}
\sum_{\g_2\asymp T_2}
\bo_{\{\<e_2,\g_1\g_2e_1\>\equiv r(\mod q) \}} 
\\
&=&
\sum_{\g_1\asymp T_1}
\sum_{\g_0\in \G(q)\bk\G}
\left[
\bo_{\{\<e_2,\g_1\g_0e_1\>\equiv r(\mod q) \}} 
\sum_{\g_2\asymp T_2\atop \g_2\equiv\g_0(\mod q)}
1
\right].
\eeann
The bracketed term may be estimated using the uniform spectral gap (see  \cite{BourgainKontorovichSarnak2010}) to give
\beann
&\ll&
\sum_{\g_1\asymp T_1}
\sum_{\g_0\in \G(q)\bk\G}
\bo_{\{\<e_2,\g_1\g_0e_1\>\equiv r(\mod q) \}} 
\left[
\frac1{q^3}
T_2^{2\gd}
+
O(T_2^{2\gd-\gT_0})
\right]
\\
&\ll&
\frac1{q}
T_2^{2\gd}
T_1^{2\gd}
+
q^2T_1^{2\gd}T_2^{2\gd-\gT_0}.
\eeann 
Here we used that $[\G:\G(q)]\asymp q^3$.
This saves $1/q$ (more than claimed) as long as $q<T_2^{\gT_0/3}$.

{\bf Case $q\ge T_2^{\gT_0/3}=T_1^{10^{8}}$:} The overview of the argument is as follows. For any fixed $\g_2$, we consider the set of $\g_1\asymp T_1$ for which $\<e_2,\g_1\g_2e_1\>\equiv r(\mod q)$. Since different integers having the same residue class mod $q$ differ by $q$, and $q$ is huge compared to $T_1$, we will show by Nullstellensatz that in fact the modular restriction can be lifted to an absolute restriction $\<e_2,\g_1\g_2e_1\>=r_*$, for some integer $r_*$ (depending on $\g_2$, which is fixed). Then we will relax the absolute restriction back down to a modular one, but with a much smaller modulus, $\<e_2,\g_1\g_2e_1\>\equiv r_*(\mod q_*)$, where $q_*\asymp T_1^{\gT_0/3}$, and apply the previous argument to save a power of $q_*$, which itself is a tiny power of $q$. 

An issue arises in the use of Nullstellensatz that was overlooked in related arguments in \cite{BourgainKontorovich2014a, BourgainKontorovich2015}. Write $\g_2e_1=(u,v)$ and $\g_1e_1=(c,d)$, so that $\<e_2,\g_1\g_2e_1\>=uc+vd$, with $|u|,|v|\le T_2$ being ``large'' and fixed, and $|c|,|d|\le T_1$ being ``small'' variables. It was claimed that, since $(u,v)=1$, we may assume that, say, $(u,q)=1$, and rewrite the modular condition as $c+v\bar u d \equiv r \bar u(\mod q)$. Unfortunately the obvious linear transformation that allows this rewrite requires changing the coefficients $c,d$ to ones of size bounded by $T_1T_2=T$, and this ruins the heights of the polynomials to be used in effective Nullstellensatz. So we need a more delicate argument to control the size of coefficients, as follows.

Suppose $q<N$ has a divisor $\tilde q\mid q$  of size $T_1<\tilde q<T_2^{\gT_0/3}$, say. Then we relax $\<e_2,\g_1\g_2e_1\>\equiv r(\mod q)$ to the same congruence mod $\tilde q$, and count as in the previous case. This saves $1/\tilde q>1/T_1$, which is a (very small) power of $N>q$, and completes the argument in this case.

Next we suppose that $q$ has no divisor in this range. Let $\tilde q$ be the largest divisor of $q$ not exceeding $T_1$, and write $q_0:=q/\tilde q$. 
We again relax the congruence restriction to $\<e_2,\g_1\g_2e_1\>\equiv r(\mod q_0)$; if we can save a small power of $q_0$, this also saves a small power of $q$.
Then any prime divisor $p$ of $q_0$ must exceed $T_2^{\gT_0/3}$, for otherwise either $p$ or $p\tilde q$ is a divisor of $q$ which is does not exceed $T_2^{\gT_0/3}$. 
Therefore $q_0$ is ``almost-prime'', that is, there are primes $p_j\ge T_2^{\gT_0/3}=N^\eta$, say, so that
$$
q_0=p_1p_2\cdots p_\ell,
$$
with $\ell<\lceil1/\eta\rceil$.

Next we consider the values of $u+\ga v$, for $\ga=1,2,\dots,\ell+1$, and claim that at least one such value is coprime to $q_0$. (The point here is that $\ell$ depends only on $q$, and is bounded only in terms of $\gT_0$, which only depends on $\G$.) 
Consider first the primes $p_j$ which divide either $u$ or $v$ (recall that $u$ and $v$ are coprime); then since $p_j>T_2^{\gT_0/3}> \ell+1$ (for $N$, and hence $T_2$, large enough), none of the $p_j$ can divide {\it any} value of $u+\ga v$. 
Now consider the $p_j$ which are coprime to $u$ and $v$. Then since $u+\ga v$ is an arithmetic progression of length $\ell+1<p_j$, at most one value $\ga_j\in\{1,2,\dots,\ell+1\}$ can satisfy $u+\ga_j \equiv0(\mod p_j)$.
Since the number of $\ga$'s exceeds the number of $p_j$'s, there is some $\ga$ so that $u+\ga v$ is coprime to all the $p_j$, and hence coprime to $q_0$. 

Again, this $\ga$ is bounded absolutely, and depends only on $q$, $u$, and $v$, and not on $c$ and $d$ (which depend on $\g_1$). Now we proceed with the Nullstellensatz argument.
Using the modulus $q_0$, we fix $\g_2$, let $(u,v)=\g_2e_1$, and consider the set 
$$
S=S_{\g_2}:=\{\g_1\in\G,\ \|\g_1\|\le T_1\ \text{ with }\ uc+vd\equiv r(\mod q_0)\},
$$
where we have set $(c,d):=\g_1e_2$. Using $\ga$ from the previous argument with $u+\ga v$ coprime to $q_0$, we write $uc+vd=(u+\ga v)c+v(d-\ga c)$, so that the congruence condition becomes
$$
c + v \overline{(u+\ga v)}(d-\ga c)\equiv r \overline{(u+\ga v)} (\mod q_0).
$$
Now consider the (linear) polynomials $P_{\g_1}\in \Z[U,V]$ given by 
$$
P_{\g_1}(U,V) := c + U (d-\ga c) - V,
$$
and consider the affine variety
$$
\cV := \bigcap_{\g_1\in S} \{P_{\g_1}=0\}.
$$

We claim that $\cV(\C)$ is nonempty. 
Note that the coefficients of $P_{\g_1}$ are bounded by $(\ell+2)T_1$.
Then if $\cV(\C)$ is empty, Hilbert's Nullstellensatz, in effective form (see, e.g., \cite[Theorem IV]{MasserWustholz1983})
gives the existence of polynomials $Q_{\g_1}\in \Z[U,V]$ and an integer $\fd\ge1$ so that
\be\label{eq:Null}
\sum_{\g_1\in S}P_{\g_1}(U,V) Q_{\g_1}(U,V) \ = \ \fd,
\ee
and with $\fd$ bounded (for $N$, and hence $T_1$, large enough) by 
$$
\fd \le \exp(8^{7}(\log T_1 + \log (\ell+2) + 8 \log 8) \le T_1^{10^7}.
$$
(``Large enough'' is in terms of  an implied constant depending only on $\G$, since $\ell$ depends only on $\gT_0$).
But if we reduce \eqref{eq:Null} mod $q_0$ and set $(U,V)\equiv( v \overline{(u+\ga v)} , r \overline{(u+\ga v)})$, we get $\fd\equiv 0 (\mod q_0)$, which is impossible since $q_0=q/\tilde q > T_1^{10^8-1}$.

Therefore $\cV(\C)$ is nonempty, and hence $\cV(\Q)$ is nonempty, and so clearing denominators, there exist coprime integers $u_*,v_*,r_*$, so that for all $\g_1\in S$,
$$
u_*c + v_* d = r_*.
$$
We have turned our congruence condition into an archimedean condition. Now we take some $q_*\asymp T_1^{\gT_0/3}$ coprime to $u_*,v_*,r_*$, relax the archimedean condition back to a modular one, 
$
u_*c + v_* d \equiv r_*(\mod q_*),
$
and count the number of $\g_1\asymp T_1$ satisfying this. As before, the spectral argument saves $1/q_*$, which is some small power of $q$.
\epf

\subsection{Singular Series Preliminaries}\

Recall that $\gcd(A,B,C,D)=1$ and $\gD=AD-BC$. Let $c_q$ denote the Ramanujan sum. (There should be no confusion between $c_q$ and bottom left element $c=\g_c$ of a matrix $\g=\mattwos abcd$.) We study here sums arising in the singular series analysis, of the form
$$
\fS_q(n):=
\frac1{|\G(q)\bk \G|}
\sum_{\g\in \G(q)\bk \G}
c_q (\ff_{\g}(x,y)-n),
$$
for fixed $x,y\in\Z$.
Note immediately that the sum, being over all $\g\in\G(q)\bk\G$, is independent of $x,y$, which we may assume are both $0$; thus $\ff_\g=Aa+Bb+Cc+Dd$. 
By the structure of $\G(q)\bk\G$ in \eqref{eq:multG}, the sum is multiplicative, so we may assume that $q=p^\ell$ is a prime power. For ``good'' primes, we have that $\G(q)\bk\G=\SL_2(q)$, and for ``bad'' primes, $\G(\mod p^k)=\{I\}$ for some ``saturation exponent'' $k$, while $\G(\mod p^\ell)$ for $\ell>k$  is the full lift to $\SL_2(p^\ell)$ of the identity in $\SL_2(p^k)$. 

\subsection{Good Primes}\

\begin{lemma}\label{lem:fSqgoodprime}
Assume that $q$ is a power of a good prime. Then we have that
$$
\fS_q(n)=
\frac1{|\G(q)\bk \G|}
\sum_{\g\in \G(q)\bk \G}
c_q (a_\g+\gD d_\g-n).
$$
\end{lemma}
\pf
Recalling that 
$$
\ff_\g(0,0)=\tr\left[\g\mattwos ACBD\right]=\tr\left[\mattwos\ga\gb\gk\gd \g \mattwos ACBD \mattwos\gd{-\gb}{-\gk}\ga \right],
$$
for any $\mattwos\ga\gb\gk\gd\in\SL_2$, and the sum being over all $\g\in\SL_2(q)$, we may simplify the expression in the following way. Assume WLOG that $(A,q)=1$; then we can rescale $\mattwos ACBD$ to $\mattwos 1{C\bar A}{BA}{DA}$, and continuing by elementary operations, we may replace $\mattwos ACBD$ by $\mattwos100\gD$, as claimed.
\epf

\subsubsection{Case $\gD\equiv 0 (\mod p)$}\

\begin{lemma}\label{lem:gD0p}
Assume $\gD\equiv0(p)$.
For $q=p$ a good prime, we have that:
$$
\fS_q(n)
=
\twocase{}
{-1\over p+1}{if $n\equiv 0(p)$,}
{1\over p^2-1}{if $n\not\equiv0(p)$.}
$$
\end{lemma}

\pf
Recall that 
$$
c_p(x)=\twocase{}{p-1}{if $x=0$,}{-1}{else.}
$$
Write $\g=\mattwos abcd$. Assume $p$ is a good prime. From \lemref{lem:fSqgoodprime}, we need to count the number of $\mattwos abcd$ with $a=n$ or $a\neq n$.

Consider the case $n\equiv0(\mod p)$. Then either $a\equiv n\equiv 0(\mod p)$ or not. In the first case, $bc\equiv 1(\mod p)$ and $d$ is free  ($p(p-1)$ matrices) and $c_p=p-1$ for a total contribution of $p(p-1)^2$. In the second case $a\neq n$, there are $p-1$ choices for $a$, then $p^2$ choices for $b,c$, and $d=(bc+1)\bar a$ is determined. This is $p^2(p-1)$ matrices with $c_p=-1$. Combining these contributions gives $(-1)p(p-1)$ when $n\equiv 0$.

Now suppose $n\not\equiv0$. Then if $a\equiv n$, then $b $ and $c$ are free (with $p^2$ choices) and $d$ is determined, with $c_p=p-1$, for a net contribution of $p^2(p-1)$. If $a\not\equiv n$, then $c_p=-1$ and we either have $a=0, bc\equiv1$ and $d$ free ($p(p-1)$ choices), or $a\neq0$ (with $p-2$ choices), and $b,c$ free and $d$ determined ($p^2$ choices). The total contribution is then $p$ when $n\not\equiv0(p)$.

The size of $\SL_2(p)=p(p-1)(p+1)$, which gives the claim.
\epf

\begin{lemma}\label{lem:pell}
Assume $\gD\equiv0(\mod p)$.
For $q=p^\ell$ a power of a good prime ($\ell\ge2$), we have that:
$$
\fS_q(n)=0.
$$
\end{lemma}
\pf
For prime powers, we have that:
$$
c_{p^\ell}(x)=\threecase{0}{if $x\not\neq 0(p^{\ell-1})$,}
{-p^{\ell-1}}{if $x\not\equiv0(p^\ell)$ but $x\equiv 0 (p^{\ell-1})$,}
{p^{\ell-1}(p-1)}{if $x\equiv 0 (p^{\ell})$.}
$$
So there is no contribution unless $a+\gD d\equiv n(\mod p^{\ell-1})$. Fix $\g_0\in\SL_2(p^{\ell-1})$ which solves $a_0+\gD d_0\equiv n(\mod p^{\ell-1})$ and $a_0d_0-b_0c_0=1 
$. Consider any lift $\g\in\SL_2(p^\ell)$ of $\g_0$, that is, $a=a_0+p^{\ell-1}a_1$, etc, with the restriction that $ad-bc\equiv1(\mod p^\ell)$, that is,
\be\label{eq:pEll}
a_1d_0+d_1a_0-c_1b_0-b_1c_0\equiv 0 
(\mod p).
\ee
(This is just the Jacobian of the determinant.) 
The above defines a 3-dimensional subspace restricting the values of $a_1,b_1,c_1,d_1$ (this is the Lie algebra). 
Assume WLOG that $a_0\neq0(\mod p)$. Then \eqref{eq:pEll} determines $d_1$ once $a_1,b_1,c_1$ are determined.
We consider two cases, $a+\gD d\equiv n(\mod p^\ell)$ or not; since $\gD\equiv0(\mod p)$, this is a restriction on $a_1$, which leaves $b_1,c_1$ free ($p^2$ choices, which is the same count either way). 
If $a_1$ is the unique value mod $p$ for which $a+\gD d\equiv n(\mod p^\ell)$, then $c_{p^\ell}=p^{\ell-1}(p-1)$.
But if $a_1$ is one of the $(p-1)$ values for which $a+\gD d\not\equiv n(\mod p^\ell)$, then $c_{p^\ell}=-p^{\ell-1}$.
The net contribution from these two cases exactly cancels.
\epf

\subsubsection{Case $\gD\not\equiv0(\mod p)$}\

\begin{lemma}\label{lem:gDneq0}
Assume $\gD\not\equiv 0(\mod p)$.
For $q=p$ a good prime, we have that:
$$
\fS_q(n)
=
{1+p\left({n^2-4\gD\over p}\right) \over p^2-1},
$$
where $\left({\cdot\over p}\right)$ is the Legendre symbol.
\end{lemma}

\pf
Again by \lemref{lem:fSqgoodprime}, we need to count the number of $\mattwos abcd$ with $a+\gD d=n$ or not.
We decompose $\SL_2(p)$ according to whether $c=0$ or not.

If $c=0$, then $\g=\mattwos a b0{\bar a}$, and we need to know whether $a+\gD \bar a\equiv n$ or not. This equation is equivalent to $a^2-na+\gD\equiv0$, which has $\left({n^2-4\gD\over p}\right)+1$ solutions for $a$ with $b$ free (with $p$ choices), each contributing $c_p=p-1$ to $\fS_q$. 
The remaining $\left(p-1 - \left({n^2-4\gD\over p}\right)-1\right)p$ solutions contribute $c_p=-1$ each.

If $c\not\equiv 0$, then for any choice of $d$, we either have $a\equiv n-\gD d$ (with one choice, contributing $c_p=p-1$) or not ($p-1$ choices contributing $c_p=-1$). Then $c$ is free ($p-1$ choices) and  $b=(ad-1)\bar c$ is determined. These two contributions exactly cancel.

On using $|\SL_2(p)|=p(p-1)(p+1)$, the net contribution is as claimed.
\epf

\begin{lemma}\label{lem:gDne0Pell}
Assume $\gD\not\equiv0(\mod p)$.
Let $p^L\| (n^2-4\gD)$.
For $q=p^\ell$ a power of a good prime ($\ell\ge2$), we have that:
\be\label{eq:gDne0Pell1}
\fS_q(n)=
\fourcase
{0}{if $L\le \ell-2$, or if $\ell$ is odd and $L\ge \ell$,}
{{p^{-(\ell-3)/2}
}
( p^2-1)^{-1}
\left({(n^2-4\gD)/p^L\over p}\right)
}{if $\ell$ is odd and  $L= \ell-1$,}
{p^{-(\ell-2)/2}({ p+1})^{-1}}{if $\ell$ is even and $L\ge\ell$,}
{-p^{-(\ell-2)/2} (p^2-1)^{-1}}{if $\ell$ is even and $L=\ell-1$.}
\ee
In any case,
\be\label{eq:gDne0Pell}
\fS_{p^\ell}(n) \ll p^{-\ell/2}.
\ee
\end{lemma}

\pf
We decompose $\SL_2$ according to the value of $\g_c$:
$$
\SL_2(p^\ell)
 \ = \
 \bigsqcup_{c\in\Z/p^\ell}
 \sC_c,
$$
where
$$
\sC_c=\{\g\in\SL_2(p^\ell):\g_c=c\}.
$$
Notice that $\G_\infty=\{n_x \ : \ x\in\Z/p^\ell\}$ acts on the left on $\sC_c$, where $n_x=\mattwos1x01$, so we may decompose $\sC_c$ into $\G_\infty$-cosets. The value $\ff_\g=a+\gD d$ changes to $\ff_{n_x\g}=a+\gD d + cx$ when $\g$ is replaced by $n_x\g$.
If $c\not\equiv0(\mod p^\ell)$, then this is an arithmetic progression as $x$ varies (otherwise, it is constant). For some $\G_\infty$-cosets, the values of this progression are never $\equiv n(\mod p^{\ell-1})$, in which case there is no contribution to $\fS_q$ since the Ramanujan value $c_{p^\ell}$ vanishes. If the progression does attain the value $n (\mod p^{\ell-1})$, then as $x$ ranges mod $p^\ell$,  this value in $\Z/p^{\ell-1}$ is attained with equal multiplicities from its $p$ lifts in $\Z/p^\ell$. Exactly one of these lifts is $\equiv n(\mod p^\ell)$, which contributes $c_{p^\ell}=p^{\ell-1}(p-1)$, and the other $(p-1)$ lifts contribute $c_{p^\ell}=-p^{\ell-1}$. The two types of contributions exactly cancel.

We are left to study the distribution of the values $a+\gD\bar a$ from $\g=\mattwos a b 0{\bar a}$ ranging in $\sC_0$. 
In particular, we need only consider the values $a+\gD\bar a\equiv n(\mod p^{\ell-1})$ and determine which of these are also $\equiv n (\mod p^\ell)$.
The equation
$$
a+\gD\bar a \equiv n(\mod p^\ell)
$$
is equivalent to
$$
a^2-na+\gD \equiv 0(\mod p^\ell),
$$
which on completing the square gives the equation:
\be\label{eq:073001}
(a-\bar 2 n)^2
\equiv 
\bar 4(n^2-4\gD)
(\mod p^\ell).
\ee
We want to consider the number of solutions to \eqref{eq:073001} as compared to the solutions to the same equation but mod $p^{\ell-1}$:
\be\label{eq:073002}
(a-\bar 2 n)^2
\equiv 
\bar 4(n^2-4\gD)
(\mod p^{\ell-1}).
\ee
Consider first solutions to \eqref{eq:073002}.
If $n^2-4\gD$ is not a square mod $p^{\ell-1}$, then \eqref{eq:073002} has no solutions, and $\fS_q=0$.
Assume henceforth that $n^2-4\gD$ is a square mod $p^{\ell-1}$.
Let $p^L\| n^2-4\gD$.

{\bf Case: } $L\le \ell-2$. 

If $L$ is odd, then \eqref{eq:073002} has no solutions. So we assume that $L=2L_1$ is even. Since $n^2-4\gD$ is a square, we can thus write
 $n^2-4\gD\equiv s^2 p^{2L_1}(\mod p^{\ell-1})$ for some $s\not\equiv0(\mod p)$. Then \eqref{eq:073002} becomes:
$$
(a-\bar 2 n
-
\bar 2 sp^{L_1})
(a-\bar 2 n
+
\bar 2 s p^{L_1})
\equiv 
0
(\mod p^{\ell-1}).
$$
This equation is equivalent to the existence of $U,V\le \ell-1$ with $U+V\ge \ell-1$ such that 
$$
a-\bar 2 n
-
\bar 2 sp^{L_1}
\equiv0(\mod p^U),
\
a-\bar 2 n
+
\bar 2 sp^{L_1}
\equiv0(\mod p^V).
$$
Assume WLOG that $U\le V$.
Then taking the difference of these equations, we have that $U\ge L_1$.
But since $s$ is invertible mod $p$, we must also have $U\le L_1$, that is,
 $U=L_1$,
 and only the equation mod $p^V$ needs to be solved, which is solved uniquely.
Thus are then  $2p^{L_1}$ solutions to 
\eqref{eq:073002},
which are all of the form:
$$
a= \bar 2 n \pm \bar 2 s p^{L_1} + k p^{\ell-1-L_1},
$$
as $k$ ranges in $\Z/p^{L_1}$.

For each such value of $a$, the question becomes: which of these  also solves \eqref{eq:073001}? Letting $k$ range in $\Z/p^{L_1+1}$ and inserting this expression for $a$ into \eqref{eq:073001}, we get:
$$
(a-\bar 2 n)^2
\equiv
\bar 4 s^2 p^{2L_1} 
\pm s k p^{\ell-1}
\overset{?}{\equiv}
\bar 4(n^2-4\gD)
(\mod p^\ell),
$$
where we used that $2(\ell-1-L_1)\ge \ell$.
Since $s\not\equiv0(\mod p)$, as $k$ ranges in $\Z/p^{L_1+1}$, 
the values of $\bar 4 s^2 p^{2L_1} 
\pm s k p^{\ell-1}$
range in an arithmetic progression of step size $p^{\ell-1}$, and so are periodic, taking each value with equal probability. 
As before, the corresponding Ramanujan values are then such that the contributions to $\fS_q$ exactly cancel.

{\bf Case $L\ge \ell-1$ and $\ell$ even:}\

In this case, \eqref{eq:073002} asks for $(a-\bar 2 n)^2\equiv 0(\mod p^{\ell-1})$. 
The solutions to this are 
$$
a=\bar 2 n + k p^{\ell/2},
$$
as $k$ ranges in $\Z/p^{\foh\ell-1}$.

To see which solutions lift to  \eqref{eq:073001}, we let $k$ range in $\Z/p^{\ell/2}$. Then  \eqref{eq:073001} asks whether
$$
(a-\bar2 n)^2\equiv
k^2 p^{\ell}
\overset{?}{\equiv}
\bar 4(n^2-4\gD)
(\mod p^\ell)
.
$$
If $n^2-4\gD\equiv0(\mod p^\ell)$, that is, $L\ge \ell$, then every solution to \eqref{eq:073002} also solves  \eqref{eq:073001}.
So there are $p^{\ell/2}$ values of $a$ in $\mattwos a b 0 {\bar a}$, and another $p^{\ell}$ values of $b$ which is free. Each such matrix has a  Ramanujan value $c_q=p^{\ell-1}(p-1)$, for a net contribution of:
$$
\fS_q = {p^{(2-\ell)/2}
\over p+1
}
,
$$
where we used that $|\SL_2(p^\ell)|=p^{3(\ell-1)}p(p+1)(p-1)$.

If $n^2-4\gD\not\equiv0(\mod p^\ell)$, that is, $L= \ell-1$, then no solution to \eqref{eq:073002} lifts to  \eqref{eq:073001}. Each matrix as above has a Ramanujan value of $c_q=-p^{\ell-1}$, for a net contribution of:
$$
\fS_q =-{ p^{(2-\ell)/2} \over p^2-1}.
$$

{\bf Case $L\ge \ell-1$ and $\ell$ odd:}\

Now the solutions to \eqref{eq:073002}  are:
$$
a=\bar 2 n + k p^{(\ell-1)/2},
$$
as $k$ ranges in $\Z/p^{(\ell-1)/2}$.

Inserting these values into  \eqref{eq:073001}  and letting $k$ range in $\Z/p^{(\ell+1)/2}$, we are asking whether
$$
(a-\bar2 n)^2\equiv
k^2 p^{(\ell-1)}
\overset{?}{\equiv}
\bar 4(n^2-4\gD)
(\mod p^\ell)
.
$$
If $L\ge \ell$, then  this equation is satisfied if and only if $k\equiv0(p),$ so again there is a balance and the contributions to $\fS_q$ cancel. 

Lastly, if $L=\ell-1$, which is even since $\ell$ is odd, then note that $(n^2-4\gD)/p^L$ is a non-zero square mod $p$ if and only if $n^2-4\gD$ is a square mod $p^\ell$.
Whether or not this holds,  there are $\left({(n^2-4\gD)/p^L\over p}\right)+1$ solutions for $k(\mod p)$, and every lift of these to $\Z/p^{(\ell+1)/2}$ solves \eqref{eq:073001}. The number of these lifts is 
$$
p^{(\ell-1)/2}\left(\left({(n^2-4\gD)/p^L\over p}\right)+1\right),
$$ 
each contributing a Ramanujan value of $c_q=p^{\ell-1}(p-1)$.
And of course the complementary number of solutions to \eqref{eq:073002} that do not lift to \eqref{eq:073001} is 
$$
p^{(\ell-1)/2}\left(p-1-\left({(n^2-4\gD)/p^L\over p}\right)\right),
$$ 
each giving a Ramanujan value of $c_q=-p^{\ell-1}$. Recalling that there are $p^\ell$ values of $b$ in $\mattwos ab0{\bar a}$, the total contribution to $\fS_q$ is then:
$$
\fS_q
=
{p^{-(\ell-3)/2}
\over p^2-1}
\left({(n^2-4\gD)/p^L\over p}\right)
.
$$
This completes the proof.
\epf

We summarize this subsection as follows.

\begin{cor}\label{cor:summary}
Let $p$ be a good prime for $\G$, and let 
$$
\fS^{(p)}(n):=1+\fS_p(n)+\fS_{p^2}(n)+\cdots
$$
be the ``local factor'' at $p$. Then for all $p$,
$$
\fS^{(p)}(n)=1+\fS_p(n)+\fS_{p^2}(n)+O(p^{-3/2}).
$$
Moreover,
\begin{itemize}
\item
If $p\mid \gD$, then
$$
\fS^{(p)}(n)= 1 + O(p^{-1}).
$$
\item 
If $p\nmid \gD$ and $p\nmid n^2-4\gD$, then
$$
\fS^{(p)}(n)= 1 + \frac1p\left({n^2-4\gD\over p}\right)+O(p^{-2}).
$$
\item 
If $p\nmid \gD$ and $p\| n^2-4\gD$, then
$$
\fS^{(p)}(n)= 1 + O(p^{-3/2}).
$$
\item 
If $p\nmid \gD$ and $p^2\| n^2-4\gD$, then
$$
\fS^{(p)}(n)= 1 + O(p^{-1}).
$$
\item 
If $p\nmid \gD$ and $p^L\| n^2-4\gD$ with $L\ge3$, then
$$
\fS^{(p)}(n)= 1 + O(p^{-3/2}).
$$
\end{itemize}
\end{cor}

\subsection{Bad Primes}\

For bad primes, our strategy is as follows. Rather than evaluating $\fS_q(n)$ explicitly, we show the following ``density formula.''
\begin{lem}\label{lem:dens}
For any $\ell\ge0$ and any prime $p$ (good or bad), we have that
\be\label{eq:dens}
1+\fS_p(n)+\cdots +\fS_{p^\ell}(n)
=
p^\ell {\#\{\g\in\G(p^\ell)\bk\G:\ff_\g\equiv n(\mod p^\ell)\} \over [\G:\G(p^\ell)]}.
\ee
\end{lem}
This will tautologically capture the condition that $n$ is admissible.
And then, for $\ell$ large enough, we claim that $\fS_{p^\ell}(n)=0$, so these probabilities stabilize.
\pf[Proof of \lemref{lem:dens}]
This follows immediately from 
$$
\fS_{p^m}(n)
=
\frac1{|\G(p^\ell)\bk \G|}
\sum_{\g\in \G(p^\ell)\bk \G}
c_{p^m} (\ff_{\g}-n),
$$
for any $0\le m\le \ell$, together with the fact that
$$
1+c_p(x)+\dots+c_{p^\ell}(x)=\bo_{\{x\equiv0(\mod p^\ell)\}}p^\ell.
$$
\epf

Finally, we show that the densities stabilize.
\begin{lem}\label{lem:badPell}
Let $p$ be a bad prime, and let $k=k_p$ be the ``saturation level'' of $p$, as in \eqref{eq:kpDef}. Let $p^L\| B$ (and recall that $B\neq0$, and that $k> L$ by \eqref{eq:kpLB}). If $\ell>2k$, then $\fS_{p^\ell}(n)=0$.
\end{lem}
\pf
Decompose $\G(p^\ell)\bk\G$ into disjoint $\G_\infty=\{n_x=\mattwos1x01\}$ cosets; here
$x$ ranges in $\Z/p^\ell$ but is restricted (by saturation) to $x\equiv 0(p^k)$. 
We claim that the Ramanujan values on each coset exactly cancel. Note that $\ff_{\g}=Aa+Bb+Cc+Dd$ changes when $\g\mapsto n_x\g$ to 
$$
\ff_{n_x\g}=\ff_\g+(Ac+Bd)x.
$$
Since $c\equiv0(p^k)$ and $d\equiv1(\mod p^k)$, and $p^L\| B$ with $k>L$, we have that 
$$
p^L \| Ac+Bd.
$$
Now as $x$ ranges over $p^\ell$ subject to $x\equiv0(\mod p^k)$, since $\ell>2k>k+L$, the values of $\ff_{n_x\g}$ range in some non-constant arithmetic progression. 
The resulting Ramanujan values cancel exactly, as claimed.
\epf

So the high powers of bad primes have vanishing $\fS_q$. For the lower powers, we give the following trivial estimate on $\fS_q$.
\begin{lem}
For any prime $p$ (good or bad) and any $\ell\ge1$, we have:
\be\label{eq:anyPanyL}
|\fS_{p^\ell}(n)| \le p^\ell.
\ee
\end{lem}
\pf
The density formula \eqref{eq:dens} gives upper and lower bounds for its left-hand side of: $p^\ell$ and $0$, respectively. Replace $\ell$ by $\ell-1$ and subtract to get the claim.
\epf

\subsection{Short Sum of $\fS_q$}\label{sec:shortSum}\

Define 
$$
\fS(n) =
\sum_{q\ge 1}\fS_q(n)
.
$$

\begin{lem}
Assume that $\gD=0$. Then the series defining $\fS(n)$ is absolutely convergent,
$$
\sum_{q<Q_0} \fS_q(n) =\fS(n) + O_\vep (Q_0^{-1}n^\vep),
$$
as $Q_0\to\infty$,
 and satisfies, for $n$ admissible,
$$
\frac1{\log \log n} \ll 
\fS(n) 
\ll 1
.
$$
\end{lem}
\pf
By \lemref{lem:gD0p} and \lemref{lem:pell}, 
we have that
$$
\sum_{q\ge Q_0}|\fS_q(n)|\ll_\vep n^{\vep}Q_0^{-1}.
$$
For $n$ admissible,
$$
\fS(n) \asymp
\prod_{p\mid n}
\left(1-\frac1p\right)
,
$$
where we also used  \lemref{lem:dens} and  \lemref{lem:badPell}. The claim follows immediately.
\epf

To prepare for the case $\gD\neq0$, we need some preliminaries.

\begin{lem}\label{lem:Burgess}
Let $\chi$ be a Dirichlet character of conductor $M$ and fix $\Pi\in\Z$. Then
$$
\left|
\sum_{q\asymp H\atop {\text{squarefree}\atop (q,\Pi)=1}}
\chi(q)
\right|
\ll_\vep
H^{1/2}
M^{3/16}
(HM\Pi)^\vep,
$$
as $H\to\infty$.
\end{lem}
\pf
To capture both the squarefree and coprime conditions, we use M\"obius inversion. 
Using $\gz(s)/\gz(2s)=\sum_{n\text{ squarefree}}1/n^s$, we have that $\mu(q)^2=\sum_{m^2\mid q}\mu(m)$.
Similarly, $\sum_{d\mid x}\mu(d)=1$ if $x=1$ and $0$ otherwise. Therefore
\beann
\sum_{q\asymp H\atop {\text{squarefree}\atop (q,\Pi)=1}}
\chi(q)
 &= &
\sum_{q\asymp H\atop (q,\Pi)=1}
\chi(q)
\sum_{m^2\mid q}
\mu(m)
=
\sum_{m\ll H^{1/2} \atop (m,\Pi)=1}
\mu(m)
\chi(m)^2
\sum_{q\asymp H/m^2\atop (q,\Pi)=1}
\chi(q)
\\
&=&
\sum_{m\ll H^{1/2} \atop (m,\Pi)=1}
\mu(m)
\chi(m)^2
\sum_{q\asymp H/m^2}
\chi(q)
\sum_{d\mid q\atop d\mid \Pi}
\mu(d)
\\
&=&
\sum_{m\ll H^{1/2} \atop (m,\Pi)=1}
\mu(m)
\chi(m)^2
\sum_{d\ll H/m^2\atop d\mid\Pi}
\mu(d)
\chi(d)
\sum_{q\asymp H/(m^2d)}
\chi(q)
.
\eeann
Applying Burgess \cite{Burgess1962} to the last sum, we have that 
\beann
\left|
\sum_{q\asymp H\atop {\text{squarefree}\atop (q,\Pi)=1}}
\chi(q)
\right|
&\ll_\vep &
\sum_{m\ll H^{1/2} }
\sum_{d\ll H/m^2\atop d\mid\Pi}
{H^{1/2}\over md^{1/2}}
M^{3/16+\vep}
\ll_\vep
\Pi^\vep
{H^{1/2+\vep}}
M^{3/16+\vep},
\eeann
as claimed. (Slightly better estimates are available today but not needed here.)
\epf

Going forward, we let
$\fB_1$ be the (finitely many) primes which are ``bad'' (for $\G$),
$\fB_2$ be the primes not in $\fB_1$ which  divide $\gD$, and
$\fB_3=\fB_3(n)$ be the primes not in $\fB_1$ or $\fB_2$ which divide $n^2-4\gD$.
Let
 $\fB=\fB(n)=\sqcup_{j}\fB_j$, and set
$$
\Pi = \Pi(\fB) := \prod_{p\in\fB}p.
$$
For all the other primes, \lemref{lem:gDneq0} gives that  
$$
\fS_p(n)=  \frac1p\left({n^2-4\gD\over p}\right)+E_p,
$$
where
$$
E_p = E_p(n) :=  \frac{p+\left({n^2-4\gD\over p}\right)}{p(p^2-1)}  \ll \frac1{p^2}.
$$
We extend $E_p$ to a multiplicative function $E_q$ supported on square-free $q$.

Note that we now do not have absolute convergence, and must be much more careful. 

We break the tail $\sum_{q\ge Q_0}$ of $\fS(n)$ into dyadic regions $\sum_{q\asymp H}$, with $Q_0\le H\to\infty$.

\begin{lem}\label{lem:Hlarge}
Assume that $\gD\neq0$. Then as $H\to\infty$,
\be\label{eq:Hlarge}
\left|
\sum_{q\asymp H}\fS_q(n) 
\right|
\ll_\vep
(nH)^\vep
n^{3/8}
{H^{-1/2}}
,
\ee
for any $\vep>0$.
\end{lem}
\pf
From
\lemref{lem:gDne0Pell},
we have that $\fS_{q}$ vanishes if $(q,\Pi)=1$ and $q$ is not square-free. For $q$ square-free and coprime to $\Pi$, we have then that
\be\label{eq:fSqEq}
\fS_q(n)= 
\prod_{p\mid q}
\left(
\frac1p\left({n^2-4\gD\over p}\right)+E_p
\right)
=
\sum_{ab= q}
  \frac1a\left({n^2-4\gD\over a}\right)E_b(n).
\ee

Write any $q$ as
$$
q=q_\fB \cdot q_1,
$$
where
$$
q_\fB=
\prod_{p^\ell \| q \atop p\mid \Pi}p^\ell,
\text{ and }
q_1 =
\prod_{p^\ell \| q \atop (p,\Pi)=1}p^\ell
$$
From multiplicativity, we have that $\fS_q=\fS_{q_\fB}\cdot \fS_{q_1}$.

Then we have
\beann
\sum_{q\asymp H}\fS_q(n) 
&=&
\sum_{q_\fB \ll H \atop p\mid q_\fB \Longrightarrow p\mid\Pi}
\fS_{q_\fB}(n)
\sum_{q_1\asymp H/q_\fB 
\atop 
{
\text{square-free}\atop
( q_1, \Pi)=1}}
\fS_{q_1}(n)
\\
&=&
\sum_{q_\fB \ll H \atop p\mid q_\fB \Longrightarrow p\mid\Pi}
\fS_{q_\fB}(n)
\sum_{q_1\asymp H/q_\fB 
\atop 
{
\text{square-free}\atop
( q_1, \Pi)=1}}
\sum_{ab=q_1}
 \frac1{a}\left({n^2-4\gD\over a}\right)
E_b
\\
&=&
\sum_{q_\fB \ll H \atop p\mid q_\fB \Longrightarrow p\mid\Pi}
\fS_{q_\fB}(n)
\sum_{b\ll H/q_\fB\atop \text{square-free},\ (b,\Pi)=1}
E_b
\sum_{a\asymp H/(q_\fB b)
\atop 
{
\text{square-free}\atop
( a, \Pi)=1}}
 \frac1{a}\left({n^2-4\gD\over a}\right)
\\
&\ll_\vep&
(nH)^\vep
n^{3/8}
{H^{-1/2}}
\sum_{q_\fB \ll H \atop p\mid q_\fB \Longrightarrow p\mid\Pi}
\fS_{q_\fB}(n)
q_\fB^{1/2}
,
\eeann
where we used \lemref{lem:Burgess}, partial summation, and $E_b\ll b^\vep/b^2$.

To deal with the remaining $q_\fB$ sum, we decompose 
$$
\fS_{q_\fB}(n)q_\fB^{1/2}=
\fS_{q_{\fB_1}}(n)q_{\fB_1}^{1/2}\cdot
\fS_{q_{\fB_2}}(n)q_{\fB_2}^{1/2}\cdot
\fS_{q_{\fB_3}}(n)q_{\fB_3}^{1/2},
$$
corresponding to $\fB=\fB_1\sqcup\fB_2\sqcup\fB_3$.

Since $\fB_1$ is a finite set of primes which are bad for $\G$, and only finitely many powers of such primes have non-vanishing $\fS_{q_{\fB_1}}$ by \lemref{lem:badPell}, the total contribution from $\fB_1$ is  bounded by a constant depending only on $\G$ and the linear form $\sL$, that is, on $A,B,C,D$.

Recall that $\fB_2$ consists of the good primes dividing $\gD$; \lemref{lem:pell} removes any non-square-free $q_{\fB_2}$ contributions, and \lemref{lem:gD0p} otherwise gives $\fS_{q_{\fB_2}}(n)\ll 1/q_{\fB_2}$. So again this contribution is bounded.

Finally, for $\fB_3$, we use \eqref{eq:gDne0Pell} to offset the factor of $q_{\fB_3}^{1/2}$, and \eqref{eq:gDne0Pell1} to kill the contribution from any powers $\ell$ of $p^\ell$ in $q_{\fB_3}$ unless  $\ell\le L+1\le2L$, where $p^L\| n^2-4\gD$. Therefore the only $q_{\fB_3}$ contributing to the sum are divisors of $(n^2-4\gD)^2$, and the number of such is $\ll n^\vep$.
This gives the claim.
\epf

\lemref{lem:Hlarge} is sufficient to show that $\fS(n)$ converges (conditionally, not absolutely), but does not allow us a good enough error estimate for the very short sum $\sum_{q<Q_0}\fS_q(n)$, since $H$ needs to be at least $n^{3/4+}$ for \eqref{eq:Hlarge} to decay. If we replaced our use of Burgess with GRH, we could get good estimates with $H$ as small as $Q_0$, which is a tiny power of $N$. Unconditionally, we can only do this on average, as follows.

\begin{thm}\label{thm:L2err}
As $H\to\infty$, we have
$$
\sum_{n\ll N}
\left|
\sum_{q\asymp H}
\fS_q(n)
\right|^2
\ \ll_\vep \
H^\vep
\left(
H
+
{N \over H^{1/2}}
\right)
,
$$
for any $\vep>0$.
\end{thm}
\pf
As before, let $\fB_1$ be the ``bad'' primes for $\G$, and $\fB_2$ be the primes not in $\fB_1$ which divide $\gD$. 
Since $\fB_3$ depends on $n$, we now have to handle it separately. We now write $\fB=\fB_1\sqcup\fB_2$ and $\Pi:=\prod_{p\in\fB}p$ as before, and decompose
$$
q=q_\fB\cdot q_1,
$$
with $(q_1,\Pi)=1$.
Furthermore, we will split off the square-full part of $q_1$, writing $q_1=q_2\cdot q_3$, where
$$
q_2 :=\prod_{p\|q_1}p,\ \ \ q_3:=q_1/q_2=\prod_{p^\ell \| q_1\atop \ell\ge2}p^\ell.
$$

With this decomposition, we open the square and reverse orders:
\beann
&&
\hskip-.5in
\sum_{n\ll N}
\left|
\sum_{q\asymp H}
\fS_q(n)
\right|^2
\\
&=&
\sum_{n\ll N}
\left|
\sum_{q_\fB\ll H \atop p|q_\fB\Longrightarrow p\in\fB}
\fS_{q_\fB}(n)
\sum_{q_3\ll H/q_\fB \atop (q_3,\Pi)=1,\ \text{square-full}}
\fS_{q_3}(n)
\sum_{q_2\asymp H/(q_\fB q_3)\atop (q_2,\Pi)=1=(q_2,q_3),\ \text{square-free}}
\fS_{q_2}(n)
\right|^2
\\
&=&
\sum_{q_\fB, q_3, q_2, q_\fB', q_3', q_2'}
\sum_{n\ll N}
\fS_{q_\fB}(n)
\fS_{q_\fB'}(n)
\fS_{q_3}(n)
\fS_{q_3'}(n)
\fS_{q_2}(n)
\fS_{q_2'}(n)
\eeann

Here, instead of using the decomposition  \eqref{eq:fSqEq}, we return to  \lemref{lem:gDneq0} and write
$$
\fS_p(n) = 
\left({n^2-4\gD\over p}\right){p \over p^2-1}
+
E_p,
$$
where
$$
E_p = {1\over p^2-1}.
$$
The crucial fact for our purposes here is that $E_p$ is now independent of $n$. Extending $E_p$ to a multiplicative function on square-frees gives
\be\label{eq:Eqbnd}
E_q \ll {1\over q^2},
\ee
where we used that $\prod_{p\mid q}(1-1/p^2)\asymp1$.
Then we can write, for any $q$ square-free and coprime to $\Pi$, that
$$
\fS_q(n)=\prod_{p\mid q}\left(\left({n^2-4\gD\over p}\right){p \over p^2-1}
+
E_p\right) 
=
\sum_{ab=q}
\psi(a)
\left({n^2-4\gD\over a}\right)
E(b)
,
$$
where $\psi$ is a multiplicative function supported on square-free numbers taking the value $\psi(p)=p/(p^2-1)$ on primes. In particular,
\be\label{eq:PsiBnd}
\psi(q) \ll \frac1q.
\ee

For $q_2$ and $q_2'$, we insert this expression to get:
\bea
\sum_{n\ll N}
\left|
\sum_{q\asymp H}
\fS_q(n)
\right|^2
&=&
\sum_{q_\fB, q_3, q_2, q_\fB', q_3', q_2'}
\sum_{n\ll N}
\fS_{q_\fB}(n)
\fS_{q_\fB'}(n)
\fS_{q_3}(n)
\fS_{q_3'}(n)
\\
\nonumber
&&
\times
\sum_{ab= q_2}
\psi(a)\left({n^2-4\gD\over a}\right)E_b
\sum_{a'b'= q_2'}
\psi(a')\left({n^2-4\gD\over a'}\right)E_{b'}
\\
\nonumber
&\le&
\sum_{q_\fB, q_3, q_2, q_\fB', q_3', q_2'}
\sum_{ab= q_2}
\sum_{a'b'= q_2'}
\psi(a)
\psi(a')
|E_b
E_{b'}|
  \\
  \nonumber
  &&
  \times
\sum_{n_0\mod \tilde q}
\left|\fS_{q_\fB}(n_0)
\fS_{q_\fB'}(n_0)
\fS_{q_3}(n_0)
\fS_{q_3'}(n_0)
\right|
\\
\label{eq:081204}
&&
\times
\left|
    \sum_{n\ll N\atop n\equiv n_0(\mod \tilde q)}
\left({n^2-4\gD\over a}\right)
\left({n^2-4\gD\over a'}\right)
\right|
,
\eea
where we have decomposed $n$ into progressions mod $\tilde q$, where 
$$
\tilde q := [q_\fB,q_\fB',q_3,q_3'].
$$
While it is clear at, say, $a$ and $q_3$ are coprime (by construction, since $a\mid q_2$), we actually also have that there is no contribution unless $a$ and $q_3'$ are coprime. Indeed, if $p\mid a$ and $p\mid q_3'$, then either $n^2-4\gD\equiv0(p)$, in which case  $\left({n^2-4\gD\over a}\right)$ vanishes, or else $\fS_{q_3'}(n)$ vanishes from \eqref{eq:gDne0Pell1} and the square-full-ness of $q_3'$.
Therefore, we may restrict the summations to 
$$
(a,\tilde q)=(a',\tilde q)=1.
$$

We first analyze the last $n$ sum. Let $\tilde a := a a'/(a,a')^2$, so that 
$$
\left({n^2-4\gD\over a}\right)
\left({n^2-4\gD\over a'}\right)
=
\left({n^2-4\gD\over \tilde a}\right),
$$
for $n$ such that $n^2-4\gD$ is coprime to $(a,a')$. (Otherwise the characters vanish.)
Breaking the $n$ sum further into residue classes mod $\tilde a$ gives:
\beann
    \sum_{n\ll N\atop {n\equiv n_0(\mod \tilde q)\atop (n^2-4\gD,a,a')=1}}
\left({n^2-4\gD\over \tilde a}\right)
&=&
\sum_{m(\mod \tilde a)}
\left({m^2-4\gD\over \tilde a}\right)
\left[
    \sum_{n\ll N\atop {n\equiv n_0(\mod \tilde q),\ n\equiv m(\mod \tilde a) \atop (n^2-4\gD,a,a')=1}}
1
\right]
\eeann
We wish to get square-root cancellation from Weil in the $m$ summation, but
the $n$ sum may be incomplete, which will give too large an error  in terms of $m$. So we separate the roles of $n$ and $m$ by completing the sum.
\bea\nonumber
    \sum_{n\ll N\atop {n\equiv n_0(\mod \tilde q)\atop (n^2-4\gD,a,a')=1}}
\left({n^2-4\gD\over \tilde a}\right)
&=&
\frac1{\tilde a}
\sum_{k(\tilde a)}
\left[
\sum_{m(\mod \tilde a)}
\left({m^2-4\gD\over \tilde a}\right)
e_{\tilde a}(-km)
\right]
\\
\label{eq:081203}
&&\times
\left[
    \sum_{n\ll N\atop {n\equiv n_0(\mod \tilde q)\atop (n^2-4\gD,a,a')=1}}
e_{\tilde a}(kn)
\right].
\eea
Now the $m$ sum is free and is bounded by $\tilde a^{1/2+\vep}$ by Weil.
We now deal with the last $n$ sum. We remove the gcd condition via M\"obius inversion.
\bea\nonumber
    \sum_{n\ll N\atop {n\equiv n_0(\mod \tilde q)\atop (n^2-4\gD,a,a')=1}}
e_{\tilde a}(kn)
&=&
\sum_{d\mid (a,a') }
\mu(d)
    \sum_{n\ll N\atop {n\equiv n_0(\mod \tilde q)\atop n^2-4\gD \equiv 0(\mod d)}}
e_{\tilde a}(kn)
\\
\label{eq:081202}
&=&
\sum_{d\mid (a,a') }
\mu(d)
\sum_{m_0(\mod d)\atop m_0^2\equiv 4\gD(\mod d)}
\left[
    \sum_{n\ll N\atop {n\equiv n_0(\mod \tilde q)\atop n \equiv m_0(\mod d)}}
e_{\tilde a}(kn)
\right]
,
\eea
where we decomposed the $n$ sum further into residue classes $m_0$ mod $d$. Note that $d\mid a$ is square-free, and for each $p\mid d$, there are at most two solutions to $m_0^2\equiv 4\gD(p)$, so the number of $m_0$ is at most $d^\vep$. The last bracketed sum restricts $n$ to a residue class $x$, say, mod $\tilde q d$ (since $(d,\tilde q)=(a,\tilde q)=1$). Changing $n\mapsto x+n\tilde q d$, the bracketed term is a geometric series, giving:
\bea
\nonumber
\left|
    \sum_{n\ll N\atop {n\equiv n_0(\mod \tilde q)\atop n \equiv m_0(\mod d)}}
e_{\tilde a}(kn)
\right|
&=&
\left|
e_{\tilde a}(kx)
    \sum_{n\ll N/(\tilde q d)}
e_{\tilde a}(k\tilde q dn)
\right|
\\
\label{eq:081201}
&\ll&
\min\left(
{N\over \tilde q }+1,
\frac1{\| {k\tilde q d\over \tilde a } \|}
\right)
,
\eea
where $\|\cdot\|$ is the distance to the nearest integer.
Inserting \eqref{eq:081201} into \eqref{eq:081202} and into \eqref{eq:081203} gives
\beann
\left|    \sum_{n\ll N\atop {n\equiv n_0(\mod \tilde q)\atop (n^2-4\gD,a,a')=1}}
\left({n^2-4\gD\over \tilde a}\right)
\right|
&\ll_\vep&
\tilde a^{1/2+\vep}
\sum_{d\mid (a,a') }
d^\vep
\frac1{\tilde a}
\sum_{k(\tilde a)}
\min\left(
{N\over \tilde q }+1,
\frac1{\| {k\tilde q d\over \tilde a } \|}
\right)
.
\eeann
Since $a$ and $a'$ are square-free, $d$ is coprime to $\tilde a$, and hence $(\tilde q d,\tilde a)=1$. So the $k$ sum is invariant under $k\mapsto k\overline{\tilde q d}$. This finally gives
\beann
\left|    \sum_{n\ll N\atop {n\equiv n_0(\mod \tilde q)\atop (n^2-4\gD,a,a')=1}}
\left({n^2-4\gD\over \tilde a}\right)
\right|
&\ll_\vep&
\tilde a^{\vep}
\left(
{N\over \tilde q \tilde a^{1/2}}
+
\tilde a^{1/2}
\right)
.
\eeann

Returning to \eqref{eq:081204}, we get that
\beann
\sum_{n\ll N}
\left|
\sum_{q\asymp H}
\fS_q(n)
\right|^2
&\ll&
H^\vep
\sum_{q_\fB, q_3, q_2, q_\fB', q_3', q_2'}
\sum_{ab= q_2}
\sum_{a'b'= q_2'}
\frac1{aa'b^2b'^2}
  \\
  &&
  \times
\sum_{n_0\mod \tilde q}
\left|\fS_{q_\fB}(n_0)
\fS_{q_\fB'}(n_0)
\fS_{q_3}(n_0)
\fS_{q_3'}(n_0)
\right|
\\
&&
\times
\left(
{N (a,a')\over \tilde q (aa')^{1/2}}
+
{( aa')^{1/2}\over (a,a')}
\right)
,
\eeann
where we used \eqref{eq:PsiBnd} and  \eqref{eq:Eqbnd}.

Next we analyze the contributions from $q_\fB,q_\fB'$. 
Combining \lemref{lem:pell} and \lemref{lem:gD0p} (for $q_{\fB_2}$) with  \lemref{lem:badPell} and \eqref{eq:anyPanyL} (for $q_{\fB_1}$), we see that in fact there is no contribution unless $q_\fB,q_\fB' \ll1$, and in this case the constribution is $\fS_{q_\fB}\fS_{q_\fB'}\ll1$. Therefore $\tilde q\asymp [q_3,q_3']$.

Recall from \eqref{eq:gDne0Pell} that $\fS_{q_3}(n_0)\ll q_3^{-1/2}$.
Finally, we analyze the number of $n_0(\mod \tilde q)$ for which $\fS_{q_3}\fS_{q_3'}$ is non-vanishing. 
Suppose that $p^m\|[q_3,q_3']$. Then since $\gD\not\equiv0(p)$,
 \eqref{eq:gDne0Pell1} shows that, if $\fS_{p^m}(n_0)\neq0$, then $n_0^2-4\gD\equiv 0(\mod p^{m-1})$.  The number of such $n_0(\mod p^m)$ is at most $2p\ll p^{m/2},$ since $m\ge2$.
 So the number of $n_0(\mod \tilde q)$ which contribute is $\ll_\vep \tilde q ^{1/2+\vep}$.

Putting everything together gives
\beann
\sum_{n\ll N}
\left|
\sum_{q\asymp H}
\fS_q(n)
\right|^2
&\ll_\vep&
H^\vep
\sum_{q_3, q_3'}
([q_3,q_3'])^{1/2}
q_3^{-1/2}
q_3'^{-1/2}
\\
&&
\times
\sum_{q_2, q_2'}
\sum_{a\mid q_2}
\sum_{a'\mid q_2'}
\frac{aa'}{q_2^2q_2'^2 }
\left(
{N (a,a')\over [q_3,q_3'] (aa')^{1/2}}
+
{( aa')^{1/2}\over (a,a')}
\right)
.
\eeann
Let $t:=(a,a')$, which is a divisor of 
$(q_2,q_2')$, and let $a_1:=a/t$ and $a_1':=a'/t$. Then
\beann
\sum_{n\ll N}
\left|
\sum_{q\asymp H}
\fS_q(n)
\right|^2
&\ll_\vep&
H^\vep
\sum_{q_3, q_3'}
([q_3,q_3'])^{1/2}
q_3^{-1/2}
q_3'^{-1/2}
\\
&&
\times
\sum_{q_2, q_2'}
\frac1{q_2^2q_2'^2 }
\sum_{t\mid (q_2,q_2')}
t^2
\sum_{a_1\mid \frac{q_2}t}
\sum_{a_1'\mid \frac{q_2'}t}
\left(
{N(a_1a'_1)^{1/2} \over [q_3,q_3'] }
+
{( a_1a'_1)^{3/2}}
\right)
\\
&\ll&
H^\vep
\sum_{q_3\ll H \atop {(q_3,\Pi)=1\atop \text{square-full}}}
\sum_{q_3'\ll H \atop {(q_3',\Pi)=1\atop \text{square-full}}}
([q_3,q_3'])^{1/2}
q_3^{-1/2}
q_3'^{-1/2}
\\
&&
\times
\left(
{H\over (q_3q_3')^{1/2}}
+
{N \over  [q_3,q_3'] }
\sum_{q_2\asymp H/q_3\atop {(q_2,\Pi)=1=(q_2,q_3)\atop \text{square-free}}}
\sum_{q_2'\asymp H/q_3'\atop {(q_2',\Pi)=1=(q_2',q_3')\atop \text{square-free}}}
{(q_2,q_2') \over  (q_2q_2')^{3/2} }
\right)
.
\eeann
Next we need some cancellation from $(q_2,q_2')$. Let $d=(q_2,q_2')$ which is a divisor of $q_2$ such that $q_2'\equiv0(d)$.
\beann
\sum_{n\ll N}
\left|
\sum_{q\asymp H}
\fS_q(n)
\right|^2
&\ll_\vep&
H^\vep
\sum_{q_3\ll H \atop {(q_3,\Pi)=1\atop \text{square-full}}}
\sum_{q_3'\ll H \atop {(q_3',\Pi)=1\atop \text{square-full}}}
([q_3,q_3'])^{1/2}
q_3^{-1/2}
q_3'^{-1/2}
\\
&&
\times
\left(
{H\over (q_3q_3')^{1/2}}
+
{N \over  [q_3,q_3'] }
\sum_{q_2\asymp H/q_3\atop {(q_2,\Pi)=1=(q_2,q_3)\atop \text{square-free}}}
\sum_{d\mid q_2}
\sum_{q_2'\asymp H/(q_3'd)\atop {(q_2',\Pi)=1=(q_2',q_3')\atop \text{square-free}}}
{d \over  (q_2q_2'd)^{3/2} }
\right)
\\
&\ll_\vep&
H^\vep
\sum_{q_3\ll H \atop {(q_3,\Pi)=1\atop \text{square-full}}}
\sum_{q_3'\ll H \atop {(q_3',\Pi)=1\atop \text{square-full}}}
\left(
H
{1
\over
(q_3
q_3')^{1/2}
(q_3,q_3')^{1/2}
}
+
{N \over H}
{(q_3,q_3')^{1/2}
\over   (q_3q_3')^{1/2} }
\right)
.
\eeann
Finally, we bound $(q_3,q_3')\ll H$ in the numerator and $(q_3,q_3') \ge1$ in the denominator. It remains to estimate a sum of the form
$$
\sum_{q_3\ll H\atop \text{square-full}}
\frac1{q_3^{1/2}}.
$$
Since $q_3$ is square-full, any such $q_3$ can be written as $q_3=k^2\ell$ where $\ell\mid k$. Then
$$
\sum_{q_3\ll H\atop \text{square-full}}
\frac1{q_3^{1/2}}
\ll
\sum_{k^2\ll H}
\sum_{\ell\mid k}
\frac1{k\ell^{1/2}}
\ll_\vep
\sum_{k^2\ll H}
k^\vep
\frac1{k}
\ll
H^\vep.
$$
The claim follows immediately.
\epf

\thmref{thm:L2err} allows us to show, for almost all $n$ (with power savings error), that the very short sum
$
\sum_{q\le Q_0}\fS_{q}(n_0) 
$
(with $Q_0= N^{\ga_0}$, $\ga_0>0$ small) is a good approximation (also with power savings error) for $\fS(n)$.

\begin{thm}\label{thm:sENfSq}
For any $\eta>0$ with $\eta<\frac16\ga_0$, there is  a set $\sE$ of ``exceptional'' $n$ of cardinality
$$
\sE\cap[1,N]\ll N^{1-\eta}
$$
such that, for all $n\asymp N$, $n\not\in\sE$, 
$$
\sum_{q\le Q_0}\fS_q(n) \ = \ \fS(n) + O(N^{-\eta}).
$$
\end{thm}
\pf
Recall that  the series $\fS(n)$ does converge (conditionally) by \lemref{lem:Hlarge}, but the error there is insufficient to approximate it to the required error in all ranges of $H$.

For $\eta>0$ fixed, let 
$$
\sE(N) := \left\{ n\in[1,N]:
\left|
\fS(n) - \sum_{q\le Q_0}\fS_q(n)\right|
\ge N^{-\eta}
\right\}
.
$$

We estimate
\beann
\#\sE(N)
&=&
\sum_{n\ll N\atop \left|
\fS(n) - \sum_{q\le Q_0}\fS_q(n)\right|
\ge N^{-\eta}
}
1
\  \ \le \ \
N^{2\eta}
\sum_{n\ll N\atop 
}
\left|\fS(n) - \sum_{q\le Q_0}\fS_q(n)\right|^2
\\
&\ll&
N^{2\eta}
\sum_{n\ll N\atop 
}
\left(
\sum_{Q_0< H<N^{4/5}\atop\text{dyadic}}
\left|
\sum_{q\asymp H}
\fS_q(n)
\right|^2
+
\sum_{H\ge N^{4/5}\atop\text{dyadic}}
\left|
\sum_{q\asymp H}
\fS_q(n)
\right|^2
\right).
\eeann
We apply \thmref{thm:L2err} in the first term, and \lemref{lem:Hlarge} (individually) in the second term.
\beann
\#\sE(N)
&\ll_\vep&
N^{2\eta+\vep}
\left(
N^{4/5}+{N\over Q_0^{1/2}}
+
\sum_{H\ge N^{4/5}\atop\text{dyadic}}
\sum_{n\ll N\atop 
}
n^{3/4}
H^{-1}
\right)
\\
&\ll&
N^{2\eta+\vep}
\left(
N^{4/5}+{N\over Q_0^{1/2}}
+
N^{7/4}
N^{-4/5}
\right)
.
\eeann
Since $7/4-4/5=19/20<1$, we have a power savings as long as $2\eta<\foh \ga_0$, where $Q_0=N^{\ga_0}$. As long as $\eta<\frac16\ga_0$, we are guaranteed to have $\sE(N)\ll N^{1-\eta}$. This completes the proof.
\epf

\begin{thm}
For all admissible $n\asymp N$, and all $\vep>0$,
$$
\fS(n)\gg_\vep n^{-\vep} L(1,\chi_n),
$$
where
\be\label{eq:L1def}
L(1,\chi_n):=\prod_p\left(1-\frac1p\left({n^2-4\gD\over p}\right)\right)^{-1}.
\ee
The implied constant is effective.
\end{thm}
\pf
By the multiplicativity of $\fS_q$, we have that 
$$
\fS(n) =
\prod_p
\left(
1+\fS_p(n)+\fS_{p^2}(n)+\cdots
\right).
$$
(Since the series only converges conditionally, we argue by considering the functions $s\mapsto \sum_{q\in\N}\fS_q(n)q^{-s}$ and $s\mapsto\prod_p(1+\fS_p(n)p^{-s}+\fS_{p^2}(n)p^{-2s}+\cdots)$; for $\Re(s)>0$, both converge absolutely and coincide, and hence their limiting values as $s\to0^+$ do too.)

For $p\in\fB_1$ a ``bad'' prime for $\G$, this is a finite sum (\lemref{lem:badPell}) which is non-vanishing only if $n$ is admissible by \lemref{lem:dens}.
For the other primes $p$, the Euler factor is $(1+\fS_p(n)+\fS_{p^2}(n))(1+O(p^{-3/2}))$.
Recall that $\fB_2$ contains the (finite list of) primes $p\mid \gD$.
By \lemref{lem:gD0p} and \lemref{lem:pell}, we have
$$
\prod_{p\in\fB_2}(1+\fS_p(n))\gg{1\over \log\log n}.
$$

For all other primes we apply \corref{cor:summary}. If $p\nmid\gD$ and $p\nmid n^2-4\gD$, we have
$$
1+\fS_p(n)+\fS_{p^2}(n) = 1+\frac1p\left({n^2-4\gD\over p}\right)+O(p^{-2}),
$$
while if $p\nmid\gD$ but $p\mid n^2-4\gD$,
$$
1+\fS_p(n)+\fS_{p^2}(n) = 1+O(p^{-1}).
$$
The product of the latter (finite set of primes) is $\gg_\vep n^{-\vep}$.
\epf

By Siegel's theorem, $L(1,\chi_n)\gg_\vep n^{-\vep}$ with an {\it ineffective} implied constant. But since we anyway only prove our result on average over $n$, we want to make this constant effective. 

\begin{thm}\label{thm:fSnSize}
There is an
 exceptional set $\sE$ with the following property.
For all admissible $n\asymp N$ outside of $\sE$,
and any $\vep>0$, we have
$$
\fS(n)\gg_\vep n^{-\vep}.
$$
Moreover
\be\label{eq:sEsize}
 \#\sE\ll_\vep N^{\vep}.
\ee
The implied constants are all effective. (But the exact determination of the exceptional set $\sE$ is ineffective!)
\end{thm}

\pf
Consider the characters $\chi_n=\left({n^2-4\gD\over \cdot}\right)$ appearing in \eqref{eq:L1def}. These need not be primitive, and are induced from characters $\left({q_n\over \cdot}\right)$, where $q_n:=\sqf(n^2-4\gD)\ll N^2$ is the square-free part of $n^2-4\gD$; that is,
\be\label{eq:n2qm2}
n^2-4\gD=q_n m^2,
\ee
for some integer $m$. Group admissible $n\asymp N$ according to the values of $q_n$; that is, for a given square-free $q\ll N^2$, let
$$
\cN_q := \{n\asymp N : \sqf(n^2-4\gD)=q\}.
$$
If $(n,m)$ is a solution to $n^2-qm^2=4\gD$, then the ideal $(n+\sqrt q m)$ in $\Z[\sqrt q]$ has norm $|4\gD|$. The prime ideals $\fp$ dividing $(n+\sqrt q m)$ and their multiplicities are bounded in terms of those of the rational primes dividing $4\gD$ (which is fixed). Therefore there are $\ll1$ inequivalent solutions to \eqref{eq:n2qm2}, and equivalent solutions grow exponentially in terms of the units in $\Z[\sqrt q]$. Therefore
\be\label{eq:maxQ}
\max_{q\ll N^2}\#\cN_q \ll_\vep N^\vep,
\ee
for any $\vep>0$ with  absolute implied constants.

By Landau's theorem (see, e.g., \cite[Theorem 5.28]{IwaniecKowalski}), there is an absolute constant $A>0$, such that for all distinct primitive real characters $\chi$, $\chi'$ of conductors $q$, $q'$ (resp.), with $L$-functions $L(s,\chi)$, $L(s,\chi')$ having largest real zeros $\gb$, $\gb'$ (resp.), we have:
$$
\max(\gb,\gb')\le 1-{A\over \log (qq')}.
$$
Therefore, there is at most a single exceptional $\fq\ll N^2$ such that, for all other square-free $q\ll N^2$ and their corresponding largest real zeros $\gb$ (if any such exist), we have
$$
\gb \le 1-{A'\over \log N},
$$
where $A'>0$ is another absolute constant.

We then define the exceptional set $\sE:=\cN_\fq$, so that the bound \eqref{eq:sEsize} is confirmed by \eqref{eq:maxQ}, again with absolute constants. (Though we cannot effectively determine the elements of $\sE$, we can effectively control their cardinality.)

Then we use   standard arguments (see, e.g., \cite{Goldfeld1974}),  and take into account the imprimitive factors, to show that $L(1,\chi_n)\gg_\vep N^{-\vep}$ with absolute implied constants, for all $n\not\in\sE$. This gives the claim.
\epf

\newpage

\section{Minor Arc Technical Estimates}\label{sec:minTech}

We collect here various lemmata needed in the analysis of the minor arcs.
We begin by defining the exponential sum
\be\label{eq:cSqDef}
\cS_q(r,k,\ell;\g)  \ := \
\frac1{q^2}
\sum_{x(q)}
\sum_{y(q)}
e_q(r \ff_\g(x,y)+kx+\ell y)
.
\ee

\begin{lem}\label{lem:cSq}
Assume that $(r,q)=1$. 
Write $q_1:=(BcP^2,q)$, $q=q_1q_2$, and $BcP^2=q_1 E$, with $E\bar E\equiv 1(\mod q_2)$. Then
\beann
\cS_q(r,k,\ell;\g) 
&=&
{(BcP^2,q)\over q}
e_q(r
(A a
+B b
+C c
+D d
))
\bo_{\{
-\ell \equiv Pr
(B  a
+D  c)
(\mod q_1)
\atop
-k
\equiv
Pr
(A c
+B  d)
(\mod q_1)
\}}
\\
&&
\times
e_{qq_1}\left(-
\bar r
\bar E
(
Pr(B  a
+D  c)
+\ell )
(Pr
(A c
+B  d
)
+k)
\right)
\eeann
Note that the last exponential term is well-defined by the congruence conditions on $\ell$ and $k$, and  independent of the lifts of $\bar r, \bar E$ to $\Z/(qq_1)$.
\end{lem}

\pf
Write $\g=\mattwos abcd$ and insert \eqref{eq:ffIs}:
\beann
\cS_q(r,k,\ell;\g)  &=&
\frac1{q^2}
e_q(r
(A a
+B b
+C c
+D d
))
\sum_{x(q)}
e_q(r
(A c
+B  d)
Px
+kx)
\\
&&\times
\sum_{y(q)}
e_q(y
[r
((B  a
+D  c)
P
+B   c
P^2
x)
+\ell
]).
\eeann
The $y$ sum vanishes unless
\be\label{eq:xeqn}
B   c
P^2
x
\equiv  -\ell\bar r -(B  a
+D  c)
P
 \pmod q
,
\ee
in which case the sum contributes $q$.

Let $q_1=\gcd(BcP^2,q)$ and write $q=q_1q_2$ and $BcP^2=q_1 E$. Then \eqref{eq:xeqn} has a solution only if the right hand side is congruent to zero mod $q_1$.
If this is the case, then $x$ is determined mod $q_2$,
$$
x
\equiv
x_0:=-
\overline{E}
\frac{  \ell\bar r +(B  a
+D  c)
P}{q_1}
 \pmod {q_2}
 .
$$
So $x\equiv x_0+q_2 x'$, where $x'\in\Z/q_1$.
Thus  we have:
\beann
\cS_q(r,k,\ell;\g)  &=&
e_q(r
(A a
+B b
+C c
+D d
))
\bo_{\{
-\ell \equiv Pr
(B  a
+D  c)
(\mod (BcP^2,q))\}}
\\
&&\times
\frac1{q}
e_q(x_0(r
(A c
+B  d)
P
+k))
\sum_{x'(q_1)}
e_{q_1}(x'(r
(A c
+B  d)
P
+k))
.
\eeann
The $x'$ sum vanishes unless 
$$
-k
\equiv
Pr
(A c
+B  d)
(\mod q_1),
$$
in which case it contributes $q_1$.
\epf

Next we need cancellation over the $r$ sum on the product of two such. A preliminary calculation is the following.

\begin{lem}\label{lem:cSq2}
Assume that $(r,q)=1$ as before, and also use $\mattwos{a'}{b'}{c'}{d'}=\g'$.
Write $q_1:=(BcP^2,q)$, $q=q_1q_2$, and $BcP^2=q_1 E$, with $E\bar E\equiv 1(\mod q_2)$.
Similarly, set $q_1':=(Bc'P^2,q)$, $q=q_1'q_2'$, and $Bc'P^2=q_1' E'$, with $E'\bar E'\equiv 1(\mod q_2')$.
 Then
\bea\nonumber
\sum_{r(q)}'
\cS_q(r,k,\ell;\g) 
\overline{
\cS_q(r,k',\ell';\g') }
&=&
{(BcP^2,q)\over q}
{(Bc'P^2,q)\over q}
{\phi(q)\over \phi(qq_1q_1')}
\\
\nonumber
&&
\sum_{r(qq_1q_1')}'
\bo_{\{
-\ell \equiv Pr (B  a +D  c) (\mod q_1)
\atop
-k \equiv Pr (A c +B  d) (\mod q_1)
\}}
\bo_{\{
-\ell' \equiv Pr (B  a' +D  c') (\mod q_1')
\atop
-k' \equiv Pr (A c' +B  d') (\mod q_1')
\}}
\\
\label{eq:SqSq1}
&&
\hskip1in
\times \ e_{qq_1q_1'}(r J + \bar r K + L)
,
\eea
where
\beann
J &:=&
 q_1q_1' (A a +B b +C c +D d) 
 - q_1q_1' (A a' +B b' +C c' +D d') 
\\
&&
-  q_1' \bar E  (q_1 E a A +a B^2 d P^2 +A c^2 D P^2 +q_1 E d D )
\\
&&
+  q_1 \bar E'  (q_1' E' a' A +a' B^2 d' P^2 +A c'^2 D P^2 +q_1' E' d' D )
,
\eeann

$$
K:=
-  q_1' \bar E  \ell k
+  q_1 \bar E'  \ell' k'
,
$$
and
$$
L :=
-  q_1' \bar E P (A c \ell  + B a k+B d \ell  + D c  k )
+  q_1 \bar E' P (A c' \ell'  + B a' k' + B d' \ell'  + D c'  k' ).
$$
Moreover,
\be\label{eq:Xmodq2q2p}
J   \equiv 
\bar E \bar E'
BP^4
(c-c')
(
B^2
-\gD c c'
)
\qquad
(\mod (q_2,q_2'))
,
\ee
where $\gD=AD-BC$.

Note that for every value of $r(\mod qq_1q_1')$ occurring in \eqref{eq:SqSq1}, we have that 
\be\label{eq:SqSq1a}
rJ+\bar r K + L \equiv 0(\mod q_1q_1').
\ee
\end{lem}
\pf
Inserting \lemref{lem:cSq}, and extending the $r$ sum to modulus $qq_1q_1'$ (which overcounts by a factor of $\phi(qq_1q_1')/\phi(q)$), we have that
\beann
&&
\hskip-.5in
\sum_{r(q)}'
\cS_q(r,k,\ell;\g) 
\overline{
\cS_q(r,k',\ell';\g') }
\\
&=&
{(BcP^2,q)\over q}
{(Bc'P^2,q)\over q}
{\phi(q)\over \phi(qq_1q_1')}
\sum_{r(qq_1q_1')}'
\bo_{\{
-\ell \equiv Pr (B  a +D  c) (\mod q_1)
\atop
-k \equiv Pr (A c +B  d) (\mod q_1)
\}}
\bo_{\{
-\ell' \equiv Pr (B  a' +D  c') (\mod q_1')
\atop
-k' \equiv Pr (A c' +B  d') (\mod q_1')
\}}
\\
&&
\hskip1in
e_{qq_1q_1'}(r q_1q_1' (A a +B b +C c +D d) )
\\
&&
\hskip1in
e_{qq_1q_1'}(-r q_1q_1' (A a' +B b' +C c' +D d') )
\\
&&
\hskip1in
e_{qq_1q_1'}\left(- \bar r q_1' \bar E ({Pr(B  a+D  c) +\ell}) (Pr(A c+B  d)+k) \right)
\\
&&
\hskip1in
e_{qq_1q_1'}\left( \bar r q_1 \bar E' ({Pr(B  a'+D  c') +\ell'}) (Pr(A c'+B  d')+k') \right)
.
\eeann
By the congruence restrictions on $r$, the values of $\bar E $ and $\bar E'$ are independent of their lifts to $\Z/(qq_1q_1')$.
Collecting terms gives \eqref{eq:SqSq1}.

In the modulus $qq_1q_1'$, we do not know that, for example, $E\bar E\equiv1$, since we took arbitrary lifts. But this does hold when reduced mod $q_2$, or any divisor thereof. Therefore 
to prove \eqref{eq:Xmodq2q2p}, 
we compute $J$ mod $(q_2,q_2')$, as follows.
\beann
J
&
\equiv 
&
-
q_1'\bar E P^2 ( B^2  +\gD c^2  )
+
q_1\bar E' P^2 ( B^2  +\gD c'^2  )
\\
&
\equiv 
&
\bar E \bar E'
BP^4
(c-c')
(
B^2
-\gD c c'
)
\qquad
(\mod (q_2,q_2'))
,
\eeann
where we used $ad-bc=a'd'-b'c'=1$.

To see \eqref{eq:SqSq1a}, observe that the ``true'' exponential sum has all along been $e_q$, with periodic replacement of terms like $e_q(X)$ by $e_{qq_1}(q_1X)$ etc. 
\epf

As a corollary, we record a simplified version of this lemma.

\begin{cor}\label{cor:gcdccp}
With the same notation as \lemref{lem:cSq2}, we have:
$$
\left|
\sum_{r(q)}'
\cS_q(r,k,\ell;\g) 
\overline{
\cS_q(r,k',\ell';\g') }
\right|
\ll
{(q_1,q_1')\over q}
.
$$
\end{cor}
\pf
Returning to \lemref{lem:cSq2}, we estimate
\beann
\left|
\sum_{r(q)}'
\cS_q(r,k,\ell;\g) 
\overline{
\cS_q(r,k',\ell';\g') }
\right|
&\le &
{1\over q^2}
\sum_{r(qq_1q_1')}'
\bo_{\{
-\ell \equiv Pr (B  a +D  c) (\mod q_1)
\atop
-\ell' \equiv Pr (B  a' +D  c') (\mod q_1')
\}}
,
\eeann
where we used that $\phi(q)/\phi(qq_1q_1')=1/(q_1q_1')$, since every prime dividing $q_1q_1'$ also divides $q$.

Since $B,P$ are fixed throughout, consider the condition 
$$
-\ell \equiv Pr (B  a +D  c) (\mod q_1), \ \Longrightarrow
-\ell \equiv r P B  a (\mod (q_1,c)).
$$
Since $\det\g=1$, we have that $a$ is invertible mod $c$, so $r$ is restricted to a bounded (in terms of $B,P$) number of residues mod $q_1$.
Similarly, $r$ is also restricted to a bounded number of residue classes mod $q_1'$. Therefore the total number of $r$ mod $qq_1q_1'$ satisfying the congruence is at most $qq_1q_1'/[q_1,q_1']$.
This gives the claim.
\epf

%
Next we record a Kloosterman-type estimate necessary in what follows.
\begin{lem}\label{lem:Kloost}
Fix any $J,K,L\in\Z$ and let $q_0\mid q$. 
Then
$$
\left|
\sum_{r(q)\atop r\equiv r_0(\mod q_0)}'
e_q\left(
{Jr + K\bar r + L}
\right)
\right| 
\ll_\vep
\min\left(\frac q{q_0},
q^{3/4+\vep}
\frac {1}{q_0^{1/4}}
\gcd(q/q_0,J,K)^{1/4}
\right).
$$
\end{lem}
\pf
If $q_0=1$, this is just Kloosterman's estimate, so we assume $q_0>1$.
The first bound in the minimum is just the trivial bound, and  is sometimes better than the second bound.

Following Kloosterman's method, we take the fourth moment, and consider 
\be\label{eq:Kloo1}
\cU \ := \
\sum_{J',K'\mod q}
\left|
\sum_{r(q)\atop r\equiv r_0(\mod q_0)}'
e_{q}\left(
{J'r + K'\bar r + L}
\right)
\right|^4.
\ee
We open the power and evaluate.
\be\label{eq:Kloo2}
\cU 
\ = \ 
\sum_{J',K'\mod q}
\sum_{r_1,r_2,r_3,r_4(q)\atop r_j\equiv r_0(\mod q_0)}'
e_{ q}\left(
J'(r_1+r_2-r_3-r_4)
+
K'(\bar r_1+\bar r_2-\bar r_3-\bar r_4)
\right)
\ee
The $J',K'$ sum is a complete sum over all of $\Z/q$,
which vanishes unless
$$
{r_1+r_2-r_3-r_4}\equiv 0(\mod  q), \ \ 
{\bar r_1+\bar r_2-\bar r_3-\bar r_4}\equiv 0 (\mod q),
$$
in which case they contribute $q$ each. So we have that
\beann
\cU 
&=&
q^2
\sum_{r_1,r_2,r_3,r_4(q)\atop r_j\equiv r_0(\mod q_0)}'
\bo_{\{
{r_1+r_2-r_3-r_4}\equiv 0(\mod  q) \atop
{\bar r_1+\bar r_2-\bar r_3-\bar r_4}\equiv 0 (\mod  q)
\}}
\eeann
We need 
to count the number of $r_j$ contributing to the remaining sum. 
The count is multiplicative, so we may assume that $q$ is a prime power.
Set  $\widetilde q:= q/q_0$, and 
let $R\in\Z/\widetilde q$ be defined by: $ r_1-r_3=Rq_0$; the first condition above is that we also have $r_4-r_2\equiv Rq_0$. For the second condition on the $r_j$, we multiply through by $r_1r_2r_3r_4$, getting the condition:
$$
\left(r_2-r_3\right) R \left(q_0 R+r_2+r_3\right) \equiv0(\mod \widetilde q).
$$
Recall that $q_0>1$, and notice that $q_0 R+r_2+r_3\equiv 2r_0(\mod q_0)$ is then invertible $\mod q$ (except perhaps when $2\mid q$, in which case an extra constant factor  contributes to the estimate below).
We now evaluate the count as follows. First sum  over divisors $\fq\mid \widetilde q$, then over those $R$ with $(R,\widetilde q)=  \fq$. The above condition becomes
$$
r_3\equiv r_2(\mod \widetilde q/\fq).
$$
Then  $r_3$ has $\ll q \fq/\widetilde q=q_0\fq$ possible values. In total, we have:
$$
\sum_{r_1,r_2,r_3,r_4(q)\atop r_j\equiv r_0(\mod q_0)}'
\bo_{\{
{r_1+r_2-r_3-r_4}\equiv 0(\mod  q) \atop
{\bar r_1+\bar r_2-\bar r_3-\bar r_4}\equiv 0 (\mod q)
\}}
\ll
\sum_{\fq\mid \widetilde q}
\sum_{R (\mod \widetilde q)\atop (R,\widetilde q)=\fq}
\sum_{r_2(\mod q)\atop r_2\equiv r_0(\mod q_0)}
\sum_{r_3(\mod q)\atop r_3\equiv r_0(\mod q_0)}
\bo_{\{
r_3\equiv r_2(\mod \widetilde q/\fq)
\}}
$$
$$
\ll
\sum_{\fq\mid \widetilde q}
{\widetilde q\over \fq}
\widetilde q
q_0\fq
\ll_\vep
{q^{2+\vep}\over q_0}.
$$
In summary, we obtain the following estimate:
\beann
\cU 
&\ll_\vep&
q^2
\frac {q^{2+\vep}}{q_0}
.
\eeann
Next we determine the multiplicity of the size of the original sum (that is, when $(J',K')=(J,K)$) contributing to $\cU$. Any change of variables $r\mapsto rs$ with $s\in(\Z/q)^\times$ and $s\equiv 1(\mod q_0)$ corresponds to a change in the coefficients $\iota_s:(J,K)\mapsto(Js,K\bar s)$.
Another invariance comes from the map $\gs_{u,v}:(J,K)\mapsto (J+u\widetilde q,K+v\widetilde q)$, because
$$
\sum_{r(q)\atop r\equiv r_0(\mod q_0)}'
e_{q}\left(
{(J+u\widetilde q) r + (K+v\widetilde q)\bar r + L}
\right)
=
e_{q_0}\left(
{u r_0 + v \bar r_0}
\right)
\sum_{r(q)\atop r\equiv r_0(\mod q_0)}'
e_{q}\left(
{J r + K\bar r + L}
\right)
,
$$
with both sides having the same magnitude. 

Next we must determine the number of distinct $(J',K')$ obtained by the above transformations which contribute the same magnitude to $\cU$ as $(J,K)$. Assume that $\gcd(J,q)\le\gcd(K,q)$. We use $\iota_s$ to produce as many values of $J'$ as possible, and for each such, we use $\gs_{0,v}$ to construct distinct $K'$s. 

The $s\in(\Z/q)^\times$ with $s\equiv1(\mod q_0)$ which give distinct values of $Js(\mod q)$ are determined by solving
$$
J\equiv Js(q).
$$
The number of distinct values of $Js(\mod q)$ is then  $q/\gcd(q,q_0J)$.
For each such value of $sJ$, applying $\gs_{0,v}$ produces a distinct pair $(J',K')$ where $v$ ranges in $\Z/q_0$. 

  In total, we have that:
  $$
{
q 
\over 
\gcd(q,q_0J,q_0K)
}
q_0
\left|
\sum_{r(q)\atop r\equiv r_0(\mod q_0)}'
e_q\left(
{Jr + K\bar r + L}
\right)
\right|  ^4
\le
\cU
\ll_\vep
q^2
\frac {q^{2+\vep}}{q_0}
,
  $$
  from which the claim follows.
\epf

\begin{lem}\label{lem:cSq3}
With the same notation as \lemref{lem:cSq2}, we have that
\beann
\sum_{r(q)}'
\cS_q(r,k,\ell;\g) 
\overline{
\cS_q(r,k',\ell';\g') }
&\ll_\vep&
q^{-5/4+\vep}
(q_1q_1')^{1/2}(q_1,q_1')^{1/4}
\gcd\left(q(q_1,q_1'),J, K\right)^{1/4}
,
\eeann
for any $\vep>0$.
\end{lem}
\pf
Applying \eqref{eq:SqSq1}, we decompose the sum on $r$ mod $qq_1q_1'$ into residue classes mod $q_0:=[q_1,q_1']$ to catch the indicator functions. 
\bea
\nonumber
&&
\hskip-.5in
\sum_{r(q)}'
\cS_q(r,k,\ell;\g) 
\overline{
\cS_q(r,k',\ell';\g') }
\\
\nonumber
&=&
{(BcP^2,q)\over q}
{(Bc'P^2,q)\over q}
{\phi(q)\over \phi(qq_1q_1')}
\sum_{r_0(q_0)}'
\bo_{\{
-\ell \equiv Pr_0 (B  a +D  c) (\mod q_1)
\atop
-k \equiv Pr_0(A c +B  d) (\mod q_1)
\}}
\bo_{\{
-\ell' \equiv Pr_0 (B  a' +D  c') (\mod q_1')
\atop
-k' \equiv Pr_0 (A c' +B  d') (\mod q_1')
\}}
\\
\label{eq:cSest}
&&
\hskip1in
\sum_{r(qq_1q_1')\atop r\equiv r_0(q_0)}'
e_{qq_1q_1'}(r J + \bar r K +L)
.
\eea
On the last summation, we apply \lemref{lem:Kloost}.
\beann
&\ll_\vep&
{(BcP^2,q)\over q}
{(Bc'P^2,q)\over q}
{\phi(q)\over \phi(qq_1q_1')}
\sum_{r_0(q_0)}'
\bo_{\{
-\ell \equiv Pr_0 (B  a +D  c) (\mod q_1)
\atop
-k \equiv Pr_0(A c +B  d) (\mod q_1)
\}}
\bo_{\{
-\ell' \equiv Pr_0 (B  a' +D  c') (\mod q_1')
\atop
-k' \equiv Pr_0 (A c' +B  d') (\mod q_1')
\}}
\\
&&
\hskip1in
(qq_1q_1')^{3/4+\vep} \frac1{q_0^{1/4}} \gcd(qq_1q_1'/q_0,J, K)^{1/4}
.
\eeann
Finally, we estimate the number of $r_0$ contributing.
Recall that $q_1=(BcP^2,q)$. Since $B,P$ are fixed throughout, consider the condition 
$$
-\ell \equiv Pr_0 (B  a +D  c) (\mod q_1), \ \Longrightarrow
-\ell \equiv r_0 P B  a (\mod (q_1,c)).
$$
Since $\det\g=1$, we have that $a$ is invertible mod $c$, so $r_0$ is restricted to a bounded (in terms of $B,P$) number of residues mod $q_1$.
Similarly, $r_0$ is bounded mod $q_1'$, and hence the sum on $r_0$ has a bounded number of contributions.
The claim follows immediately.
\epf

If $c\neq c'$, this analysis will suffice. But we need more work if $c=c'$, since then $J$ in \eqref{eq:Xmodq2q2p} will be $0$ mod $q_2$. (Note here that in this case, $q_1=q_1', \ q_2=q_2'$, and $E=E'$.)

\begin{lemma}\label{lem:cS4}
With notation as in  \lemref{lem:cSq2} and assuming $c=c'$, we have that:
$$
\sum_{r(q)}'
\cS_q(r,k,\ell;\g) 
\overline{
\cS_q(r,k',\ell';\g') }
\ \ll_\vep\ 
q^{-5/4+\vep}
q_1^{2}
\gcd\left(q,    \ell' k' - \ell k\right)^{1/4},
$$
\end{lemma}
\pf
Applying \lemref{lem:cSq3} gives
\beann
\sum_{r(q)}'
\cS_q(r,k,\ell;\g) 
\overline{
\cS_q(r,k',\ell';\g') }
\ \ll_\vep\ 
q^{-5/4+\vep}
q_1^{3/2}
\gcd\left(q, K\right)^{1/4},
\eeann
where
$$
K=
  q_1 \bar E (  \ell' k' - \ell k),
$$
which gives the claim on bounding $\gcd(q,\bar E)$ by $q_1$, since $E$ is invertible mod $q_2$.
\epf

\

This suffices as long as $k\ell\neq k'\ell'$. In the final case that both $c=c'$ and $k\ell=k'\ell'$, we have
\begin{lem}\label{lem:cS6}
Assume that $c=c'$ and $k\ell=k'\ell'$. Then
$$
\sum_{r(q)}'
\cS_q(r,k,\ell;\g) 
\overline{
\cS_q(r,k',\ell';\g') }
\ll
\left({(BcP^2,q)\over q}\right)^2
\sum_{r(q)}'
\bo_{\{
-\ell \equiv Pr
(B  a
+D  c)
(\mod q_1)
\atop
-k
\equiv
Pr
(A c
+B  d)
(\mod q_1)
\}}
.
$$
\end{lem}
\pf
Returning to \lemref{lem:cSq}, we estimate the $r$ sum trivially:
\beann
\sum_{r(q)}'
\left|
\cS_q(r,k,\ell;\g) 
\overline{
\cS_q(r,k',\ell';\g') }
\right|
&\le&
\left({(BcP^2,q)\over q}\right)^2
\sum_{r(q)}'
\bo_{\{
-\ell \equiv Pr
(B  a
+D  c)
(\mod q_1)
\atop
-k
\equiv
Pr
(A c
+B  d)
(\mod q_1)
\}}
\bo_{\{
-\ell' \equiv Pr
(B  a'
+D  c)
(\mod q_1)
\atop
-k'
\equiv
Pr
(A c
+B  d')
(\mod q_1)
\}}
,
\eeann
which gives the claim (crudely).
\epf

Now we need an estimate where we average over $q$  itself.
To this end, we  will first need the following result.

\begin{lem}\label{lem:newLem}
Given positive integers $R\ge S\ge W, \ U, V, X$, we have that 
\be\label{eq:newLem}
\sum_{(q,U)=1,q\equiv V(W)\atop {R \le q \le R +S}}
{\phi (q) 
\over q^2}
e_{q}(\bar U X)
\
\ll_\vep
\
{SU^{\vep}(U,(U,X)W)
\over R UW}
+
{S^2
\over R^2}
S^\vep
+
\frac {[U,W]\log R}R
+
{XS\over R^2 U W}
,
\ee
for every $\vep>0$.
\end{lem}
\pf
If $\gcd(V,W)$ is not coprime to $U$, then the sum is empty, whence \eqref{eq:newLem} holds trivially. So we assume that $(V,W,U)=1$.

The trivial bound is $S/(RW)$.
Since $(q,U)=1$, there exist $x,y$ with
$$
qx+Uy=1, \text{ or } \frac xU + \frac yq=\frac 1{qU},
$$
so that $\bar U\equiv y(\mod q)$. Then
$$
e_q(\bar U X)=e(\frac yq X)=e(-\frac xU X ) + O(X/(qU))=e_U(-xX ) + O(X/(qU)).
$$
We have that the left hand side of \eqref{eq:newLem} is:
\be\label{eq:LHS}
LHS =
\sum_{(q,U)=1,q\equiv V(W)\atop {R \le q \le R +S}}
{\phi (q) 
\over q^2}
e_U(-\bar q X ) + 
O\left(
{XS\over R^2 U W}
\right)
.
\ee
Leaving the last error term aside,
we break $q$ into residue classes mod $U_1 := [U,W]$.
\be\label{eq:eq062903}
LHS_1 =
\sum_{q_0(\mod U_1)\atop (q_0,U)=1, q_0\equiv V(W)}
e_U(-\bar q_0 X ) 
\left[
\sum_{q\equiv q_0(U_1)\atop {R \le q \le R +S}}
{\phi (q) 
\over q^2}
\right]
.
\ee
Our point will be that the bracketed sum is independent of $q_0$, to first order, and therefore we can get cancellation from the first $q_0$ sum.
To this end, next use M\"obius inversion in the form $\phi(n)=n\sum_{d\mid n}\mu(d)/d$, 
\be\label{eq:eq062902}
\Bigg[\cdot \Bigg] =
\sum_{q\equiv q_0(U_1)\atop {R \le q \le R +S}}
{1
\over q^2}
q\sum_{d\mid q}{\mu(d)\over d}
=
\sum_{d\le R+S}{\mu(d)\over d}
\sum_{q\equiv q_0(U_1), q\equiv0(d)\atop {R \le  q \le R +S}}
{1
\over q}
.
\ee
Introduce a parameter $0<D<R+S$ and break the sum on $d$ according to $d\le D$ or not.
We deal with the large $d$ first.
\beann
\sum_{ D<d\le R+S}{\mu(d)\over d}
\sum_{qd\equiv q_0(U_1)\atop {R \le  dq \le R +S}}
{1
\over dq}
&\ll&
\frac 1R
\sum_{ D<d\le R+S}{1\over d}
\sum_{{R \le  dq \le R +S}}
1
\\
&=&
\frac 1R
\sum_{ D<d\le R+S}{1\over d}
\left({S\over d}+1\right)
\\
&\ll &
{S\over R D}+
\frac {\log R}R
,
\eeann
which saves either $D$ or $S$ over the trivial bound. 

Next we handle small $d$'s. Observe that since $(q_0,U)=1$, we have that $(q_0,U_1)\mid W$. But we must also have $q_0\equiv V(W)$, and thus $(q_0,U_1)=(V,W)=:V_1$, say.
Since $q\equiv0(d)$, we let $t:=q/d$. Then
\beann
\sum_{d\le D}{\mu(d)\over d}
\sum_{td\equiv q_0(U_1)\atop {R \le  dt \le R +S}}
{1
\over dt}
&=&
\sum_{d\le D}{\mu(d)\over d^2}
\sum_{td\equiv q_0(U_1)\atop {\frac Rd \le  t \le \frac Rd +\frac Sd}}
{1
\over t}
.
\eeann
The condition $dt\equiv q_0(U_1)$ admits a solution in $t$ iff $d_1:=(d,U_1)$ divides $(q_0,U_1)=V_1$.
Write $d=d_1d_2$ with $(d_2,U_1/d_1)=1$.
For $d$ satisfying this condition, the restriction on $t$ becomes $t\equiv \bar d_2 (q_0/d_1) (\mod U_1/d_1),$ which of course is now uniquely determined mod $U_1/d_1$.
Thus
\bea
\nonumber
\sum_{d\le D}{\mu(d)\over d}
\sum_{td\equiv q_0(U_1)\atop {R \le  dt \le R +S}}
{1
\over dt}
&=&
\sum_{d\le D, d=d_1d_2\atop d_1=(d,U_1), d_1\mid V_1}{\mu(d)\over d^2}
\sum_{t\equiv \bar d_2 (q_0/d_1) (\mod U_1/d_1)\atop {\frac Rd \le  t \le \frac Rd +\frac Sd}}
{1
\over t}
\\
\nonumber
&=&
\sum_{d\le D, d=d_1d_2\atop d_1=(d,U_1), d_1\mid V_1}{\mu(d)\over d^2}
\sum_{t\equiv \bar d_2 (q_0/d_1) (\mod U_1/d_1)\atop {\frac Rd \le  t \le \frac Rd +\frac Sd}}
{1
\over R/d}
(1+O(\frac SR))
\\
\nonumber
&=&
{1
\over R}
\sum_{d\le D, d=d_1d_2\atop {d_1=(d,U_1), d_1\mid V_1\atop (d_2,U_1/d_1)=1}}{\mu(d)\over d}
\left(
{S d_1\over d U_1}
+O(1)
\right)
(1+O(\frac SR))
.
\\
\label{eq:eq06291}
\eea
The conditions $(d_2,U_1/d_1)=1$ and $(d_2,d_1)=1$ (from M\"obius) are together equivalent to $(d_2,U_1)=1$. This allows to separate the $d_1$ and $d_2$ sums, and extend the $d_2$ sum to infinity.
The ``main'' contribution becomes:
\beann
{1
\over R}
\sum_{d\le D, d=d_1d_2\atop {d_1=(d,U_1), d_1\mid V_1\atop (d_2,U_1/d_1)=1}}{\mu(d)\over d}
{S d_1\over d U_1}
&=&
{S
\over RU_1}
\sum_{d_1\mid V_1}
{ \mu(d_1)\over d_1}
\sum_{d_2\le D/d_1 \atop 
(d_2,U_1)=1
}{\mu(d_2)\over d_2^2}
\\
&=&
{S
\over RU_1}
\sum_{d_1\mid V_1}
{ \mu(d_1)\over d_1}
\left[
\sum_{d_2\le \infty \atop
(d_2,U_1)=1
}{\mu(d_2)\over d_2^2}
+
O(d_1/D)
\right]
\\
&=&
{S
\over RU_1}
\sum_{d_1\mid V_1}
{ \mu(d_1)\over d_1}
\left[
M_{U_1}
+
O(d_1/D)
\right]
,
\eeann
where
\be\label{eq:MUbnd}
M_{U_1} :=
 \prod_{p, p\nmid U_1}\left(1-\frac1{p^2}\right)
 \ \asymp \ 1
 ,
\ee
is $1/\gz(2)$ with the primes of $U_1$ removed.

Continuing the analysis gives
\beann
&=&
{S
\over RU_1}
\sum_{d_1\mid V_1}
 {\mu(d_1)\over d_1}
M_{U_1}
+
O_\vep\left(
{SW^\vep
\over DRU_1}
\right)
\\
&=&
{S
\over RU_1}
M_{U_1} V_2
+
O_\vep\left(
{SW^\vep
\over DRU_1}
\right)
,
\eeann
where
$$
V_2:=
\prod_{p|V_1}
\left(
1-\frac1p
\right)
$$
satisfies
\be\label{eq:V2bnd}
V_1^{-\vep} \ll_\vep V_2 \le 1.
\ee

We return to handle the error terms of \eqref{eq:eq06291}. The first is:
\beann
{1
\over R}
\sum_{d\le D, d=d_1d_2\atop {d_1=(d,U_1), d_1\mid V_1\atop (d_2,U_1/d_1)=1}}{1\over d}
\ll
{\log D
\over R}
,
\eeann
saving about $S/W$ over the trivial bound.
The second is:
\beann
{1
\over R}
\sum_{d\le D, d=d_1d_2\atop {d_1=(d,U_1), d_1\mid V_1\atop (d_2,U_1/d_1)=1}}{1\over d}
{S d_1\over d U_1}
\frac SR
\ll_\vep
{S^2
\over R^2U_1}
D^\vep
,
\eeann
which saves about $R/S$.

Putting everything together into \eqref{eq:eq062902} gives:
\beann
\eqref{eq:eq062902}
&=&
{S
\over RU_1}
M_{U_1} V_2
+
O_\vep\left(
{SW^\vep
\over DRU_1}
+
{\log D
\over R}
+
{S^2
\over R^2U_1}
D^\vep
\right)
+
O\left(
{S\over R D}+
\frac {\log R}R
\right)
\\
&=&
{S
\over RU_1}
M_{U_1} V_2
+
O_\vep\left(
{S^2
\over R^2U_1}
D^\vep
+
{S\over R D}+
\frac {\log R}R
\right)
,
\eeann
which is, at least in the main term, independent of $q_0$, as desired. 
Inserting this into \eqref{eq:eq062903} now gives
\bea
\nonumber
LHS_1 &=&
\sum_{q_0(\mod U_1)\atop {(q_0,U)=1, q_0\equiv V(W) 
}}
e_U(-\bar q_0 X ) 
\left[
{S
\over RU_1}
M_{U_1} V_2
+
O_\vep\left(
{S^2
\over R^2U_1}
D^\vep
+
{S\over R D}+
\frac {\log R}R
\right)
\right]
\\
\nonumber
 &=&
{S
\over RU_1}
M_{U_1} V_2
\left[
\sum_{q_0(\mod U_1)\atop {(q_0,U)=1, q_0\equiv V(W) 
}}
e_U(-\bar q_0 X ) 
\right]
+
O_\vep\left(
{S^2
\over R^2}
D^\vep
+
{S U_1\over R D}+
\frac {U_1\log R}R
\right).
\\
\label{eq:070101}
\eea
We analyze the bracketed summation by first decomposing it into residue classes mod $U$:
$$
\Bigg[\cdot\Bigg]=
\sum_{q\mod U \atop q\equiv V(\mod (W,U))}'
e_U(-\bar q X ) 
\sum_{q_0(\mod U_1)\atop {q_0\equiv q(U), q_0\equiv V(W) 
}}
1.
$$
Using $U_1=[U,W]$ and the compatibility condition $q\equiv V(\mod (W,U))$,  the Chinese Remainder Theorem gives that the last summation has exactly one $q_0$ contributing. Therefore only the first summation remains.
Let $X_1 := (U,X)$ and write $X=X_1 X_2$ and $U=X_1 U_2$ with $(X_2,U_2)=1$.
Then we break into residues mod $U_2$
$$
\Bigg[\cdot\Bigg]=
\sum_{q_2(U_2)\atop q_2\equiv V(\mod (W,U_2))}'
e_{U_2}(-\bar q_2 X_2) 
\sum_{q\mod U \atop q\equiv q_2(U_2), q\equiv V(\mod (W,U))}'
1
.
$$
In the last summation, we again get a unique contribution from the compatibility condition on $q_2$ which together with the Chinese Remainder Theorem determines $q$ uniquely mod $[U_2,(W,U)]$. Therefore the second summation evaluates to: $\phi(U)/\phi([U_2,(W,U)])$.
In summary, we have that
$$
\Bigg[\cdot\Bigg]=
{\phi (U)\over \phi ([U_2,(W,U)])}
\sum_{q_2(U_2)\atop q_2\equiv V(\mod (W,U_2))}'
e_{U_2}(-\bar q_2 X_2) 
.
$$
Looking locally, it is easy to see that the remaining summation either vanishes or is $1$ in absolute value.

Returning to \eqref{eq:070101} and inserting the above argument gives that:
\beann
LHS_1  &\ll_\vep&
{S
\over R[U,W]}
\left[
{\phi (U)\over \phi ([U/(U,X),(W,U)])}
\right]
+
{S^2
\over R^2}
D^\vep
+
{S U_1\over R D}+
\frac {U_1\log R}R
,
\eeann
where we used \eqref{eq:MUbnd} and \eqref{eq:V2bnd}. 
We choose $D=S$ for simplicity.
Returning all the way to the left hand side of \eqref{eq:newLem}, we have from \eqref{eq:LHS} that:
\beann
LHS  &\ll_\vep&
{S
\over R[U,W]}
\left[
{\phi (U)\over \phi ([U/(U,X),(W,U)])}
\right]
+
{S^2
\over R^2}
S^\vep
+
\frac {[U,W]\log R}R
+
{XS\over R^2 U W}
.
\eeann
Finally we 
note that 
$$
\phi ([U/(U,X),(W,U)])\gg_\vep
U^{-\vep}
 ([U/(U,X),(W,U)])
 $$
from which the claim follows.
\epf

Now we can give the final estimate, as follows.

\begin{lemma}\label{lem:cS7}
Assume that $c=c'$ and $k\ell=k'\ell'$. Then for parameters $Q\ge V\to\infty$, and any $\vep>0$, we have that:
\beann
&&
\hskip-.5in
\left|
\sum_{Q \le q \le Q + V}
\sum_{r(q)}'
\cS_q(r,k,\ell;\g)
\overline{\cS_q(r,k',\ell';\g')}
\right|
\\
&&
\ll
Q^\vep
\sum_{q_1\mid BcP^2\atop E=BcP^2/q_1}
\bo_{\{d\ell\equiv ak\equiv d'\ell' \equiv a'k'(\mod (q_1,c))\}}
\cN_{q_1} 
\sum_{Q_1\mid q_1\atop {p\mid Q_1 \Longrightarrow (p^\infty,q_1)\mid Q_1 \atop (E,q_1/Q_1)=1}}
\\
&&\times
\left[
{V
(EQ_1,Z)
\over QEQ_1}
+
{V^2c
\over Q^2}
+
\frac {c^3}Q
+
{V|Z|\over Q^2}
\right]
,
\eeann
where
\be\label{eq:Zdef}
Z=Z(k,\ell,k',\ell',\g,\g'):= 
     (A c \ell  + B a k+B d \ell  + D c  k )
-    (A c \ell'  + B a' k' + B d' \ell'  + D c  k' ),
\ee
and
\be\label{eq:cNm}
\cN_m=\cN_m(k,\ell,k',\ell',\g,\g') :=\#\{r\in(\Z/m)^\times:{-\ell \equiv Pr (B  a +D  c),-k \equiv Pr (A c +B  d),\atop -\ell' \equiv Pr (B  a' +D  c), -k' \equiv Pr (A c +B  d')} \}
.
\ee
\end{lemma}
\pf
We apply \lemref{lem:cSq2}, but with the special condition $q_1=q_1'$ and $K=0$:
\bea
\nonumber
&&
\sum_{Q \le q \le Q + V}
\sum_{r(q)}'
\cS_q(r,k,\ell;\g)
\overline{\cS_q(r,k',\ell';\g')}
\\
\nonumber
&=&
\sum_{Q \le q \le Q + V}
{(BcP^2,q)^2\over q^2}
{\phi(q)\over \phi(qq_1)}
\sum_{r(qq_1)}'
\bo_{\{
-\ell \equiv Pr (B  a +D  c) (\mod q_1)
\atop
-k \equiv Pr (A c +B  d) (\mod q_1)
\}}
\bo_{\{
-\ell' \equiv Pr (B  a' +D  c) (\mod q_1)
\atop
-k' \equiv Pr (A c +B  d') (\mod q_1)
\}}
e_{qq_1}(r J + L)
,
\\
\label{eq:071201}
\eea
where now
\beann
J &=&
 q_1 (A a +B b 
 +D d) 
 - q_1 (A a' +B b' 
 +D d') 
\\
&&
-   \bar E  (q_1 E a A +a B^2 d P^2 
+q_1 E d D )
\\
&&
+   \bar E  (q_1 E a' A +a' B^2 d' P^2 
+q_1 E d' D )
,
\eeann
and
\be\label{eq:Lis}
L =
-   \bar E P (A c \ell  + B a k+B d \ell  + D c  k )
+  \bar E P (A c \ell'  + B a' k' + B d' \ell'  + D c  k' ).
\ee
The analogue of \eqref{eq:Xmodq2q2p} here becomes:
$$
J \equiv 0
(\mod q_2).
$$
Returning to \eqref{eq:Lis}, we have that
$$
L \equiv
-   \bar E P B( a k-a'k'+ d \ell-d'\ell'   )
(\mod (q_1,c)).
$$
But the restrictions on $r$ in \eqref{eq:071201} require that 
$$
Pr B\equiv   
-d\ell \equiv   
-ak \equiv   
-d'\ell' 
\equiv    
-a'k'   (\mod (q_1,c)),
$$
where we used that $ad\equiv a' d' \equiv 1 (\mod c)$. Therefore 
\be\label{eq:Lmodq1}
L\equiv 0(\mod (q_1,c)).
\ee

The sum on $r$ is multiplicative with respect to the modulus, but $q_1$ and $q_2$ are not necessarily coprime.
To fix this, introduce a new parameter
$$
Q_2 := \prod_{p^u \| q\atop p|q_2}p^u,
$$
so that $Q_2\mid q$, and $q_2\mid Q_2$.
Then let $Q_1$ be defined by 
$$
q=Q_1Q_2,
$$
and it is easy to see that $(Q_1,Q_2)=1$ and $Q_1\mid q_1$. Moreover, if $p\mid Q_1$, then $(p^\infty,q_1)\mid Q_1$; that is, the largest prime power of any prime dividing $Q_1$ occurs in $Q_1$. 
It will be convenient to define 
$$
m=(q_1,Q_2).
$$
Decompose $q_1$ further as
$$
q_1=(Q_1,q_1)(Q_2,q_1) = Q_1m .
$$
Note that $(E,q_2)=1$, and $m\mid q_2$, so $(E,q_1/Q_1)=1$.
Then
$$
qq_1=Q_1^2 Q_2 m.
$$
Let
$$
Q_2' := Q_2 m\qquad\text{and}\qquad Q_1'=Q_1^2,
$$
so that
$$
qq_1=Q_1'Q_2',
$$
with $(Q_1',Q_2')=1$.
The lift $\bar E$ can be chosen so that $E \bar E \equiv 1 (\mod Q_2')$. Observe for later use that
$$
(q_1,Q_1')=Q_1.
$$
Then the same calculation leading to \eqref{eq:Xmodq2q2p} gives
$$
J\equiv 0 (\mod Q_2').
$$
Now we split the $r$ sum according to these moduli. Let $s_1,s_2$ be determined by:
$$
{s_1\over Q_1'}+{s_2\over Q_2'}
=
{1\over Q_1'Q_2'}
=
{1\over qq_1},
$$
that is, $s_1 Q_2'\equiv 1(\mod Q_1')$, and $s_2\equiv \overline{Q_1'}\mod Q_2'$, which implies that
\be\label{eq:s2}
s_2\equiv \overline{Q_1}^2\mod Q_2
.
\ee
We will also need the basic fact that, if $a\mid b$, then $\phi(ab)=\phi(b)\cdot a$.
Then we can write
\bea
\nonumber
&&
\hskip-.5in
\sum_{Q \le q \le Q + V}
\sum_{r(q)}'
\cS_q(r,k,\ell;\g)
\overline{\cS_q(r,k',\ell';\g')}
\\
\nonumber
&=&
\sum_{E\mid BcP^2\atop q_1=BcP^2/E}
\sum_{Q \le q \le Q + V\atop (q,BcP^2)=q_1,\ qq_1=Q_1'Q_2'}
\bo_{\{d\ell\equiv ak\equiv d'\ell' \equiv a'k'(\mod (q_1,c))\}}
{q_1^2\over q^2}
{\phi(q)\over \phi(q)\cdot q_1}
\\
\nonumber
&&\times
\left(
\sum_{r(Q_1')}'
\bo_{\{
-\ell \equiv Pr (B  a +D  c) (\mod (q_1,Q_1'))
\atop
-k \equiv Pr (A c +B  d) (\mod (q_1,Q_1'))
\}}
\bo_{\{
-\ell' \equiv Pr (B  a' +D  c) (\mod (q_1,Q_1'))
\atop
-k' \equiv Pr (A c +B  d') (\mod (q_1,Q_1'))
\}}
e_{Q_1'}(s_1(r J + L))
\right)
\\
\nonumber
&&
\times
\left(
e_{Q_2'}(s_2 L)
\sum_{r(Q_2')}'
\bo_{\{
-\ell \equiv Pr (B  a +D  c) (\mod (q_1,Q_2'))
\atop
-k \equiv Pr (A c +B  d) (\mod (q_1,Q_2'))
\}}
\bo_{\{
-\ell' \equiv Pr (B  a' +D  c) (\mod (q_1,Q_2'))
\atop
-k' \equiv Pr (A c +B  d') (\mod (q_1,Q_2'))
\}}
\right)
\\
\nonumber
&=&
\sum_{E\mid BcP^2\atop q_1=BcP^2/E}
\bo_{\{d\ell\equiv ak\equiv d'\ell' \equiv a'k'(\mod (q_1,c))\}}
\sum_{Q_1\mid q_1\atop {p\mid Q_1 \Longrightarrow (p^\infty,q_1)\mid Q_1\atop {Q_1'=Q_1^2 \atop (E,q_1/Q_1)=1}}}
{q_1
\over Q_1^2}
\sum_{(Q_2,E)=(Q_2,Q_1)=1,Q_1Q_2\equiv0(q_1)\atop {Q \le Q_1Q_2 \le Q + V\atop {q=Q_1Q_2\atop {q_2=Q_1Q_2/q_1\atop Q_2'=qq_1/Q_1'}}}}
{1\over Q_2^2}
\\
\label{eq:070601}
&&
\times
\left(
\sum_{r_1(Q_1')}'
\bo_{\{
-\ell \equiv Pr_1 (B  a +D  c) (\mod Q_1)
\atop
-k \equiv Pr_1 (A c +B  d) (\mod Q_1)
\}}
\bo_{\{
-\ell' \equiv Pr_1 (B  a' +D  c) (\mod Q_1)
\atop
-k' \equiv Pr_1 (A c +B  d') (\mod Q_1)
\}}
e_{Q_1'}(\overline{Q_2'}(r_1 J + L))
\right)
\\
\nonumber
&&
\times
\left(
e_{Q_2'}(s_2 L)
\sum_{r(Q_2')}'
\bo_{\{
-\ell \equiv Pr (B  a +D  c) (\mod (q_1,Q_2'))
\atop
-k \equiv Pr (A c +B  d) (\mod (q_1,Q_2'))
\}}
\bo_{\{
-\ell' \equiv Pr (B  a' +D  c) (\mod (q_1,Q_2'))
\atop
-k' \equiv Pr (A c +B  d') (\mod (q_1,Q_2'))
\}}
\right)
.
\eea
Next we make the following two claims: $(i)$ that the first sum on $r_1$ only depends only on the value of $Q_2$ modulo $q_1$; and $(ii)$,
that we can count the number of solutions in the $r$ sum mod $Q_2'$.
We first work on $(ii)$. Observe that
$$
m=(q_1,Q_2)=(q_1,Q_2') = {q_1\over Q_1},
$$ 
which is a divisor of $Q_2'$. That is, $(m,Q_1)=1$. 
Recalling the definition \eqref{eq:cNm} of $\cN_m$, 
 the sum on $r(Q_2')$ clearly contributes
$$
\sum_{r(Q_2')} ' = \cN_m {\phi (Q_2')\over \phi(m)}
=
\phi (Q_2) \cN_m { (q_1,Q_2)\over \phi(
(q_1,Q_2)
)}
.
$$
Next we argue $(i)$. Recall from the analogue of \eqref{eq:SqSq1a} in this setting that any $r_1$ occurring in the first summation satisfies: 
$$
r_1 J + L \equiv 0 (Q_1).
$$
Therefore the $r_1$ summation in \eqref{eq:070601} is:
$$
\left(
Q_1
\sum_{r_1(Q_1)}'
\bo_{\{
-\ell \equiv Pr_1 (B  a +D  c) (\mod Q_1)
\atop
-k \equiv Pr_1 (A c +B  d) (\mod Q_1)
\}}
\bo_{\{
-\ell' \equiv Pr_1 (B  a' +D  c) (\mod Q_1)
\atop
-k' \equiv Pr_1 (A c +B  d') (\mod Q_1)
\}}
e_{Q_1}(\overline{Q_2}\overline{(Q_2,q_1)}{r_1 J + L\over Q_1})
\right)
.
$$
The term $\overline Q_2$ only depends on the residue class, $Q_2^0$, say, mod $Q_1$, so we break the sum according to these residue classes.
Returning to the original expression, we have:
\beann
&&
\sum_{Q \le q \le Q + V}
\sum_{r(q)}'
\cS_q(r,k,\ell;\g)
\overline{\cS_q(r,k',\ell';\g')}
\\
&=&
\sum_{E\mid BcP^2\atop q_1=BcP^2/E}
\bo_{\{d\ell\equiv ak\equiv d'\ell' \equiv a'k'(\mod (q_1,c))\}}
\sum_{Q_1\mid q_1\atop {p\mid Q_1 \Longrightarrow (p^\infty,q_1)\mid Q_1\atop {Q_1'=Q_1^2 \atop {m=q_1/Q_1, \ (m,Q_1)=1\atop  (E,q_1/Q_1)=1}}}}
{q_1
\over Q_1^2}
\cN_m { m\over \phi(m)}
\\
&&
\times
\sum_{Q_2^0 (\mod Q_1)'}
\left(
Q_1
\sum_{r_1(Q_1)}'
\bo_{\{
-\ell \equiv Pr_1 (B  a +D  c) (\mod Q_1)
\atop
-k \equiv Pr_1 (A c +B  d) (\mod Q_1)
\}}
\bo_{\{
-\ell' \equiv Pr_1 (B  a' +D  c) (\mod Q_1)
\atop
-k' \equiv Pr_1 (A c +B  d') (\mod Q_1)
\}}
e_{Q_1}(\overline{Q_2^0}\overline{m}{r_1 J + L\over Q_1})
\right)
\\
&&
\times
\left[
\sum_{(Q_2,EQ_1)=1,Q_2\equiv0(m), Q_2\equiv Q_2^0(\mod Q_1)\atop {Q \le Q_1Q_2 \le Q + V\atop {q=Q_1Q_2\atop {q_2=Q_1Q_2/q_1\atop Q_2'=qq_1/Q_1'}}}}
{\phi (Q_2) 
\over Q_2^2}
e_{Q_2'}(s_2 L)
\right]
.
\eeann

Let
$$
Q_1'' := (Q_1,c),
$$
with 
$Q_1/Q_1'' \asymp 1.$
(Recall here that $P$ is a constant depending only on the group $\G$, and implied constants may depend on $\G$ and the fixed parameters $A,B,C,D$.)
From \eqref{eq:Lmodq1}, we have that  $L\equiv0(\mod Q_1'')$, since $Q_1''\mid(q_1,c)$.
 Let 
 $$m_1:=(m, c),
 $$ 
 so that $m/m_1 \asymp 1$. 
 By the same argument, we also have that $L\equiv 0 (\mod m_1)$, and since $(m_1,Q_1'')=1$, we have that $L\equiv 0(\mod Q_1''m_1)$.
 
  So using \eqref{eq:s2} the last sum can be written as 
\beann
\Bigg[\cdot\Bigg]&=&
\sum_{(Q_2,EQ_1)=1,Q_2\equiv0(m), Q_2\equiv Q_2^0(\mod Q_1)\atop {Q \le Q_1Q_2 \le Q + V}}
{\phi (Q_2) 
\over Q_2^2}
e_{Q_2}(s_2 {L\over m})
\\
&=&
\sum_{(Q_2,U)=1,Q_2\equiv0(m), Q_2\equiv Q_2^0(\mod Q_1)\atop {Q \le Q_1Q_2 \le Q + V}}
{\phi (Q_2) 
\over Q_2^2}
e_{Q_2}(\bar U X)
,
\eeann
where 
$U=EQ_1(Q_1/Q_1'')(m/m_1)\asymp EQ_1$, and 
$$
X=
{-    P (A c \ell  + B a k+B d \ell  + D c  k )
+   P (A c \ell'  + B a' k' + B d' \ell'  + D c  k' )\over Q_1''m_1},
$$
by \eqref{eq:Lis}.
Note that 
$$
X={PZ\over Q_1''m_1},
$$ 
where 
$Z$ is defined as in \eqref{eq:Zdef}.
It is at this point that we apply \lemref{lem:newLem}, with 
$R=Q/Q_1$, $S=V/Q_1$, and $W=[m,Q_1]=q_1$.
Note that
$$
(U,X)=K
(EQ_1,X)=K
\frac1{mQ_1''}
(EQ_1mQ_1'',PZ)
=K
\frac1{q_1}
(Eq_1Q_1,Z)
.
$$
Here $K$ is an absolute constant, not the same in each occurrence. 
Then
$$
(U,(U,X)W)=
K
(EQ_1,EQ_1q_1,Z)
\ll
(EQ_1,Z).
$$

Now applying \lemref{lem:newLem} gives:
\beann
\Bigg[\cdot\Bigg]
&\ll_\vep&
Q^\vep
\left(
{V (EQ_1,Z)
\over Q {EQ_1 }q_1}
+
{V^2
\over Q^2}
+
\frac {Q_1c}Q
+
{V|X|\over Q^2 c}
\right)
.
\eeann

Returning to the original summation, we have:
\beann
&&
\hskip-.5in
\left|
\sum_{Q \le q \le Q + V}
\sum_{r(q)}'
\cS_q(r,k,\ell;\g)
\overline{\cS_q(r,k',\ell';\g')}
\right|
\\
&\ll_\vep&
Q^\vep
\sum_{E\mid BcP^2\atop q_1=BcP^2/E}
\bo_{\{d\ell\equiv ak\equiv d'\ell' \equiv a'k'(\mod (q_1,c))\}}
\sum_{Q_1\mid q_1\atop {p\mid Q_1 \Longrightarrow (p^\infty,q_1)\mid Q_1\atop {Q_1'=Q_1^2 \atop {m=q_1/Q_1, \ (m,Q_1)=1, \atop  (E,q_1/Q_1)=1}}}}
q_1
\cN_m { m\over \phi(m)}
\\
&&
\times
\left(
\sum_{r_1(Q_1)}'
\bo_{\{
-\ell \equiv Pr_1 (B  a +D  c) (\mod Q_1)
\atop
-k \equiv Pr_1 (A c +B  d) (\mod Q_1)
\}}
\bo_{\{
-\ell' \equiv Pr_1 (B  a' +D  c) (\mod Q_1)
\atop
-k' \equiv Pr_1 (A c +B  d') (\mod Q_1)
\}}
\right)
\left[
{V
(EQ_1,Z)
\over QEQ_1q_1}
+
{V^2
\over Q^2}
+
\frac {Q_1c}Q
+
{V|X|\over Q^2 c}
\right]
.
\eeann
Note that  the $r_1(\mod Q_1)$ summation is exactly $\cN_{Q_1}$, and since $(m,Q_1)=1$, we have that $\cN_m\cdot \cN_{Q_1}=\cN_{q_1}$. Now we have, crudely, that
\beann
&&
\hskip-1in
\left|
\sum_{Q \le q \le Q + V}
\sum_{r(q)}'
\cS_q(r,k,\ell;\g)
\overline{\cS_q(r,k',\ell';\g')}
\right|
\\
&\ll_\vep&
Q^\vep
\sum_{E\mid BcP^2\atop q_1=BcP^2/E}
\bo_{\{d\ell\equiv ak\equiv d'\ell' \equiv a'k'(\mod (q_1,c))\}}
\sum_{Q_1\mid q_1\atop {p\mid Q_1 \Longrightarrow (p^\infty,q_1)\mid Q_1 \atop  (E,q_1/Q_1)=1}}
\cN_{q_1}
\\
&&
\times
\left[
{V
(EQ_1,Z)
\over QEQ_1}
+
{V^2c
\over Q^2}
+
\frac {c^3}Q
+
{V|Z|\over Q^2 }
\right]
,
\eeann
from which is the claim.
\epf

Next we analyze the size of $\cN_m$ in \eqref{eq:cNm}.
\begin{lemma}
Let $m\mid BcP^2$. Then
\be\label{eq:cNmbnd1}
\cN_m\ll 1.
\ee
\end{lemma}
\pf
Our goal is to count the number of $r$ in $(\Z/m)^\times$ satisfying 
$-\ell \equiv Pr (B  a +D  c)$, $-k \equiv Pr (A c +B  d),$ and also  $-\ell' \equiv Pr (B  a' +D  c),$ $-k' \equiv Pr (A c +B  d')$.
Let 
$$
m_1:=(m,c)
,
$$
so that $m_1\mid c$, and $m\ll m_1$ (since $B$ and $P$ are fixed). Reducing the moduli mod $m_1$ gives the equations:
$$
-\bar a\ell \equiv r P B  ,\  -\bar d k \equiv rP B , \ -\bar a'\ell' \equiv rP B , \ -\bar d'k' \equiv rP B  (\mod m_1).
$$
Here we used that $(a,c)=1$ since $\g\in\SL_2(\Z)$, etc.
There are clearly a bounded number of solutions in $r$ to the above, which
gives the claim.
\epf

And lastly, we analyze the number of elements in $\SL_2(\Z)$ with a given value of $c$ and satisfying a congruence in $Z$ in \eqref{eq:Zdef}.

\begin{lemma}
Let $\g\in\SL_2(\Z)\cap B_T$ be given with $\g_c=c$, and fix a divisor $Z_1\mid c$ with $c\ll T$. 
Also fix  $k,\ell,k',\ell'$ with $k\ell=k'\ell'$.
Then the number of $\g'\in\SL_2(\Z)$ with $\g_c'=c$ and $|a'|,|b'|,|d'|\asymp T$ and satisfying $Z\equiv 0(\mod Z_1)$ is bounded by
 \be\label{eq:Z1bnd}
\ll { T\over Z_1} (Z_1,\ell d - k a),
 \ee
 as $T\to\infty$.
\end{lemma}
\pf
In the variables $a', d'$,
we have the pair of equations: $a'd'\equiv 1(\mod Z_1)$ and
$$
Z=
     (A c \ell  + B a k+B d \ell  + D c  k )
-    (A c \ell'  + B a' k' + B d' \ell'  + D c  k' )
\equiv 0(Z_1),
$$
or
$$
Bk'  a'^2 
-(B a k+B d \ell)a'
+ B \ell'
\equiv 0(Z_1).
$$
Let $B_1 := (B,Z_1)$ and set $Z_2:=Z_1/B_1$ and $B_2:=B/B_1$, with $(B_2,Z_2)=1$. Dividing through by $B_1$, the equation reduces to 
$$
k'  a'^2 
-(a k+d \ell)a'
+ \ell'
\equiv 0(Z_2).
$$
Let
$$
\tilde Z := (k',\ell',Z_2).
$$
Since $k\ell=k'\ell'$, we have that $\tilde Z^2\mid k\ell$. Working locally, suppose that $\tilde Z\mid \ell$. Then reducing mod $\tilde Z$ gives the equation
$$
-( a k)a'
\equiv 0(\tilde Z).
$$
But $(a,\tilde Z)=(a,c)=1$ and $a'$ is also coprime to $\tilde Z$, which implies that $k\equiv 0(\tilde Z)$. Therefore there are no solutions unless both $\tilde Z\mid k$ and $\tilde Z\mid \ell$.
In this case, we can divide the whole equation by $\tilde Z$. Set $k_1:=k/\tilde Z,\ \dots, \ \ell_1':=\ell'/\tilde Z$ and $Z_3:=Z_2/\tilde Z$. Then the equation becomes
$$
k_1'  a'^2 
-(a k_1+d \ell_1)a'
+ \ell_1'
\equiv 0(Z_3),
$$
with $(k_1',\ell_1',Z_3)=1$. Working locally, we may assume that $(k_1',Z_3)=1$. Now we can simply solve the equation. Assuming for simplicity that $Z_3$ is odd (with minor modifications otherwise), we have that
$$
(  a' 
-\bar 2\bar k_1'(a k_1+d \ell_1))^2
\equiv - \bar k_1' \ell_1'
+\bar 4 \bar k_1'^2(a k_1+d \ell_1)^2
\equiv
(\bar 2 \bar k_1' (ak_1-d\ell_1))^2
(\mod Z_3),
$$
where we again used that $k\ell=k'\ell'$ and $ad\equiv1(\mod Z_3)$.
The equation is now a difference of squares, so
$$
(  a' 
-\bar k_1' \ell_1d)
(
  a' 
-\bar k_1' k_1a
)
\equiv 0
(\mod Z_3).
$$
Again working locally, suppose that $Z_3=p^U$. Then for $V+W=U$, we get a solution  
$$
a'\equiv \bar k_1'\ell_1d(\mod p^V)\text{ and }a'\equiv \bar k_1'k_1 a(\mod p^W).
$$
Assume WLOG that $V\le W$.
Then for there to be any solutions, it must be the case that $\ell_1d- k_1 a\equiv0 (\mod p^{V})$, and if this is the case, then there are $p^{U-W}=p^{V}$ solutions for $a'$.
Let
$$
Z_4: = (Z_3, \ell_1 d - k_1a).
$$
By the above discussion, the number of solutions for $a'(\mod Z_3)$ is at most 
$$
\min(Z_4, \sqrt{Z_3}).
$$
So once the value of $a'$ mod $Z_3$ is fixed, the total number of such $a'\asymp T$ is 
$$
\ll \frac T{Z_3}
\ll { T\tilde Z \over Z_2}
\ll { T(k',\ell',Z_1) \over Z_1}.
$$
Thus we can bound the total number of values of $a'$ by 
$$
Z_4 T (k',\ell',Z_1)/Z_1 \ll T (Z_3,\ell_1 d - k_1 a) (k',\ell',Z_1)/Z_1
\ll T (Z_1,\ell d - k a)  /Z_1.
$$
With $a'$ and $c$ fixed, the number of $d', b'$ is $\ll1$, since they are all of order $T$ and $a'd'-b' c=1$.
This completes the proof.
\epf

\newpage

\section{Major Arc Analysis}\label{sec:major}

\begin{thm}\label{thm:majArcs}
There is an $\eta>0$ and a  set $\sE\subset\Z$ of ``exceptional'' $n$, of zero density, 
$$
\frac1N\#(\sE \cap [-N,N]) = O(N^{-\eta}),
$$ 
such that, for $n\notin \sE$, $n\asymp N$, we have that
$$
\cM_N(n) \ \gg \ \fS(n) {\widehat{\cR_N}(0)\over N} + O\left({\widehat{\cR_N}(0)\over N} N^{-\eta}\right),
$$
as $N\to\infty$,
where, for admissible $n\notin\sE$ and for any $\vep>0$, the ``singular series'' $\fS(n)$ satisfies
$$
\fS(n) \gg_\vep |n|^{-\vep}.
$$
The implied constants are absolute. 
\end{thm}

\pf
We begin with \eqref{eq:cMdef}:
\beann
\cM_N(n) & = &
\int_0^1
\fM(\gt)
\widehat{\cR_N}(\gt)
e(-n\gt)
d\gt
\\
&=&
\int_0^1
\sum_{q<Q_0}
\sum_{r(q)}'
\sum_{m\in\Z}
\psi((\gb+m)\tfrac{N}{K_0})
\widehat{\cR_N}(\tfrac rq+\gb)
e(-n(\tfrac rq+\gb))
d\gb
\\
&=
&
\int_\R
\sum_{q<Q_0}
\sum_{r(q)}'
\widehat{\cR_N}(\tfrac rq+\gb)
e(-n(\tfrac rq+\gb))
\psi(\gb\tfrac{N}{K_0})
d\gb
.
\eeann
We have that
$$
\cR_N(\tfrac rq+\gb)
=
\sum_{x,y\in\Z}
\gU\left({x\over X}\right)
\gU\left({y\over X}\right)
\sum_{\g_0\in \G(q)\bk \G}
e_q(r \ff_{\g_0}(x,y))
\left[
\sum_{\g\in \sF_T\atop \g\equiv \g_0(\mod q)}
e(\gb \ff_\g(x,y))
\right]
.
$$
For the bracketed term, we apply \lemref{lem:Expansion}, together with $|\gb|<K_0/N<1/X^2$.
\beann
\cR_N(\tfrac rq+\gb)
&=&
\sum_{x,y\in\Z}
\gU\left({x\over X}\right)
\gU\left({y\over X}\right)
\sum_{\g_0\in \G(q)\bk \G}
e_q(r \ff_{\g_0}(x,y))
\\
&&
\left[
\frac1{[\G:\G(q)]}
\sum_{\g\in \sF_T}
e(\gb \ff_\g(x,y))
+
O(|\sF_T|N^{-\gT})
\right]
.
\eeann
Inserting this into $\cM_N$ gives:
\beann
\cM_N(n) & = &
\sum_{x,y\in\Z}
\gU\left({x\over X}\right)
\gU\left({y\over X}\right)
\left[
\sum_{q<Q_0}
\frac1{[\G:\G(q)]}
\sum_{\g_0\in \G(q)\bk \G}
\sum_{r(q)}'
e_q(r (\ff_{\g_0}(x,y)-n))
\right]
\\
&&
\times
\left[
\sum_{\g\in \sF_T}
\int_\R
\psi(\gb\tfrac{N}{K_0})
e(\gb (\ff_\g(x,y)-n))
d\gb
\right]
\\
&&
+
O\left(
{\widehat{\cR_N}(0)\over N}
N^{-\gT}
{Q_0^5K_0 }
\right)
,
\eeann
where we have split into modular and archimedean components. The proof then follows on applying  \thmref{thm:sENfSq} and \thmref{thm:fSnSize} to the  modular component, and
\lemref{lem:Arch} to the archimedean part, together with the choice of parameters in  \eqref{eq:ga0gk0}.
\epf

\newpage

\section{Minor Arc Analysis}\label{sec:minor}

Our goal is to estimate $\cE_N$ in $\ell^2$, or what is the same (by Parseval), to bound
$$
\|\cE_N\|^2 \ = \ 
\|\widehat {\cE_N}\|^2 \ = \
\int_0^1
|1-\fM(\gt)|^2
|\widehat{\cR_N}(\gt)|^2
d\gt.
$$
The main result of this section is the following.

\begin{thm}\label{thm:cEmain}
There exists some $\eta>0$ so that, as $N\to\infty$,
$$
\|\cE_N\|^2 \ll 
{|\widehat{\cR_N}(0)|^2\over N}
N^{-\eta}.
$$
\end{thm}

A standard argument concludes the main \thmref{thm:main2} (and hence \thmref{thm:main1}) from \thmref{thm:majArcs} and \thmref{thm:cEmain}.
We begin the analysis as follows.
In Dirichlet's approximation theorem, we choose the level 
$$
M \ := \ TX,
$$ 
so that for every $\gt\in[0,1]$, there is a $q<M$ and $(r,q)=1$ so that $\gt=\frac rq+\gb$ with 
$$
\left|
\gt-\frac rq
\right|
=
|\gb|
<
\frac1{qM}.
$$
Now we decompose the circle into dyadic  regions of the form
$$
W_{Q} \ := \
\left\{
\gt=\frac rq +\gb \ : \
q\asymp Q, (r,q)=1, |\gb|\ll  \frac 1{QM}
\right\}
,
$$
so that
$$
\|\cE_N\|^2 \ \ll \ 
\sum_{Q<M\atop{dyadic}}
\int_{W_{Q}}
|1-\fM(\gt)|^2
|\widehat{\cR_N}(\gt)|^2
d\gt
.
$$
This decomposes further into three ranges, according to whether $Q$ satisfies: $Q<Q_0$, or $Q_0<Q<X/Y$, or $X/Y<Q<XT=M$.
Here we have set
\be\label{eq:Ydef}
Y:=T^{(2\gd-1)/10}=N^y,\qquad y>0,
\ee
to be a small power of $N$.

 On the latter two ranges, the weight  $|1-\fM(\gt)|^2$ is exactly $1$.
To keep track, we define the integrals:
$$
\cI_{Q_0,K_0} \ : = \ \int_{\gt=\frac rq+\gb\atop q<Q_0, |\gb|<K_0/N} \left|\gb{N\over K_0}\right|^2 |\widehat{\cR_N}(\gt)|^2 d\gt,
$$
$$
\cI_{Q_0} \ : = \ \int_{\gt=\frac rq+\gb\atop q<Q_0, K_0/N\le |\gb|<1/(qM)}  |\widehat{\cR_N}(\gt)|^2 d\gt,
$$
\be\label{eq:cIQdef}
\cI_Q \ := \ 
\int_{\gt = r/q + \gb \atop Q\le q<2Q, (r,q)=1,|\gb|<1/(QM)}
\left|
\widehat{\cR_N}(\gt)
\right|^2
d\gt.
\ee

\subsection{Preliminaries}\

We first estimate $\widehat {\cR_N}(\gt)$ for $\gt=\frac rq+\gb$ as follows.
We begin by decomposing $x$ and $y$ according to their residue classes mod $q$ and applying Poisson summation in $x$ and $y$ gives:
\bea
\nonumber
\cR_N(\gt)
& = &
\sum_{\g\in \sF_T}
\sum_{x,y\in\Z}
\gU\left({x\over X}\right)
\gU\left({y\over X}\right)
e((\frac rq+\gb) \ff_\g(x,y))
\\
\label{eq:Poisson}
& = &
X^2
\sum_{\g\in \sF_T}
\sum_{k,\ell\in\Z}
\cS_q(r,k,\ell;\g)
\cJ_X(\gb, k ,\ell,q;\g)
,
\eea
where $\cS_q$ is given in \eqref{eq:cSqDef}
and
$$
\cJ_X(\gb, k ,\ell, q;\g) \ := \
\iint_{x,y\in\R}
\gU\left({x}\right)
\gU\left({y}\right)
e(\gb \ff_\g(xX,yX)
-\tfrac Xq(kx+\ell y)
)
dx dy
.
$$

\begin{lem}\label{lem:cJX}
Suppose that $q<X/Y$. Then for any $L<\infty$, we have that
\be\label{eq:cJnonStationary}
|\cJ_X(\gb, k ,\ell,q;\g)| \ \ll_L  \ X^{-L},
\ee
 unless $|k|,|\ell|\ll1$, in which case, we have:
\be\label{eq:cJstationary}
|\cJ_X(\gb, k ,\ell,q;\g)| \ll 
\min(1,\frac1{N|\gb|}).
\ee

Alternatively, if $X/Y\le q$, then we have the same arbitrary cancellation \eqref{eq:cJnonStationary}, unless $|k|,|\ell|\ll  Y \frac qX$ (in which case, we only need the trivial bound).
\end{lem}
\pf
The phase of $\cJ_X$ can be written as $e(g)$ where
$$
g(x,y)
=
 \gb \ff_\g(xX,yX)
-\frac Xq
(kx+\ell y)
.
$$
Inputting \eqref{eq:ffIs}, the partial derivatives of $g$ are
\be\label{eq:gx}
\dd_x g(x,y)
=
\frac Xq
\left(
\gb q
B   c_{\gamma }
P^2
X  y
-
k
\right)
+
O(|\gb |
T X) 
\ee
and similarly
\be\label{eq:gy}
\dd_y g(x,y)
=
\frac Xq
\left(
\gb q
B   c_{\gamma }
P^2
Xx
-
\ell
\right)
+
O(|\gb |
T X) 
.
\ee
First consider the case that $q<X/Y$.
Recalling that $|\gb|\ll 1/(QTX)$, $B,P,x,y\asymp1$, and $c_y\asymp T$, we have that
 $\nabla g\neq0$ unless
 $$
 |k|,|\ell|\ll
| \gb| q
T
X  
\ll 1
.
 $$
Outside of this range, we apply non-stationary phase, giving \eqref{eq:cJnonStationary}.

Now suppose that the pair $(k,\ell)$ is such that $\nabla g$ does vanish at some point $p$ in the support of $\gU\times\gU$.
In this case, 
we apply stationary phase to show that
$$
|\cJ_X|
\ll
\min(1,
\gD_p^{-1/2}
)
,
$$
where $\gD_p$ is absolute determinant of the Hessian of $g$ at $p$. Since
$$
\gD_p=
|\det
(\dd_{i,j}(g))| 
= 
(
|\gb|
B   c_{\gamma }
P^2
X^2)^2
\asymp
(
|\gb|
T
X^2
)^2,
$$
we have 
that 
$$
\gD_p^{-1/2}
\ll
\frac 1{TX^2|\gb|}
=
\frac1{N|\gb|}
,
$$
which gives \eqref{eq:cJstationary}. (In fact, since the form is bilinear in the variables, it is possible to evaluate the integrals explicitly, though this is not needed here.)

In the case that $X/Y\le q$, we can only apply non-stationary phase if the phase is actually growing, which is the case if $\max(|k|,|\ell|)>Y \frac q X$. This completes the proof.
\epf

\subsection{Minor Arcs I: Case $q<Q_0$}\

\begin{prop}
Assume that $q<X$. Then
\be\label{eq:LinfBnd}
\left|\widehat{\cR_N}\left(\frac rq+\gb\right)\right|
\ll
{X^2|\sF_T|\over N|\gb|}.
\ee
\end{prop}
\pf
Inserting \lemref{lem:cJX} gives
$$
\widehat{\cR_N}\left(\frac rq+\gb\right)
\ll
X^2
\sum_{\g\in \sF_T}
\sum_{k,\ell\ll1}
|\cS_q(r,k,\ell;\g)|
\frac 1 {N|\gb|},
$$
which gives the result on  trivially estimating $|\cS_q|\le1$.
\epf

\begin{cor}
We have:
$$
\cI_{Q_0,K_0} \ \ll \ 
{|\widehat{\cR_N}(0)|^2\over N}
 {Q_0^2\over K_0}
,
$$
and
$$
\cI_{Q_0} \ \ll \ 
{|\widehat{\cR_N}(0)|^2\over N} {Q_0^2\over K_0}.
$$
\end{cor}
\pf
Inserting the $L^\infty$ bound \eqref{eq:LinfBnd} gives:
\beann
\cI_{Q_0,K_0} &\ll & 
\int_{\gt=\frac rq+\gb\atop q<Q_0, |\gb|<K_0/N} \left|\gb{N\over K_0}\right|^2 {X^4|\sF_T|^2\over N^2|\gb|^2} d\gt 
\ \ll\
{|\widehat{\cR_N}(0)|^2\over K_0^2}
Q_0^2 {K_0\over N}
,
\eeann
giving the claim. Similarly,
\beann
\cI_{Q_0} &\ll &
 \int_{\gt=\frac rq+\gb\atop q<Q_0, K_0/N<|\gb|} \left|{X^2|\sF_T|\over N|\gb|}\right|^2 d\gt 
 \ll
{|\widehat{\cR_N}(0)|^2\over N^2} {Q_0^2N\over K_0},
\eeann
as claimed.
\epf

The choice of parameters \eqref{eq:ga03gk0} ensures that these are power savings, as required in \thmref{thm:cEmain}.

\subsection{Minor Arcs II: Case $Q_0\le Q<X/Y$}\

Next take the intermediate range, where $Q_0\le Q<X/Y$.
We need to estimate \eqref{eq:cIQdef}.
Inserting \eqref{eq:Poisson} and opening the square gives
\beann
\cI_Q 
&=&
\int_{\gt = r/q + \gb \atop Q\le q<2Q, (r,q)=1,|\gb|<1/(QM)}
\left|
X^2
\sum_{\g\in \sF_T}
\sum_{k,\ell\in\Z}
\cS_q(r,k,\ell;\g)
\cJ_X(\gb, k ,\ell,q;\g)
\right|^2
d\gt
\\
&=&
X^4
\sum_{\g,\g'\in \sF_T}
\sum_{k,\ell,k',\ell'\in\Z}
\sum_{q\asymp Q}
\left[
\sum_{r(q)}'
\cS_q(r,k,\ell;\g)
\overline{\cS_q(r,k',\ell';\g')}
\right]
\\
&&
\hskip.5in
\int_{|\gb|<1/(QM)}
\cJ_X(\gb, k ,\ell,q;\g)
\overline{\cJ_X(\gb, k' ,\ell',q;\g')}
d\gb
.
\eeann
We apply  \lemref{lem:cJX} to get
\be\label{eq:cIQbnd1}
\cI_Q 
\ \ll\ 
{X^4\over N}
\sum_{\g,\g'\in \sF_T}
\sum_{k,\ell,k',\ell'\ll 1}
\sum_{q\asymp Q}
\left|
\sum_{r(q)}'
\cS_q(r,k,\ell;\g)
\overline{\cS_q(r,k',\ell';\g')}
\right|
.
\ee
%
%
Now we introduce a parameter 
\be\label{eq:Hdef}
H:=Q_0^{\eta_0/4}
,
\ee
where $\eta_0$ is the constant in  \eqref{eq:sFTqBnd},
 and decompose 
$$
\cI_Q  \ \le  \ \cI_Q^{(<)} +\cI_Q^{(\ge)} ,
$$ 
according to whether $\gcd(c,c')<H$ or $\gcd(c,c')\ge H$.
We first deal with the large gcd.

\begin{prop}
There exists some $\eta>0$ so that:
$$
\cI_Q^{(\ge)} \ll 
{|\widehat{\cR_N}(0)|^2\over N}
N^{-\eta}
.
$$
\end{prop}
\pf
Let $h\ge H$ be the gcd of $c$ and $c'$.
Then applying 
\corref{cor:gcdccp} (and notation therein) to \eqref{eq:cIQbnd1},
 we have
\beann
\cI_Q ^{(\ge)}
& \ll&
{X^4\over N}
\sum_{\g\in \sF_T}
\sum_{h\mid c_\g\atop h\ge H}
\sum_{q_1\mid BcP^2\atop q_1\ll Q}
\sum_{\g'\in \sF_T\atop c'\equiv 0(h)}
\sum_{q_1'\mid Bc'P^2\atop q_1'\ll Q}
\sum_{q\asymp Q\atop q\equiv 0((q_1,q_1'))}
{(q_1,q_1')\over q},
\eeann
where in the last sum, we weakened the condition that $q$ is divisible by both $q_1$ and $q_1'$ to just being divisible by their gcd.
Now estimating divisor sums and  applying  Nullstellensatz  \eqref{eq:sFTqBnd} in the $\g'$ summation gives
\beann
\cI_Q ^{(\ge)}
& \ll_\vep &
N^\vep {X^4\over N}
H^{-\eta_0}
|\sF_T|^2,
\eeann
from which the claim follows.
\epf

Next we handle the small gcd.
\begin{prop}
There exists some $\eta>0$ so that:
$$
\cI_Q^{(<)} \ll 
{|\widehat{\cR_N}(0)|^2\over N}
N^{-\eta}
.
$$
\end{prop}
\pf
We begin with the observation that $(c,c')<H$ implies that $c\neq c'$, since  $c\asymp T$ and $H=o(T)$.
We apply
\lemref{lem:cSq3},
 estimate
$
\gcd\left(q(q_1,q_1'),J\right)
$
by
$
\gcd\left((q_2,q_2'),J\right)
{q(q_1,q_1')\over (q_2,q_2')},
$
and apply \eqref{eq:Xmodq2q2p}, giving
\bea\nonumber
\cI_Q ^{(<)}
&\ll_\vep&
{X^4\over N}
\sum_{\g,\g'\in \sF_T\atop (c, c')<H}
\sum_{q\asymp Q}
q^{-1+\vep}
{1\over (q_2,q_2')^{1/4}}
(q_1q_1')^{1/2}(q_1,q_1')^{1/2}
\\ \nonumber
&&
\hskip.1in
\left(\gcd\left((q_2,q_2'),
(c-c')
(
B^2
-\gD c c'
)
\right)
\right)^{1/4}
.
\eea
Since $q=q_1q_2=q_1'q_2'$, we have that $(q_2,q_2')=q/[q_1,q_1']$. We crudely estimate
\be
\cI_Q ^{(<)}
\ \ll_\vep\
N^\vep
{X^4\over N}
\sum_{q\asymp Q}
q^{-5/4}
\sum_{\g,\g'\in \sF_T\atop {c\neq c' \atop  {q_1=(BcP^2,q)\atop q_1'=(Bc'P^2,q)}}}
(q_1q_1')
\left(\gcd\left(q,
BP^2(c-c')
(
B^2
-\gD c c'
)
\right)
\right)^{1/4}
.
\ee

Let $h=(q_1,q_1')$, then $h\mid BP^2(c,c')$, and $h\ll \min(Q,H)$. 
Then we can write $q_1=hg$, $q_1'=hg'$, with $(g,g')=1$. Also note that $h\mid (q,BP^2(c-c')(
B^2
-\gD c c'
)
)$, 
so we can write 
$$
(q,BP^2(c-c')(
B^2
-\gD c c'
)
)
=h\tilde g.
$$
Then it follows that 
$$
(\tilde g, g)\ll (BP^2(c-c')/h,g) \cdot (BP^2(B^2-\gD cc'), g).
$$
The first gcd on the right hand side above is $1$, since a factor of both should have been included in $h$.
Since $g\mid BP^2 \gD c c'$, the second gcd is equal to $(BP^2B^2,g)\ll1$.
A similar argument shows that 
$$
(\tilde g, g')\ll1,
$$
and hence  we have the bound
$$
[hg,hg',h\tilde g]
\gg gg'\tilde g
$$
on their least common multiple. (Since $h$ may be small, it will not help us in this estimate.)

Similarly, we note that
$$
(\tilde g,c) \ll 
(BP^2(c-c'),c) \cdot (B^2-\gD cc',c).
$$
The second gcd is again bounded, while the first is bounded by 
$
H,
$
by the definition of $\cI_Q^{(<)}$.

Thus we can estimate:
\beann
\cI_Q ^{(<)}
& \ll_\vep&
N^\vep
{X^4\over N}
Q^{-5/4}
\sum_{\g\in \sF_T}
\sum_{\g'\in \sF_T\atop {c\neq c', (c,c')< H}}
\sum_{h\mid BP^2(c,c')\atop h\ll Q}
\sum_{g\mid BP^2 c\atop g\ll Q}
\sum_{g'\mid BP^2 c'\atop g'\ll Q}
\\
&&\times
\sum_{\tilde g\mid BP^2(c-c')(B^2-\gD cc')\atop gg'\tilde g\ll Q, \ (\tilde g,c)\ll H}
(hghg')
\left(
h\tilde g
\right)^{1/4}
\sum_{q\asymp Q\atop q\equiv0([hg,hg',h\tilde g])}
1
\\
& \ll_\vep&
N^\vep
{X^4\over N}
Q^{-5/4}
H^{9/4}
\sum_{\g\in \sF_T}
\sum_{
\tilde g\ll Q \atop  (\tilde g,c)\ll H}
{Q\over \tilde g^{3/4}}
\sum_{\g'\in \sF_T\atop 
BP^2(c-c')(B^2-\gD cc')\equiv 0(\tilde g)}
1
.
\eeann
For fixed $c$, let $\cR=\cR(\tilde g)\subset \Z/\tilde g\Z$ denote the set of roots mod $\tilde g$ of the polynomial $BP^2(c-x)(B^2-\gD cx)$.
If $\gD=0$, then this is a linear polynomial with leading coefficient $BP^2$, so $\#\cR\ll1$. If $\gD\neq0$, then this is a reducible quadratic polynomial with leading coefficient $BP^2\gD c$; since $(\tilde g,c)\ll H$, it follows that $\#\cR\ll_\vep H^{1+\vep}$.
Finally, we have that 
\beann
\cI_Q ^{(<)}
& \ll_\vep&
N^\vep
{X^4\over N}
Q^{-5/4}
H^{9/4}
\sum_{\g\in \sF_T}
\sum_{
\tilde g\ll Q \atop  (\tilde g,c)\ll H}
{Q\over \tilde g^{3/4}}
\sum_{\ga\in\cR(\tilde g)}
\sum_{\g'\in \sF_T\atop 
c'\equiv \ga(\tilde g)}
1
.
\eeann
We apply  Nullstellensatz  \eqref{eq:sFTqBnd}  in the last summation to obtain:
\beann
\cI_Q ^{(<)}
& \ll_\vep&
N^\vep
{X^4\over N}
Q^{-5/4}
H^{9/4}
\sum_{\g\in \sF_T}
\sum_{
\tilde g\ll Q } 
{Q\over \tilde g^{3/4}}
H
{1\over \tilde g^{\eta_0}}
|\sF_T| 
\\
& \ll_\vep&
N^\vep
{X^4\over N}
Q^{-1/4}
H^{13/4}
 Q^{1/4-\eta_0}
|\sF_T| ^2
.
\eeann
The choice of the parameter $H$ in \eqref{eq:Hdef} ensures that we have saved a power of $Q\ge Q_0$. The claim follows since $Q_0$ is a power of $N$ (by \eqref{eq:Q0K0def}).
\epf

These two propositions establish  \thmref{thm:cEmain} in the intermediate range of $Q$.

\newpage

\subsection{Minor Arcs III: Case $X/Y \le Q<M$}\

In this largest range, we return to the exact evaluation:
\beann
\cI_Q 
&=&
X^4
\sum_{\g,\g'\in \sF_T}
\sum_{k,\ell,k',\ell'\in\Z}
\sum_{q\asymp Q}
\left[
\sum_{r(q)}'
\cS_q(r,k,\ell;\g)
\overline{\cS_q(r,k',\ell';\g')}
\right]
\\
&&
\hskip.5in
\int_{|\gb|<1/(QM)}
\cJ_X(\gb, k ,\ell,q;\g)
\overline{\cJ_X(\gb, k' ,\ell',q;\g')}
d\gb
.
\eeann

Now we break
$$
\cI_Q\le \cI_Q^= + \cI_Q^{\neq}
$$
depending on whether $c=c'$ or not. 
We first handle the latter case.

\begin{prop}
There is an $\eta>0$ so that
\beann
\cI_Q^{\neq}
&\ll&
{|\widehat\cR_N(0)|^2\over N}
N^{-\eta}
,
\eeann
as $N\to\infty$.
\end{prop}
\pf
To begin, we can use the last part of \lemref{lem:cJX} to show that:
\be\label{eq:cIQneq}
\cI_Q^{\neq}
\ll
{X^4\over QM}
\sum_{\g,\g'\in \sF_T}
\sum_{|k|,|\ell|,|k'|,|\ell'|\ll {Y q\over X}} \
\sum_{q\asymp Q}
\left|
\sum_{r(q)}'
\cS_q(r,k,\ell;\g)
\overline{\cS_q(r,k',\ell';\g')}
\right|
\ + \
O(N^{-100})
.
\ee
We will omit this last term henceforth.
In the case $c\neq c'$,  we estimate using
 \lemref{lem:cSq3},  crudely (e.g., $q_1\ll T$, etc) giving:
\beann
\cI_Q^{\neq}
&\ll&
{X^4\over QM}
\sum_{q\asymp Q}
\sum_{\g,\g'\in \sF_T\atop q_1=(BcP^2,q),q_1'=(Bc'P^2,q)}
\\
&&\times
\sum_{|k|,|\ell|,|k'|,|\ell'|\ll {Y q\over X}} \
q^{-5/4+\vep}
(q_1q_1')^{1/2}(q_1,q_1')^{1/4}
\gcd\left(q(q_1,q_1'),J, K\right)^{1/4}
\\
&\ll&
{N^\vep X^4T^{3/2}
\over Q^{9/4}M}
\sum_{q\asymp Q}
\sum_{\g,\g'\in \sF_T}
\sum_{|k|,|\ell|,|k'|,|\ell'|\ll {Y q\over X}} \
\gcd\left(q,J, K\right)^{1/4}
.
\eeann
Recall from \lemref{lem:cSq2} that $J$ does not depend on $k,k',\ell,\ell'$ but $K$ does.
Then
$$
\cI_Q^{\neq}
\ll
{Y^4 X^4T^{3/2}
\over Q^{9/4}M}
\left[{Q^4\over X^4}+1\right]
\sum_{q\asymp Q}
\sum_{\g,\g'\in \sF_T \atop c\neq c'}
\gcd\left(q,J\right)^{1/4}
.
$$
By \eqref{eq:Xmodq2q2p}, we have that:
$$
\gcd\left(q,J\right)
\ll
T^2
\gcd\left((q_2,q_2'),J\right)
\ll
T^2
\gcd\left((q_2,q_2'),(c-c')
(
B^2
-\gD c c'
)
\right)
\ll
T^5,
$$
since $B^2
-\gD c c'
$
is never zero.
Then
\beann
\cI_Q^{\neq}
&\ll&
{N^\vep X^4T^{3/2}
\over Q^{9/4}M}
\left[{Q^4\over X^4}+1\right]
Q
|\sF_T|^2
T^{5/4}
\ll
Y^4
\left[
{T^{3/2}}M^{7/4}+
{X^4T^{3/2}
\over X^{5/4}M}
\right]
|\sF_T|^2
T^{5/4}
\\
&\ll&
Y^4
|\sF_T|^2
X^{2}
{T^{5}\over X^{1/4}}
\ll
Y^4
{|\widehat {\cR_N(0)}|^2\over N}
{T^{6}\over X^{1/4}}
.
\eeann
The claim again follows due to the large power of $X$ savings (relative to the small loss of powers of $T$ and $Y$); see  \eqref{eq:TX2eqN} and  \eqref{eq:Ydef}.
\epf

Next we analyze the case that $c=c'$. At this stage, we decompose $\cI_Q^=$ further according to whether $k\ell=k'\ell'$ or not,
$$
\cI_Q^= \ll
\cI_Q^{=,=}
+
\cI_Q^{=,\neq}.
$$
We first analyze the case that $k\ell\neq k'\ell'$. 

 \begin{prop}
There is an $\eta>0$ so that
\beann
\cI_Q^{=,\neq}
&\ll&
{|\widehat\cR_N(0)|^2\over N}
N^{-\eta}
,
\eeann
as $N\to\infty$.
\end{prop}
\pf
We again apply \lemref{lem:cJX} as in \eqref{eq:cIQneq}.
In this case, we then apply \lemref{lem:cS4}, which gives:
\beann
\cI_Q^{=,\neq} 
&\ll&
{X^4\over QM}
\sum_{\g,\g'\in \sF_T\atop c=c'}
\sum_{|k|,|\ell|,|k'|,|\ell'|\ll {Y q\over X} \atop k\ell \neq k' \ell'} \
\sum_{q\asymp Q}
\left[
q^{-5/4+\vep}
q_1^{2}
\gcd\left(q,    \ell' k' - \ell k\right)^{1/4}
\right]
.
\eeann
The $\gcd$ is bounded by $|k\ell|\ll (Y q/X)^2$, giving:
\beann
\cI_Q^{=,\neq} 
&\ll&
{Y^{1/2}X^4\over QM}
\sum_{\g,\g'\in \sF_T\atop c=c'}
\sum_{|k|,|\ell|,|k'|,|\ell'|\ll {Y q\over X} \atop k\ell \neq k' \ell'} \
\sum_{q\asymp Q}
\left[
q^{-5/4+\vep}
q_1^{2}
(q/X)^{1/2}
\right]
\\
&\ll&
Y^5{X^4\over QM}
|\sF_T|^2
\left[ {Q^4\over X^4}+1 \right]
Q
Q^{-5/4}
T^{2}
Q^{1/2} X^{-1/2}
\\
&\ll&
Y^5
{|\widehat{\cR_N}(0)|^2\over N}
{T^{6}\over X^{1/4}}
.
\eeann
The claim again follows due to the large power of $X$ savings.
\epf

The last case is when $c=c'$ and $k\ell=k'\ell'$; here we will fight not for powers of $Q$ but powers of the much smaller parameter $T$.
We can save a factor of $T^{2\gd-1}$ from the fact that $c=c'$ (and hence there are only $T$ values for $\g'$, not $T^{2\gd}\asymp |\sF_T|$). But this is insufficient for a power gain in the end.
 So new ideas are needed.

 \begin{prop}\label{prop:eqeq}
There is an $\eta>0$ so that
\beann
\cI_Q^{=,=}
&\ll&
{|\widehat\cR_N(0)|^2\over N}
N^{-\eta}
,
\eeann
as $N\to\infty$.
\end{prop}

Before beginning the proof, we return to the original formulation:
\beann
\cI_Q^{=,=}
&=&
X^4
\sum_{\g,\g'\in \sF_T\atop c=c'}
\sum_{k,\ell,k',\ell'\in\Z \atop k\ell=k'\ell'}
\sum_{q\asymp Q}
\left[
\sum_{r(q)}'
\cS_q(r,k,\ell;\g)
\overline{\cS_q(r,k',\ell';\g')}
\right]
\\
&&
\hskip.5in
\int_{|\gb|<1/(QM)}
\cJ_X(\gb, k ,\ell,q;\g)
\overline{\cJ_X(\gb, k' ,\ell',q;\g')}
d\gb
.
\eeann
We apply the last part of \lemref{lem:cJX} to truncate the $k,\ell,k',\ell'$ range:
\beann
\cI_Q^{=,=}
&=&
X^4
\sum_{\g,\g'\in \sF_T\atop c=c'}
\sum_{|k|,|\ell|,|k'|,|\ell'| \ll {Y Q\over X} \atop k\ell=k'\ell'}
\sum_{q\asymp Q}
\left[
\sum_{r(q)}'
\cS_q(r,k,\ell;\g)
\overline{\cS_q(r,k',\ell';\g')}
\right]
\\
&&
\hskip.5in
\int_{|\gb|<1/(QM)}
\cJ_X(\gb, k ,\ell,q;\g)
\overline{\cJ_X(\gb, k' ,\ell',q;\g')}
d\gb
.
\eeann
Over this range of $k,\ell,k',\ell'$, we need to level out the $q$ dependence from the archimedean component $\cJ$. But the range of $q\asymp Q$ is too long for this purpose, so we decompose the sum into $U$ intervals, where $U$ is a parameter chosen to be
$$
U=Q^{1/2}.
$$
Each interval is of length $Q/U=Q^{1/2}$, which is much larger than $T$. Let $Q_1$ range over  the starting points of these intervals. Then we have:
\beann
\cI_Q^{=,=}
&=&
X^4
\sum_{\g,\g'\in \sF_T\atop c=c'}
\sum_{|k|,|\ell|,|k'|,|\ell'| \ll {Y Q\over X} \atop k\ell=k'\ell'}
\sum_{j\in\{1,\dots,U\}\atop Q_1=Q+j \frac QU}
\sum_{Q_1 \le q \le Q_1 + \frac Q U}
\left[
\sum_{r(q)}'
\cS_q(r,k,\ell;\g)
\overline{\cS_q(r,k',\ell';\g')}
\right]
\\
&&
\hskip.5in
\int_{|\gb|<1/(QM)}
\cJ_X(\gb, k ,\ell,q;\g)
\overline{\cJ_X(\gb, k' ,\ell',q;\g')}
d\gb
.
\eeann
On each of these sub-intervals, we will replace $q$ in $\cJ$ by $Q_1$, 
thereby freeing the $q$ variable for a purely modular analysis, as follows.
\begin{lem}\label{lem:levelJq}
For any $\vep>0$, we have that:
\beann
\cI_Q^{=,=}
&\ll_\vep&
\frac1{QM}
X^4
\sum_{\g,\g'\in \sF_T\atop c=c'}
\sum_{|k|,|\ell|,|k'|,|\ell'| \ll {Y Q\over X} \atop k\ell=k'\ell'}
\sum_{j\in\{1,\dots,U\}\atop Q_1=Q+j \frac QU}
\left|
\sum_{Q_1 \le q \le Q_1 + \frac Q U}
\sum_{r(q)}'
\cS_q(r,k,\ell;\g)
\overline{\cS_q(r,k',\ell';\g')}
\right|
\\
&&
\hskip.5in
\ + \
Y^4
{|\widehat\cR_N(0)|^2\over N}
\frac TU
.
\eeann
\end{lem}
\pf
Returning to the definition of $\cJ$, we see that 
\beann
&&
\hskip-.3in
\left|
\cJ_X(\gb, k ,\ell,q;\g)
-
\cJ_X(\gb, k ,\ell,Q_1;\g)
\right|
\ll
 \frac YU,
\eeann
since $k,\ell\ll Y Q/X$. We replace the two appearances of $q$ in $\cJ$ one at a time.

Each time, we apply \lemref{lem:cS6} to the resulting difference, which is bounded by
\beann
&\ll&
Y^2
\frac 1U
{X^4\over QM}
\sum_{\g,\g'\in \sF_T\atop c=c'}
\sum_{|k|,|\ell|,|k'|,|\ell'|\ll {Y q\over X} \atop k\ell = k' \ell'} \
\sum_{q\asymp Q\atop q_1=(c,q)}
\left[
{(c,q)^2\over q^2}
\sum_{r(q)}'
\bo_{\{
-\ell \equiv Pr
(B  a
+D  c)
(\mod q_1)
\atop
-k
\equiv
Pr
(A c
+B  d)
(\mod q_1)
\}}
\right]
\\
&\ll&
Y^2
\frac 1U
{X^4\over Q^3M}
\sum_{\g,\g'\in \sF_T\atop c=c'}
\sum_{ q_1\mid c\atop  q_1\ll Q}
q_1^2
\sum_{q\asymp Q\atop q\equiv 0 (\mod q_1)}
\sum_{|k|,|\ell|\ll {Y q\over X}} \
\sum_{r(q)}'
\bo_{\{
-\ell \equiv Pr
B  a
(\mod q_1)
\atop
-k
\equiv
Pr
B  d
(\mod q_1)
\}}
,
\eeann
where we used that $c\equiv0(q_1)$. Since $(a,c)=(c,d)=1$,
for given $\ell$, the value of $r$ is determined up to constants mod $q_1$, so there are $\ll q/q_1$ values of $r$ contributing. For each value of $(\ell,r)$, we have that $k$ is determined mod $q_1$, but we won't use this.
In total, we bound the difference by
\beann
&\ll&
Y^2
\frac 1U
{X^4\over Q^3M}
\sum_{\g,\g'\in \sF_T\atop c=c'}
\sum_{ q_1\mid c\atop  q_1\ll Q}
q_1^2
\sum_{q\asymp Q\atop q\equiv 0 (\mod q_1)}
{Q\over X}
{Q\over q_1}
{Q\over X}
\\
&\ll&
Y^2
\frac 1U
{X^4\over Q^3M}
\sum_{\g,\g'\in \sF_T\atop c=c'}
\sum_{ q_1\mid c\atop  q_1\ll Q}
q_1^2
{Q\over q_1}
{Q\over X}
{Q\over q_1}
{Q\over X}
\\
&\ll&
Y^2
{|\widehat\cR_N(0)|^2\over N}
\frac TU
.
\eeann

Finally, we estimate trivially that
$$
\int_{|\gb|<1/(QM)}
\cJ_X(\gb, k ,\ell,Q_1;\g)
\overline{\cJ_X(\gb, k' ,\ell',Q_1;\g')}
d\gb
\
\ll
\
\frac1{QM},
$$
whence the claim follows, since $Q\ge X/Y$.
\epf

Now we have leveled out the sum, and are in position to apply the crucial \lemref{lem:cS7}.
\pf[Proof of \propref{prop:eqeq}]
Inserting \lemref{lem:cS7} into \lemref{lem:levelJq} gives:
\bea
\nonumber
\cI_Q^{=,=}
&\ll_\vep&
\frac1{QM}
X^4
\sum_{\g,\g'\in \sF_T\atop c=c'}
\sum_{|k|,|\ell|,|k'|,|\ell'| \ll {Y Q\over X} \atop k\ell=k'\ell'}
U
\\
\nonumber
&&
\times
Q^\vep
\sum_{q_1\mid BcP^2\atop E=BcP^2/q_1}
\bo_{\{d\ell\equiv ak\equiv d'\ell' \equiv a'k'(\mod (q_1,c))\}}
\cN_{q_1} 
\sum_{Q_1\mid q_1\atop {p\mid Q_1 \Longrightarrow (p^\infty,q_1)\mid Q_1 \atop  (E,q_1/Q_1)=1}}
\\
\label{eq:072001}
&&\times
\left[
{
(EQ_1,Z)
\over UEQ_1}
+
{c
\over U^2}
+
\frac {c^3}Q
+
{|Z|\over UQ}
\right]
\\
\nonumber
&&
\hskip.5in
\ + \
N^\vep
{|\widehat\cR_N(0)|^2\over N}
\frac TU
,
\eea
where $Z$ and $\cN_{q_1}$ are as defined in the statement of \lemref{lem:cS7}.

We first handle the contribution from the latter three terms in \eqref{eq:072001}
\beann
&&
\frac{Q^\vep}{QM}
X^4
\sum_{\g,\g'\in \sF_T\atop c=c'}
\sum_{|k|,|\ell|,|k'|,|\ell'| \ll {Y Q\over X} \atop k\ell=k'\ell'}
\left[
{c
\over U}
+
\frac {c^3U}Q
+
{|Z|\over Q}
\right]
\\
&&
\ll
\frac{N^\vep Y^2}{QM}
X^4
|\sF_T|^2
\left(\frac QX\right)^2
\left[
{T
\over U}
+
\frac {T^3U}Q
+
{T^2\over Q}
\right]
\\
&&
\ll
N^\vep Y^2
{|\cR_N(0)|^2\over N}
\left[
{T^2
\over U}
+
\frac {T^4U}X
+
{T^3\over X}
\right]
,
\eeann
where we bounded $Z\ll T Y Q/X\ll Y T^2$ from \eqref{eq:Zdef} and used \eqref{eq:cNmbnd1}.
With $U=Q^{1/2}$, this is sufficient savings if $T$ is small enough relative to $Q\ge X/Y$; these are all power savings.
Only the first term of \eqref{eq:072001} remains to be handled.
\beann
&&\frac{Q^\vep}{QM}
X^4
\sum_{\g,\g'\in \sF_T\atop c=c'}
\sum_{|k|,|\ell|,|k'|,|\ell'| \ll {Y Q\over X} \atop k\ell=k'\ell'}
\sum_{q_1\mid BcP^2\atop E=BcP^2/q_1}
\bo_{\{d\ell\equiv ak\equiv d'\ell' \equiv a'k'(\mod (q_1,c))\}}
\sum_{Q_1\mid q_1\atop {p\mid Q_1 \Longrightarrow (p^\infty,q_1)\mid Q_1 \atop  (E,q_1/Q_1)=1}}
{(EQ_1,Z)
\over EQ_1}
\\
&&
\ll
\frac{Q^\vep}{QM}
X^4
\sum_{\g\in \sF_T}
\sum_{|k|,|\ell|,|k'|,|\ell'| \ll {Y Q\over X} \atop {k\ell=k'\ell' }}
\sum_{q_1\mid BcP^2\atop E=BcP^2/q_1}
\bo_{\{d\ell\equiv ak\equiv d'\ell' \equiv a'k'(\mod (q_1,c))\}}
\\
&&
\hskip.5in
\sum_{Q_1\mid q_1\atop {p\mid Q_1 \Longrightarrow (p^\infty,q_1)\mid Q_1 \atop  (E,q_1/Q_1)=1}}
{1
\over EQ_1}
\sum_{Z_1 \mid EQ_1}
Z_1
\sum_{\g'\in \sF_T\atop c=c'}
\bo_{\{Z(\g')\equiv 0 (\mod Z_1)\}}
.
\eeann
Now we apply \eqref{eq:Z1bnd} to the last summation, expanding $\g'\in\sF_T$ to all of $\SL_2(\Z)$ (recalling that in $\sF_T$, all entries are $\asymp T$). This gives
\beann
&&
\frac{Q^\vep}{QM}
X^4
\sum_{\g\in \sF_T}
\sum_{|k|,|\ell|,|k'|,|\ell'| \ll {Y Q\over X} \atop {k\ell=k'\ell' }}
\sum_{q_1\mid BcP^2\atop E=BcP^2/q_1}
\bo_{\{d\ell\equiv ak\equiv d'\ell' \equiv a'k'(\mod (q_1,c))\}}
\\
&&
\hskip.5in
\sum_{Q_1\mid q_1\atop {p\mid Q_1 \Longrightarrow (p^\infty,q_1)\mid Q_1 \atop  (E,q_1/Q_1)=1}}
{1
\over EQ_1}
\sum_{Z_1 \mid EQ_1}
T
(Z_1,\ell d - k a)
\\
&&
\ll
\frac{Q^\vep}{QM}
TX^4
\sum_{\g\in \sF_T}
\sum_{q_1\mid BcP^2\atop E=BcP^2/q_1}
\sum_{Q_1\mid q_1\atop {p\mid Q_1 \Longrightarrow (p^\infty,q_1)\mid Q_1 \atop  (E,q_1/Q_1)=1}}
{1
\over EQ_1}
\sum_{Z_1 \mid EQ_1}
\sum_{Z_2 \mid Z_1}
Z_2
\\
&&
\hskip.5in
\sum_{|k|,|\ell| \ll {Y Q\over X} \atop {d\ell\equiv ak (\mod (q_1,c)) \atop d\ell \equiv ak (\mod Z_2)}}
\sum_{|k'|,|\ell'| \ll {Y Q\over X} \atop {k\ell=k'\ell' }}
1
.
\eeann
The $k',\ell'$ sum is a divisor sum. Note that $(q_1,c)\asymp q_1$. Replace the condition $d\ell\equiv ak\mod (q_1,c)$ by $d\ell\equiv ak\mod (q_1/Q_1,c)$. Note that $q_1/Q_1$ is coprime to $EQ_1$. So with $k$ fixed, $\ell$ is restricted to a residue class mod $(q_1/Q_1,c)Z_2$.  This gives
\beann
&&
\frac{Y^4}{QM}
TX^4
\sum_{\g\in \sF_T}
\sum_{q_1\mid BcP^2\atop E=BcP^2/q_1}
\sum_{Q_1\mid q_1\atop {p\mid Q_1 \Longrightarrow (p^\infty,q_1)\mid Q_1 \atop  (E,q_1/Q_1)=1}}
{1
\over EQ_1}
\sum_{Z_2 \mid EQ_1}
Z_2
\\
&&
\hskip.5in
{Q\over X}
\left(
{Q\over X (q_1/Q_1,c)Z_2}
+1
\right)
\\
&&
\ll
{Y^2}
X^2
|\sF_T|
\ 
\ll
\
{Y^2}
{|\widehat{\cR_N}(0)|^2\over N}
{1\over T^{2\gd-1}}
.
\eeann
This gives the claim, by the choice of $Y$ in \eqref{eq:Ydef}.
\epf

\thmref{thm:cEmain} has now been established in all ranges of $Q$, thus completing the proof \thmref{thm:main1}.

\newpage

\bibliographystyle{alpha}

\bibliography{AKbibliog}

\end{document}